\documentclass[final,3p]{elsarticle}

\usepackage{amssymb}

\usepackage{graphics}

\usepackage{amsmath}
\usepackage{mathtools}
\usepackage{graphicx}
\usepackage{booktabs}
\usepackage{subcaption}
\usepackage{enumitem}
\usepackage{framed}
\usepackage{floatrow}
\usepackage{multirow}
\usepackage{setspace}
\usepackage{color}

\makeatletter
\def\ps@pprintTitle{%
 \let\@oddhead\@empty
 \let\@evenhead\@empty
 \def\@oddfoot{}%
 \let\@evenfoot\@oddfoot}
\makeatother

\begin{document}

\begin{frontmatter}



\title{Simulation of non-linear structural { elastodynamic} and impact problems using minimum energy
and simultaneous diagonalization high-order bases}

\author[1]{Dias, A. P. C.}
\author[2,3]{Suzuki, J. L.}
\author[1]{Valente, G. L.}
\author[1]{Bittencourt, M. L.\corref{cor1}}

\cortext[cor1]{Corresponding author: mlb@fem.unicamp.br}

\address[1]{Department of Integrated System, School of Mechanical Engineering, University of Campinas, Rua Mendeleyev 200, Campinas, SP, Brazil, Zip Code 13083-860}
\address[2]{Department of Mechanical Engineering, Michigan State University, 428 S Shaw Ln, East Lansing, MI 48824, USA}
\address[3]{Department of Computational Mathematics, Science and Engineering, Michigan State University, 428 S Shaw Ln, East Lansing, MI 48824, USA}


\begin{abstract}
We present the application of simultaneous diagonalization and minimum energy (SDME) high-order
finite element { modal} bases for simulation of transient non-linear
{ elastodynamic problem, including impact cases with neohookean hyperelastic materials.}
The bases are constructed using procedures for simultaneous diagonalization of the internal
modes and Schur complement of the boundary modes from the standard nodal and modal bases,
constructed using Lagrange and Jacobi polynomials, respectively. The implementation of these
bases in a high-order finite element code is straightforward, since the procedure is applied only
to the one-dimensional expansion bases. Non-linear transient structural problems with large deformation,
hyperelastic materials and impact are solved using the obtained bases with explicit and implicit time
integration procedures. { Iterative solutions based on preconditioned conjugate gradient methods
are considered.} The performance of the proposed bases in terms of the number of
iterations of pre-conditioned conjugate gradient methods and computational time are compared
with the standard nodal and modal bases. Our numerical tests obtained speedups up to 41 using the
considered bases when compared to the standard ones.
\end{abstract}


\end{frontmatter}


\section{Introduction}

The high-order finite element method (HOFEM) corresponds to the $p$-version of the standard
finite elements and the convergence of the approximate solution is achieved by increasing the
polynomial order of the basis functions \cite{GE_Karniadakis2005,ML_Bittencourt2014}.

The construction of appropriate basis functions is critical for the HOFEM due to the larger
condition numbers of the element matrices and increasing number of non-zero coefficients
as the polynomial order increases. The use of tensor-product bases is also very important
to reduce the computational complexity and therefore to improve the performance to calculate
the element operators and save memory space. Expansion bases for structured and non-structured
high-order elements are presented in
\cite{I_Babuska1981,P_Carnevali1993,S_Adjerid2001,BEUCHLER2001,GE_Karniadakis2005,Vej2007,ML_Bittencourt2007,Shen2007,WEBB2008}.
Orthogonality properties of the polynomial bases were utilized in most of these works
to obtain local matrices with better conditioning and sparsity.

In \cite{P_Carnevali1993}, basis functions and solution procedures for the $p$-version finite element analysis were described for linear elastostatics and dynamics. A hierarchical triangular element was developed
in \cite{WEBB2008} with the basis functions  constructed from orthogonal Jacobi polynomials which
made possible to achieve better conditioned matrices. A hierarchical basis for the $p$-version in two and three dimensions was presented in \cite{S_Adjerid2001}. The corresponding stiffness matrices had good sparsity properties and better conditioning than those generated from existing hierarchical bases.
In \cite{BEUCHLER2001}, the discretization of a 3D elliptic boundary value problem (BVP)  by means of the $hp$-version using a mesh of tetrahedra was investigated and several bases based on integrated Jacobi polynomials presented. Orthogonalization was used for a Legendre-Galerkin spectral method in \cite{Shen2007} to make the mass and stiffness matrices simultaneously diagonal.
A new class of higher-order finite elements based on generalized eigenfunctions of the Laplace operator was presented in \cite{Vej2007}. In \cite{ABDULRAHMAN2007115}, a set of hierarchical high-order basis functions for triangles was constructed using a systematic orthogonalization approach that yields better conditioning.
High-order bases have been also developed for the mixed finite element methods as in
\cite{CASTRO2016241,DEVLOO2019952}.

The conditioning and sparsity of the resultant system matrix after discretization directly influence the
numerical efficiency and performance of the solvers in numerical methods.
{ Consequently, the use of direct methods for the  solution of the systems of equations
become very costly due to the larger fill in of the factorization procedure.
The increasing of  condition numbers of element and global matrices for higher polynomial
orders could also limit the use of iterative methods as the number of iterations for convergence
within a given tolerance depends on the condition numbers of matrices. }
This aspect has stimulated the development of numerical preconditioners.
Many preconditioners for the $p$-version of the FEM based
on the topology of matrices and related to domain decomposition methods were exhaustively studied
theoretical and numerically
\cite{I_Babuska1989,I_Babuska1991,J_Mandel1990b,J_Mandel1990c,J_Mandel1994,MA_Casarin1997}.
The main idea used was to apply the Schur complement to the internal modes of the element matrices
and use the low-order (linear or quadratic) shape functions to construct the preconditioning matrices.
The condensation procedure computes the Schur complement in each element making interesting
the use of parallelization \cite{I_Babuska1989,J_Mandel1990c,MA_Casarin1997,VG_Korneev1997}.
Different block diagonal matrices may be constructed from the basic method generating different
versions of this preconditioner  \cite{J_Mandel1990a,S_Jensen1997}.

In this work, we apply the Schur complement on the basis level for the boundary modes using an
appropriate norm ($L_2$, energy or Helmholtz norms).

A hybrid preconditioning scheme employing
a nonorthogonal basis that combines global and locally accelerated preconditioners for rapid
iterative diagonalization of generalized eigenvalue problems in
electronic structure calculations was proposed in \cite{CAI201316}.
Numerical preconditioners have been also developed for other numerical methods in recent years,
as can be seen in
\cite{DEPRENTER2017297,DEPRENTER2019604,AGATHOS2019673,ABDULRAHMAN2007115,CF_Rodrigues2014}.
In \cite{DEPRENTER2017297} was established a scaling relation between the condition number
of the system matrices  and the smallest cell volume fraction for the Finite Cell Method. An algebraic
preconditioning technique was developed. Detailed numerical investigation about the effectiveness of the
preconditioner in improving the conditioning, convergence speed and accuracy of iterative solvers was
presented for the Poisson problem and for two- and three-dimensional problems of linear elasticity.
A dedicated Additive-Schwarz preconditioner that targets the underlying mechanism causing the ill-conditioning of immersed finite element methods  was presented in \cite{DEPRENTER2019604}. A detailed numerical investigation of the effectiveness of the preconditioner for a range of mesh sizes, isogeometric discretization orders, and {PDEs}, among which the Navier-Stokes equations, was presented. In \cite{AGATHOS2019673}, a combination of techniques to improve the convergence and conditioning properties of partition of unity enriched finite element methods was presented. The method was applied to discontinuous and singular enrichment functions keeping condition number growth rates similar to the standard finite elements.
Explicit analysis for structural and impact problems using
moderate high-order elements has been used in \cite{Bejar2016,Browning2020}. Critical time step
sizes for explicit time integration for quadratic bricks, thin plates, tetrahedra  and wedges  are
discussed in \cite{Bejar2016}. Application of quadratic elements for lumped-mass explicit
analysis of impact problems is presented in \cite{Browning2020}. In \cite{DANIELSON201663},
a twenty-one node wedge element is presented and used in the transition interfaces of
hexahedral-dominant meshes for problems of non-linear solid mechanics. The Barlow's method is
applied in \cite{DANIELSON201884} to determine superconvergent points for higher-order finite elements and for transverse stresses.

The HOFEM  has been extensively used  in
structural mechanics \cite{AlbertoNogueira2007,DongYosibash2009,JCP2015}.
showing the advantages of applying high-order methods for such problems. In
\cite{AlbertoNogueira2007,DongYosibash2009}, the modal basis presented in \cite{GE_Karniadakis2005}
was applied to large deformation problems using meshes of hexahedra and tetrahedra.
Mesh locking due to geometric properties and material  incompressibility are bypassed with
the HOFEM only by increasing the polynomial order above four of the mesh elements as presented
in \cite{ML_Bittencourt2014,yu,Suzuki2016_Incom}. The use of the HOFEM for the analysis of
a phase field model for fracture, damage and fatigue is discussed in
\cite{boldrini2016non,chiarelli2017comparison,haveroth2017lajss} and Mortar contact
finite elements are presented in \cite{DIAS201519,Dias2019}. There have also been applications on capturing the instability waves arising in near-wall flow interactions \cite{akhavan2020anomalous}. The advantages of
high-order approximations for non-linear structural problems, in terms of computational
costs and quality of the solutions, are clearly presented even when the standard nodal basis is
used. In this work, the proposed bases will obtain large speedups for the solution of similar
problems when compared with the standard modal bases.

Recently, simultaneous diagonalization has been used in high-order finite elements in
\cite{JCP2015,Dong2011,CAI201316,suzuki2017aspects}. This concept has also been widely employed in
the automatic control community \cite{LAUB1104549,moore1981}.
In \cite{Dong2011}, simultaneous diagonalization was used for the construction of the 1D interior
modes that made the element mass and stiffness submatrices of the interior modes simultaneously
diagonal. A Gram-Schmidt orthogonalization procedure was used to make orthogonal the vertex and
interior modes. In \cite{JCP2015}, PDE-specific high-order bases for squares and hexahedra,
based on simultaneous diagonalization for the internal modes and minimum energy
techniques for the vertex modes (SDME bases), were constructed and applied to linear
transient elastic problems.

{The condition numbers of the preconditioned matrices may be still large for
higher polynomials orders and the number of iterations for convergence increases quickly. The
element matrices calculated using the SDME bases have very low increasing of the condition numbers
with higher polynomial orders and consequently fewer iterations are required for the convergence of conjugate
gradient methods.}

Despite the extensive amount of works involving the application of the HOFEM for PDEs, to the
authors best knowledge, high-order SDME bases applied to impact problems have not been
found in the literature. In this paper, we apply the high-order SDME finite element bases for
transient structural problems {with geometrical, material, and boundary/interface nonlinearities,
here respectively considered through} large deformations, hyperelastic materials and impact
problems. We perform time-integration through explicit (central difference) and implicit
(Newmark) methods.
The bases are obtained from nodal and modal bases constructed with Lagrange and
Jacobi polynomials using procedures for simultaneous diagonalization of the internal
modes and Schur complement of the boundary modes. The performance of the proposed
bases are compared with the results obtained for standard nodal and modal bases.

\section{Construction of  high-order bases} \label{s2}

{
The HOFEM uses nodal and modal bases constructed from Lagrange and Jacobi polynomials,
respectively, to develop approximation solutions.

Consider a set of $P+1$ nodal or collocation points on the standard one-dimensional element in the
interval $-1\leqslant \xi_1\leqslant 1$, as illustrated in Figure \ref{fig-elem-nodalpoints}.
Lagrange polynomial of degree $P$ associated to an arbitrary node $a$, denoted as
$L_{a}^{P}(\xi_1)$, is given by
\begin{equation}
L_{a}^{(P)}(\xi_1) = \frac{ \prod_{b=0,a\neq b}^{P} \left (\xi_1 - \xi_{1_b} \right ) }
{ \prod_{b=0,a\neq b}^{P} \left ( \xi_{1_a} - \xi_{1_b} \right ) },
\label{eq:lagrangian}
\end{equation}
where $L_{a}^{(0)} \left (\xi_1 \right ) = 1$. The Lagrange polynomials have the
collocation property $L_{a}^{(P)}(\xi_{1_b})=\delta_{ab}$, where $\delta_{ab}$
is the Kronecker's delta. Gauss-Lobatto-Legendre collocation points are in general used
to avoid very oscillatory behavior of the Lagrange polynomials and improve
the conditioning of the finite element matrices to be calculated.

\begin{figure}[!thp]
  \center
  \includegraphics[height=1.6cm]{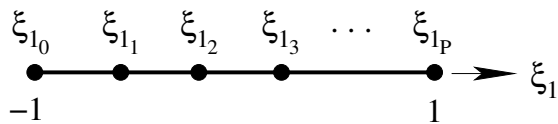}
  \caption{Nodal points on the standard coordinate system $\xi_1$ of the line element  \cite{ML_Bittencourt2014}.}
  \label{fig-elem-nodalpoints}
\end{figure}

The one-dimensional nodal standard basis (ST), denoted as $\psi_{p}(\xi_1)$,  is given by
Lagrange polynomials and indicated by
\begin{equation}
\psi_{p}(\xi_1)  =
\left\{\begin{array}{lll}
L_{0}^{(P)}(\xi_1), & & p = 0,\\
L_{P}^{(P)}(\xi_1), & & p = P,\\
L_{p}^{(P)}(\xi_1), & & 0 < p <  P,
\end{array}\right. .
\label{eq:shape_functions_lag}
\end{equation}
The shape functions are commonly associated to the element topological entities. In the case of
line element and nodal basis, the topological entities are the vertices and body, which
corresponds to the node indices $(p = 0$ and $p = P)$ and $(0 < p < P)$, respectively.
}

The one-dimensional modal standard basis of order $P$ is defined in the local
coordinate system $\xi_{1}$ as \cite{GE_Karniadakis2005,ML_Bittencourt2007}
\begin{equation}\label{BasePadrao1D}
\psi_{p}(\xi_{1})= \left\{ \begin{array}{lll}
\frac{1}{2}(1-\xi_{1}), & &p = 0,\\
\frac{1}{2}(1+\xi_{1}), & & p = 1,\\
\frac{1}{4}(1-\xi_{1})(1+\xi_{1})\mathcal{P}_{p-2}^{\alpha, \beta}
(\xi_{1}), & & 2 \leq p \leq P,
\end{array} . \right.
\end{equation}
where $\mathcal{P}_{p}^{\alpha,\beta}(\xi_{1})$ indicates the Jacobi
orthogonal polynomials of order $p$ and weights $(\alpha,\beta)$.
The vertex or boundary functions correspond to the indices $p = 0$ and $p = 1$;
$2\leq p\leq P$ for the internal functions.

{
In the HOFEM, nodal bases are used in general with collocation integration which results in diagonal
or spectral mass matrices. Modal bases are also used in general with consistent numerical integration
and does not result in spectral mass matrices. There are advantages and disadvantages of using
both approaches as stated in \cite{GE_Karniadakis2005,Canuto2006}. In this work, we consider
hierarchical modal bases.
}

\begin{figure}[!htbp]
        \centering
        \begin{subfigure}[b]{0.25\textwidth}
				\includegraphics[width=\columnwidth]{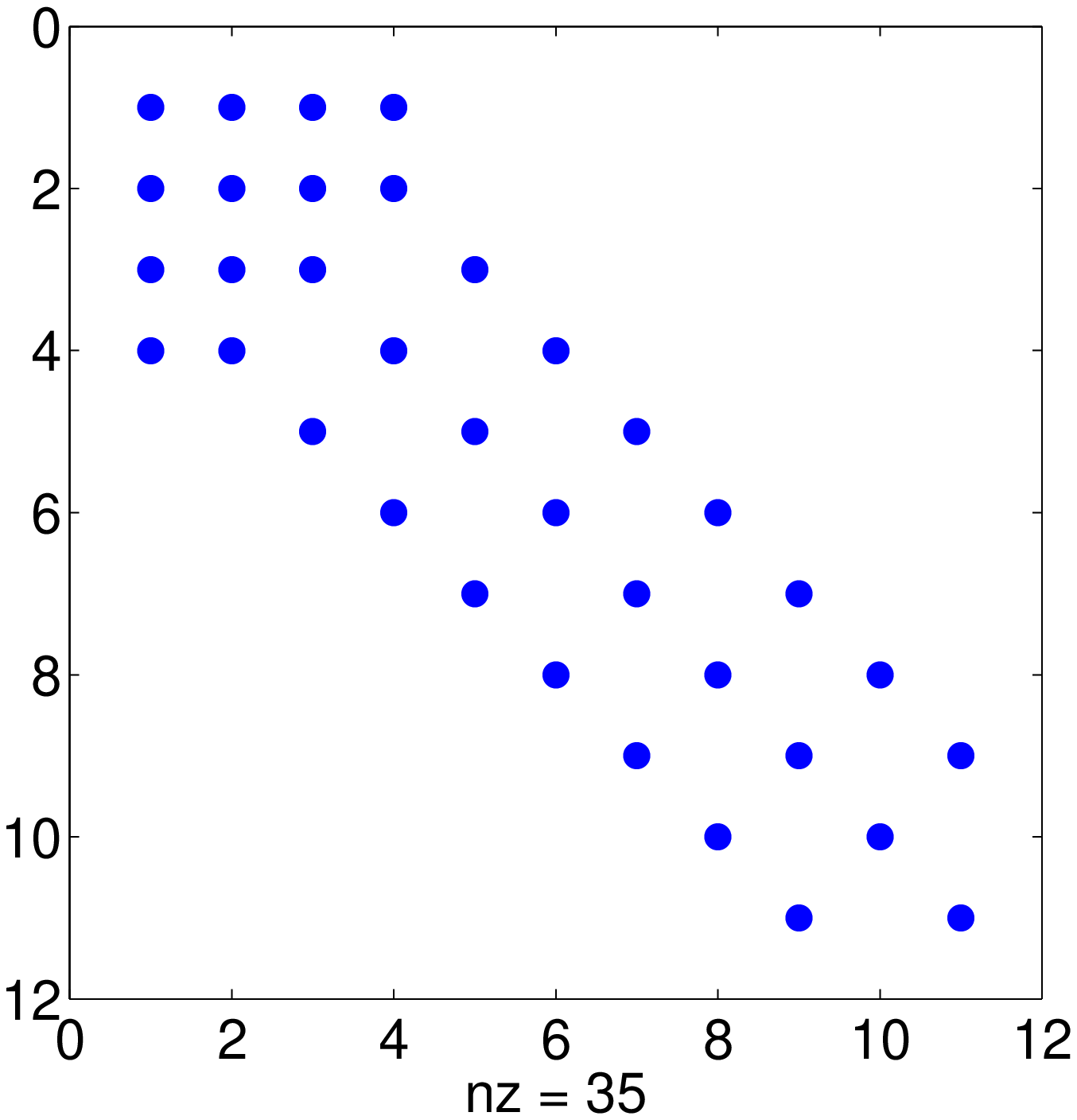}
				\caption{$\left[ M \right]$, ST basis.}
        \end{subfigure}%
        ~
        \begin{subfigure}[b]{0.25\textwidth}
                \includegraphics[width=\columnwidth]{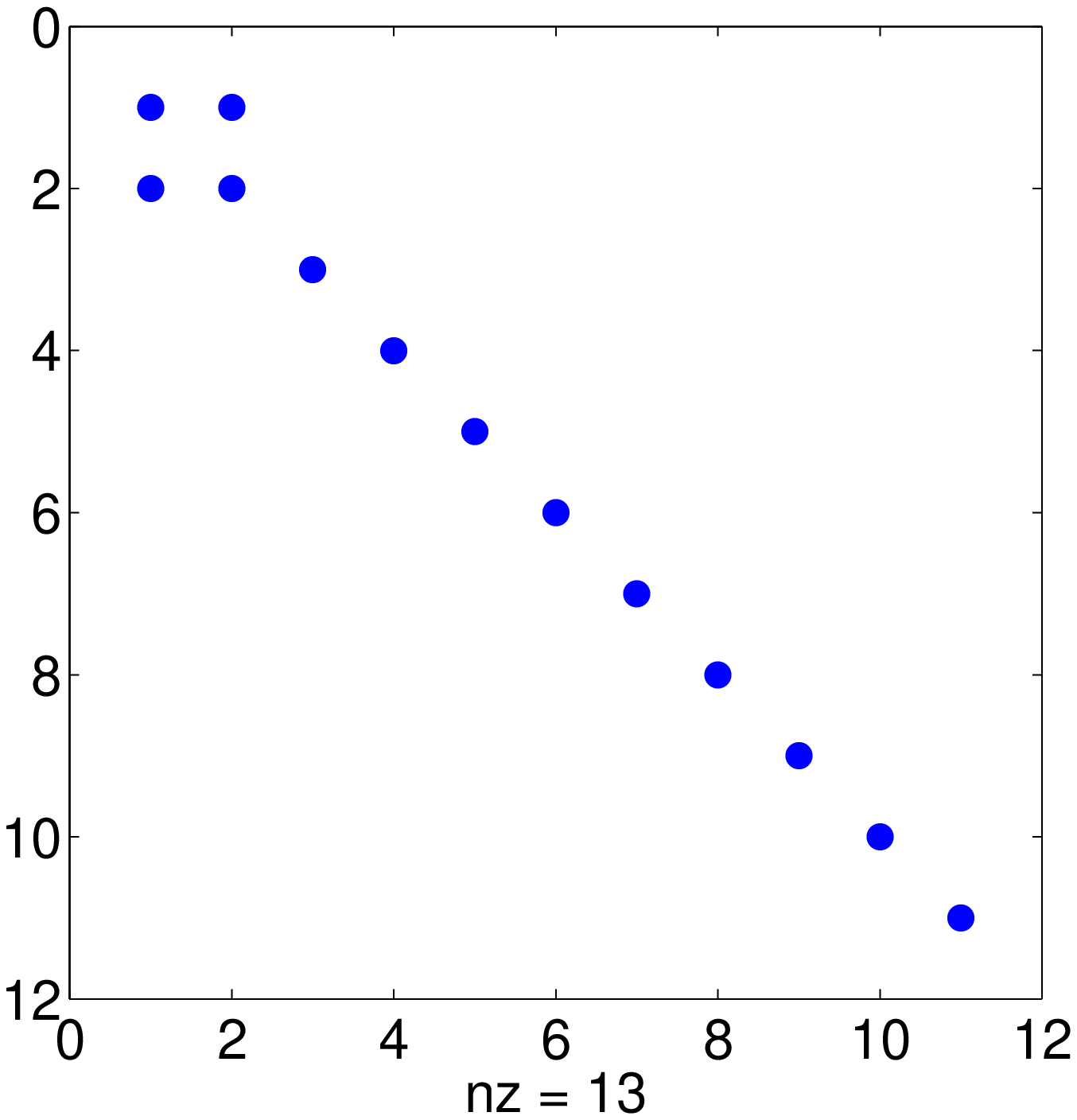}
                \caption{$\left[ K \right]$, ST basis.}
        \end{subfigure}
        ~
        \begin{subfigure}[b]{0.25\textwidth}
                \includegraphics[width=\columnwidth]{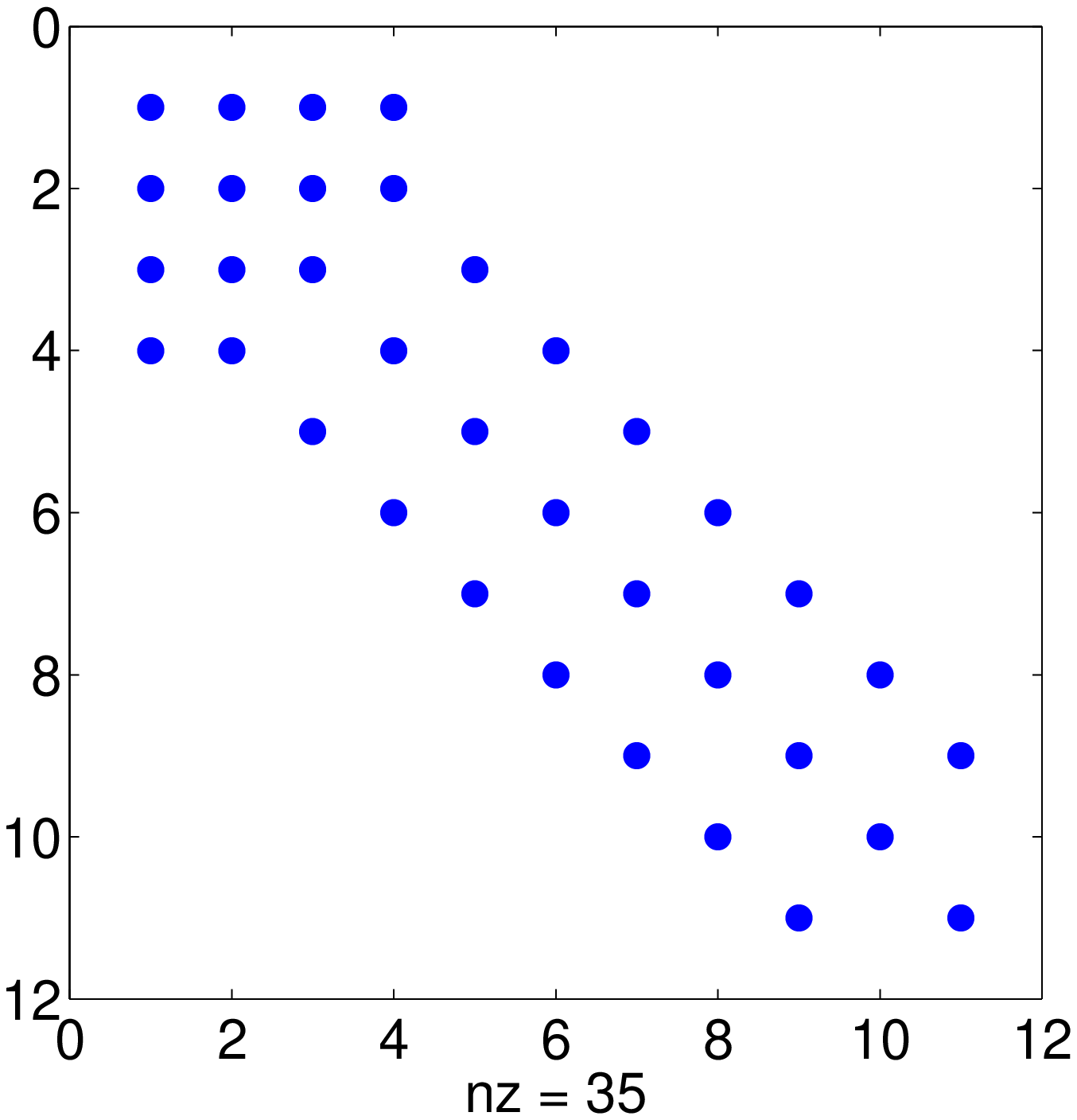}
                \caption{$\left[\hat{K}\right]$, ST basis.}
        \end{subfigure}
         \begin{subfigure}[b]{0.25\textwidth}
                \includegraphics[width=\columnwidth]{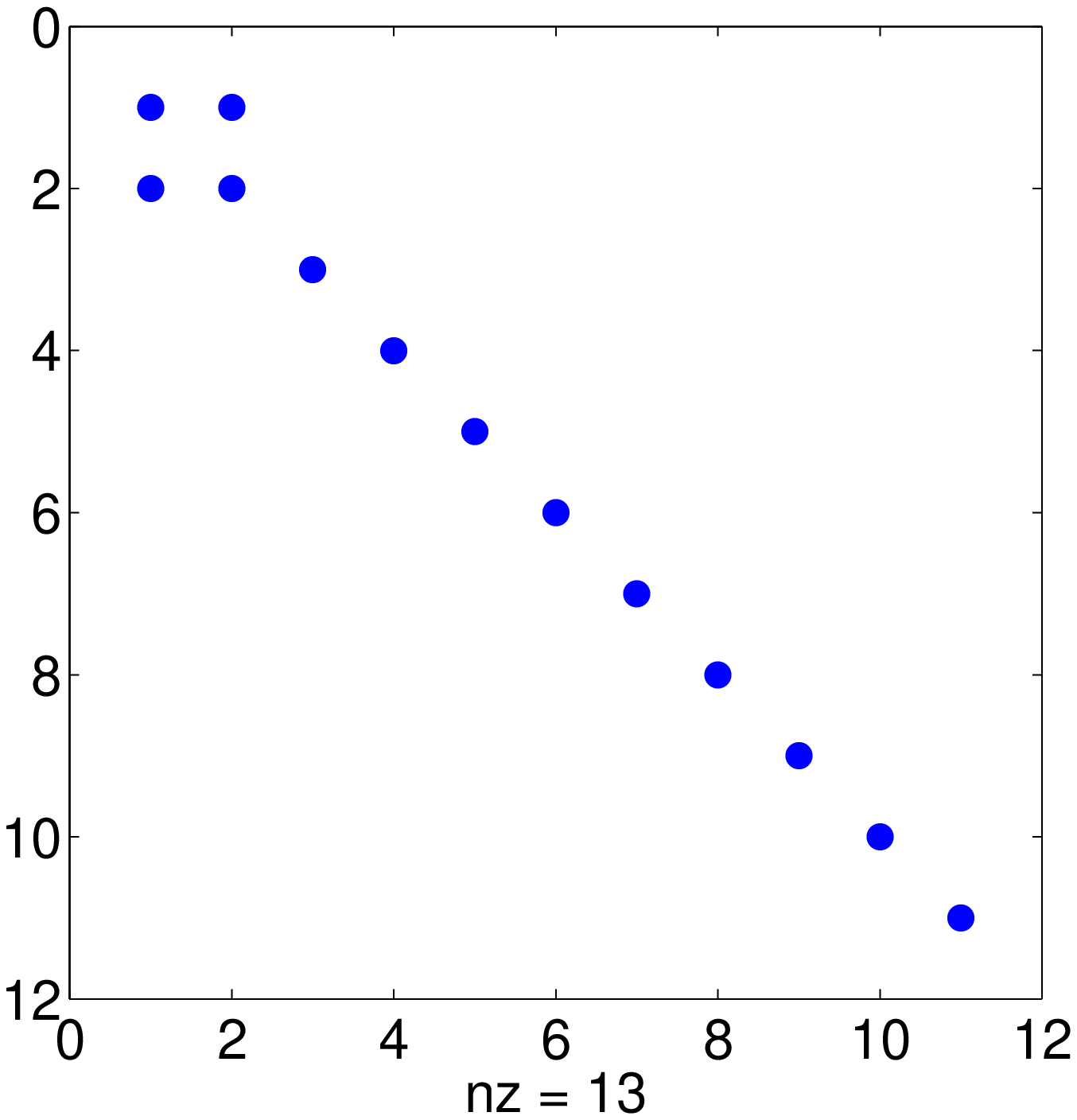}
                \caption{$\left[ M \right]$, SDME-M basis.}
        \end{subfigure}
        ~
        \begin{subfigure}[b]{0.25\textwidth}
                \includegraphics[width=\columnwidth]{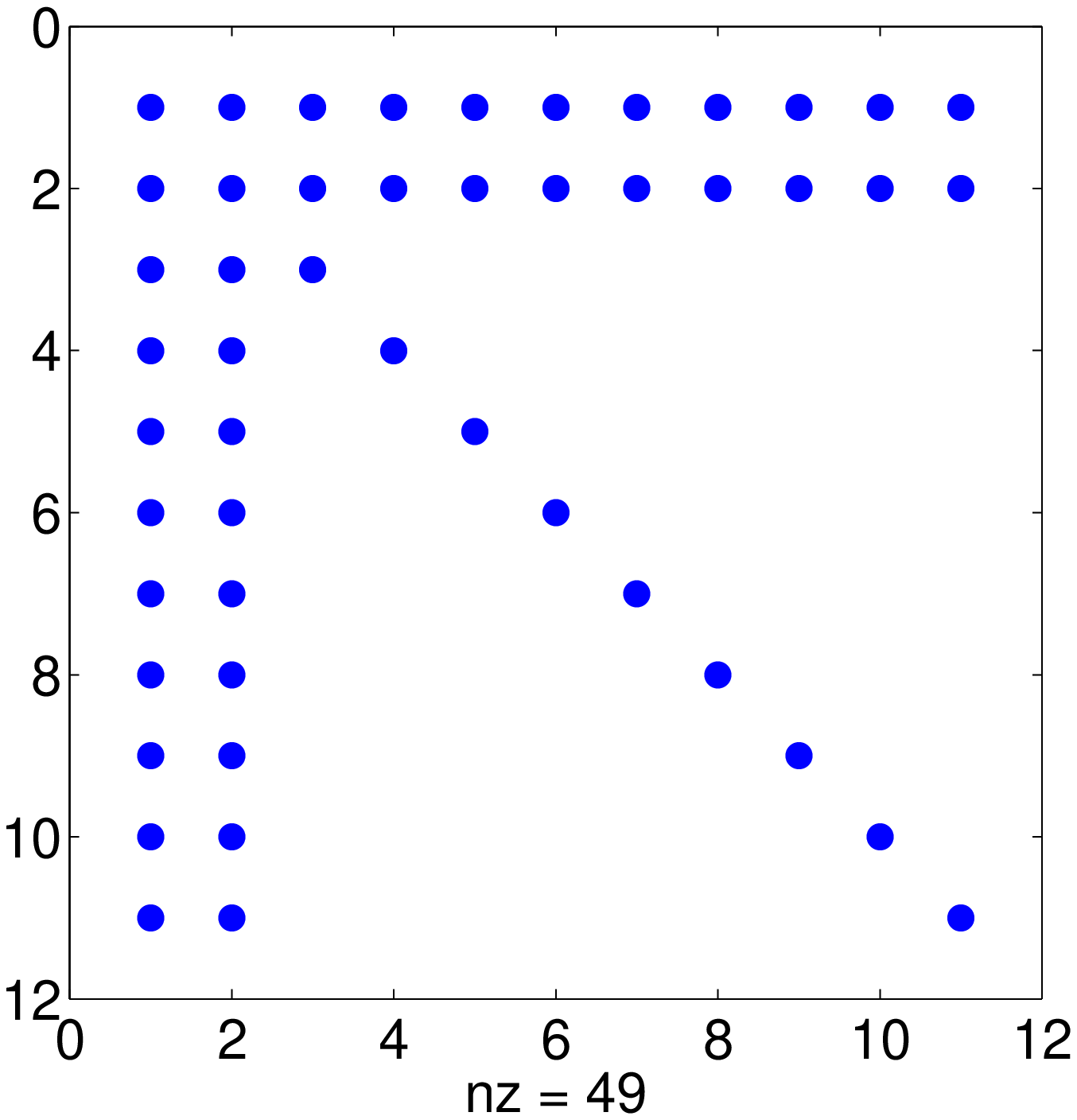}
                \caption{$\left[ K \right]$, SDME-M basis.}
        \end{subfigure}
        ~
        \begin{subfigure}[b]{0.25\textwidth}
                \includegraphics[width=\columnwidth]{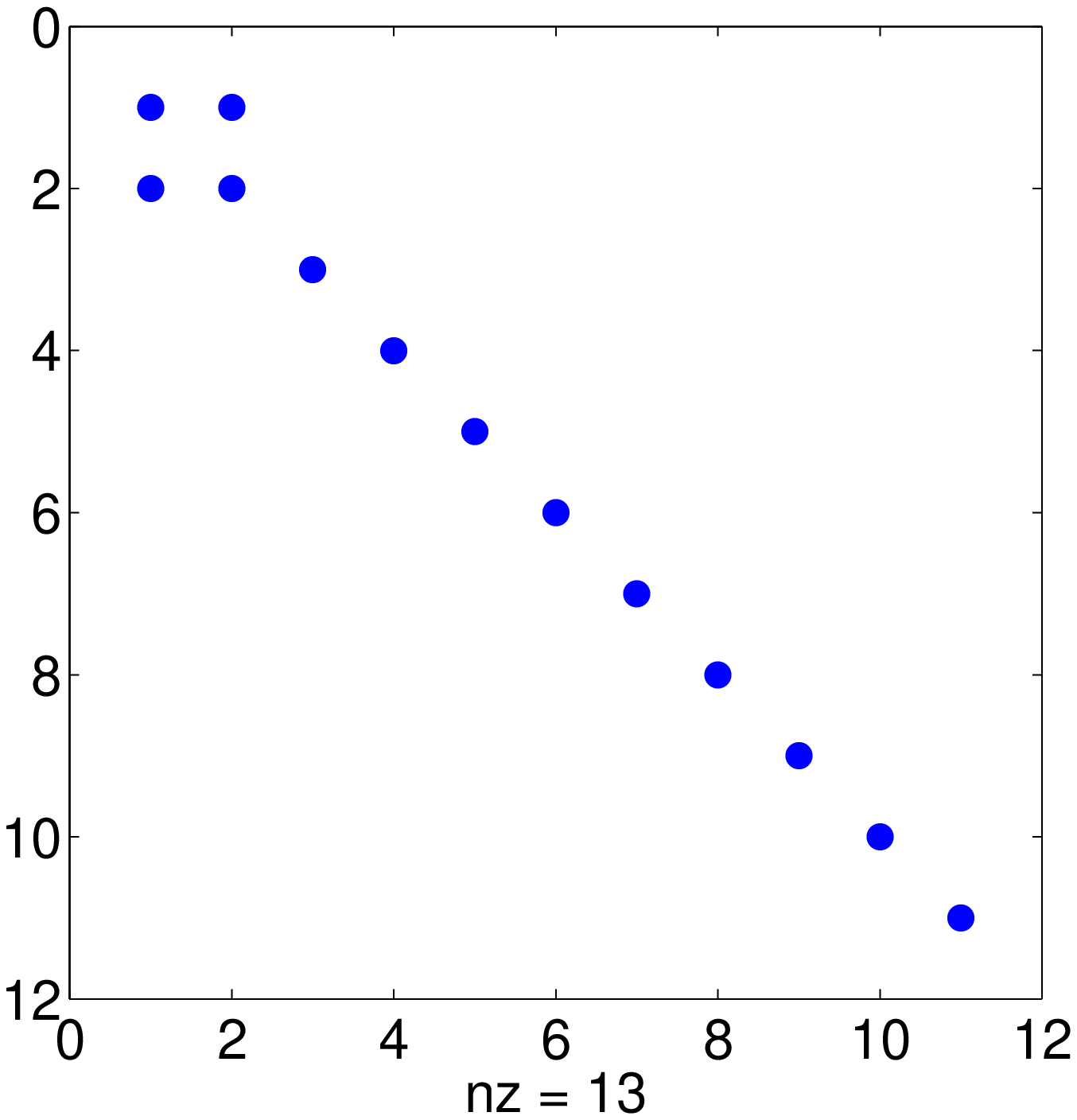}
                \caption{$\left[\hat{K}\right]$, SDME-H basis.}
        \end{subfigure}
\caption{Sparsity profiles of the mass, stiffness and effective
stiffness one-dimensional local matrices for $P = 10$ using the
standard (ST) and SDME bases with $\lambda = 1$.}
\label{RigidezMassaHelmholtz.1D.estruturado}
\end{figure}

The local coefficients of the one-dimensional mass and
stiffness matrices are respectively given by
\begin{equation}\label{Matrizesmassa}
M_{pq} = \int^{1}_{-1}\psi_{p}(\xi_{1}) \psi_{q}(\xi_{1}) d\xi_{1},
\end{equation}
\begin{equation}\label{Matrizesrigidez}
K_{pq} = \int^{1}_{-1} \psi_{p,_{\xi_{1}}} (\xi_{1})
\psi_{q,_{\xi_{1}}}(\xi_{1}) d\xi_{1},
\end{equation} where $0 \leq p,q \leq P$ and $\psi_{p,_{\xi_{1}}}$
is the derivative of $\psi_{p}$ with respect to $\xi_{1}$. Figs.
\ref{RigidezMassaHelmholtz.1D.estruturado}a and \ref{RigidezMassaHelmholtz.1D.estruturado}b
illustrate the sparsity profiles of the mass and stiffness matrices obtained with the modal
basis for $P=10$ and $\alpha=\beta=1$.

The previous mass and stiffness element matrices can be partitioned in terms
of the vertex and internal modes as
\begin{equation} \label{BlocosMassa}
[M]=\left[\begin{array}{cc}
 \left[ M_{vv} \right]  & \left[ M_{vi}\right]  \\
\left[ M_{vi}\right]^{T}  & \left[ M_{ii} \right]
\end{array}\right] \;\;\;\;\mbox{ and }\;\;\;\;
[K]=\left[\begin{array}{cc}
 \left[ K_{vv} \right]  & \left[ K_{vi}\right]  \\
\left[ K_{vi}\right]^{T}  & \left[ K_{ii} \right]
\end{array}\right].
\end{equation}

{
The bases presented here will modify the standard modal ones using two procedures.
Th first one is the simultaneous diagonalization (SD) of the internal blocks of the
element local mass and stiffness matrices. The second one is the minimum energy (ME)
which orthogonalizes the boundary and internal modes using one of the
following  norms:
\begin{itemize}
\item for the $L_2$ norm, the coupling block $\left[ M_{vi}\right] $ of the local mass
matrix is zeroed and the obtained basis is denominated SDME-M;
\item for the energy norm, the coupling block $\left[ K_{vi}\right] $ of the local
stiffness matrix is zeroed and the obtained basis is denominated SDME-K;
\item for a modified $H_1$ norm (see Equation (\ref{H1norm})), the coupling block
$\left[ \hat{K}_{vi}\right] $ of the local equivalent stiffness of the Newmark
method is zeroed and the obtained basis is denominated SDME-H.
\end{itemize}
In this work, we are interested in the SDME-M and SDME-H bases which will be used
to construct approximations for  explicit and implicit transient analyses.
}

The internal modes for the basis given in Eq.(\ref{BasePadrao1D})
will be transformed according to \cite{JCP2015,Shen2007,Vej2007} as
\begin{equation}\label{Diagonalizacao1}
\varphi_{p}(\xi_{1}) = \sum^{P}_{q=2} y_{pq}\psi_{q}(\xi_{1}).
\end{equation}
The coefficients $y_{pq}$ are entries of matrix $\left[ Y \right]$ such that the internal modes of the
new mass $\left[ M^{'}_{ii}\right] $ and stiffness $\left[ K^{'}_{ii}\right] $ matrices related to the
internal modes are given respectively by
\begin{equation}\label{MatrizInterna}
\left[ M^{'}_{ii}\right]  = \left[ Y \right] \left[ M_{ii}\right] \left[ Y\right]^{T} \;\; \mbox{and} \;\;
\left[ K^{'}_{ii}\right]  = \left[ Y\right] \left[ K_{ii} \right]  \left[ Y\right] ^{T}.
\end{equation}

The standard internal mass matrix $\left[ M_{ii}\right]$  can be made diagonal using the eigenvalue decomposition
\begin{equation}\label{diagMassa}
\left[ X\right]^{T} \left[ M_{ii}\right] \left[ X \right] = \left[ \Lambda_{M}\right],
\end{equation}
where $\left[ X\right] $ is the eigenvector matrix of $\left[ M_{ii}\right] $ and
$ \left[ \Lambda_{M}\right] $ is the diagonal matrix with the eigenvalues of $\left[ M_{ii}\right] $.
Based on that, we can define the matrix
\begin{equation}\label{L}
\left[ L \right] = \left( \left[ X\right] \left[ \Lambda_{M}^{-\frac{1}{2}}\right] \right) ^{T}
\left[ K_{ii}\right] \left( \left[ X\right] \left[ \Lambda_{M}^{-\frac{1}{2}}\right] \right),
\end{equation}
which is also symmetric and positive-definite and can be diagonalized  as
\begin{equation} \label{diagL}
\left[ Z\right]^{T}\left[ L\right] \left[ Z\right]  = \left[ \Lambda_{S}\right],
\end{equation}
in which $\left[ Z\right] $ denotes the matrix of the eigenvectors and $\left[ \Lambda_{S}\right] $
represents the diagonal matrix with eigenvalues of $\left[ L\right]$. Therefore, $\left[ Y\right] $ is then defined as
\begin{equation} \label{Y}
\left[ Y \right] = \left( \left[ X\right] \left[ \Lambda_{M}^{-\frac{1}{2}}\right] \left[ Z\right] \left[ \Lambda_{S}^{-\frac{k}{2}}\right] \right)^{T},
\end{equation}
where $k \in \left[ 0, 1 \right]$ is a parameter that influences the condition number of the matrices related to
the internal modes.

Substituting $\left[ Y\right]$ from (\ref{Y}) into (\ref{MatrizInterna}) yields
\begin{equation}\label{ParametroK}
\left[ M^{'}_{ii}\right] = \left[ Y\right] \left[ M_{ii}\right] \left[ Y \right]^{T} = \left[ \Lambda_{S}^{-k} \right]
\;\; \mbox{and} \;\;
\left[ K^{'}_{ii}\right] = \left[ Y\right] \left[ K_{ii}\right] \left[ Y \right]^{T} = \left[ \Lambda_{S}^{1-k} \right] .
\end{equation}
For $k = 0$, the internal block of the mass matrix is the identity
matrix and the condition number is $1$ for any polynomial order.
Analogously for the stiffness matrix with $k = 1$. The same condition number of the
internal mass and stiffness matrices is obtained for $k = \frac 12$.

The construction of minimum energy bases is equivalent to apply the Schur complement
for the vertex modes. The minimum energy extension of the standard basis
is computed as  \cite{Vej2007}:
\begin{equation}
\varphi_{k} = \psi_{k}^{v} - \sum_{j=2}^{P} \alpha_{kj} \psi_{j}^{i}, \;\; k = 0,1,
\label{ExtensaoMEmassa}
\end{equation}
where the coefficients $\alpha_{kj}$ are defined according to an appropriate norm. For instance, $\alpha^{M}_{kj}$
denotes the coefficients using the $L_2$ (or mass) norm and are uniquely determined as \cite{JCP2015}
\begin{equation} \label{DetermiALpha}
\langle \psi^{v}_{k},\psi^{i}_{l} \rangle_{L_2} -
\sum_{j=2}^{P}\alpha_{kj}^{M} \langle \psi^{i}_{j},\psi^{i}_{l}
\rangle_{L_2} = 0,  \;\; \forall \psi^{i}_{l} \in V^{i},
\end{equation}
which results in the following matrix for the coefficients $\alpha_{kj}^{M}$:
\begin{equation} \label{alphaMassa}
\left[ \alpha^{M} \right] = \left[ M_{vi}\right] \left[ M_{ii}\right] ^{-1}.
\end{equation}
We consider the simultaneous diagonalization (SD) of the internal
blocks of the mass and stiffness matrices to construct the
one-dimensional internal modes and the minimum energy (ME)
orthogonalization for the boundary modes based on the
choice of the appropriate norm according to the considered problem \cite{JCP2015}.
The obtained bases are labeled SDME. Specifically, when using $\left[ \alpha^{M} \right]$,
we denote the basis as SDME-M.
{Figures} \ref{RigidezMassaHelmholtz.1D.estruturado}d and \ref{RigidezMassaHelmholtz.1D.estruturado}e show
the sparsity patterns of the local one-dimensional mass and stiffness matrices for the standard basis
and the SDME-M basis.

We can also write Eq.(\ref{DetermiALpha}) in terms of the energy norm as {follows:}
\begin{equation} \label{DetermiALphaK}
\langle \psi^{v}_{k},\psi^{i}_{l} \rangle_{E} -
\sum_{j=2}^{P}\alpha_{kj}^{K} \langle \psi^{i}_{j},\psi^{i}_{l}
\rangle_{E}=0,  \;\; \forall \psi^{i}_{l} \in V^{i},
\end{equation}
which results in the following matrix for the coefficients $\alpha_{kj}^{K}$  \cite{JCP2015}:
\begin{equation} \label{alphaRigidez}
\left[ \alpha^{K}\right]  = \left[ K_{vi} \right]
\left[ K_{ii}\right] ^{-1}.
\end{equation}
We observe that matrices $ \left[ \alpha \right] $ influence the coupling blocks $\left[ M_{vi} \right]$ and
$\left[ K_{vi} \right]$ of the mass and stiffness matrices.  For $\left[ \alpha^{M}\right]$, the
basis does not decouple the internal and boundary modes of the one-dimensional stiffness matrix.
However, the one-dimensional mass matrix has the internal and boundary blocks uncoupled.

Particularly, when using the implicit Newmark scheme for time integration, an effective stiffness
matrix $[\hat{K}]$ of the following form arises:
\begin{equation} \label{EffectiveStiffness_P2}
[\hat{K}] =  \left[ K \right] + a_{0}\left[ M \right],
\end{equation}
where $a_0=\frac{1}{4 \Delta t^2}$ and $\Delta t$ {represents} the time increment. For the construction of the one-dimensional basis, we can associate the coefficient $a_0$ with the parameter $\lambda$, such that:
\begin{equation} \label{EffectiveStiffness_lambda}
[\hat{K}] =  \left[ K \right] + \lambda\left[ M \right],
\end{equation}

The matrix $[\hat{K}]$ can be expressed in terms
of vertex and internal modes. Considering the
minimum energy procedure for the energy norm of function
$u$ given by
\begin{equation}
\Vert u \Vert^{2}_{E} = \langle u', u'\rangle_{L_2} + \lambda \langle u, u \rangle_{L_2},
\label{H1norm}
\end{equation}
we obtain the following coefficients for the matrix $[\alpha^{\hat{K}}]$:
\begin{equation} \label{alphaKN}
\left[ \alpha^{\hat{K}}\right]  = \left[ \hat{K}_{vi} \right] \left[ \hat{K}_{ii}\right] ^{-1}.
\end{equation}
The boundary and internal blocks are uncoupled for $[{\hat{K}}]$
as illustrated in Fig.\ref{RigidezMassaHelmholtz.1D.estruturado}f.
The obtained basis is labeled as SDME-H.

\begin{figure}[!htbp]
        \centering
        \begin{subfigure}[b]{0.325\textwidth}
				\includegraphics[width=\columnwidth]{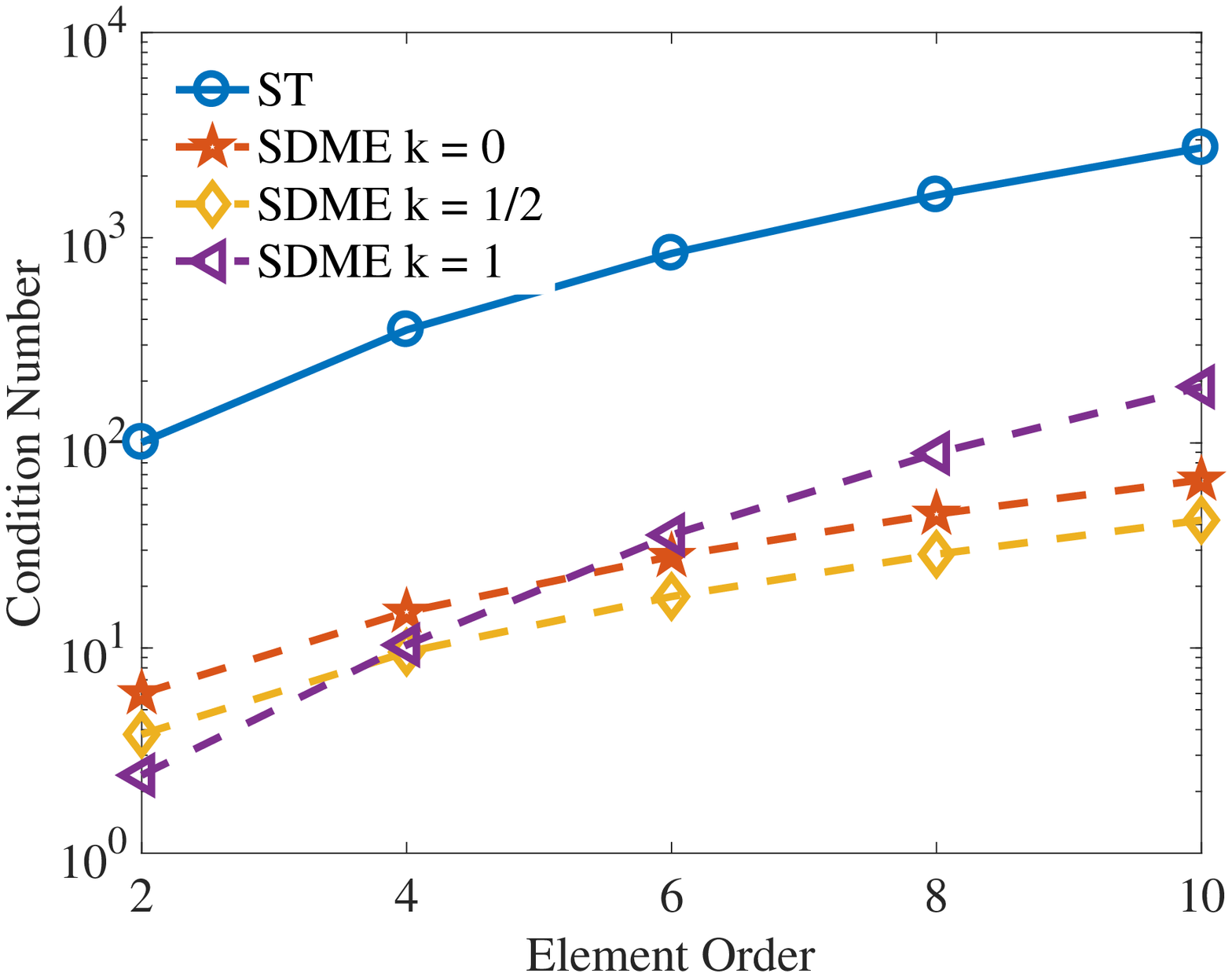}
				\caption{$\left[ M \right]$, ST and SDME-M bases.}
        \end{subfigure}%
        ~
        \begin{subfigure}[b]{0.325\textwidth}
                \includegraphics[width=\columnwidth]{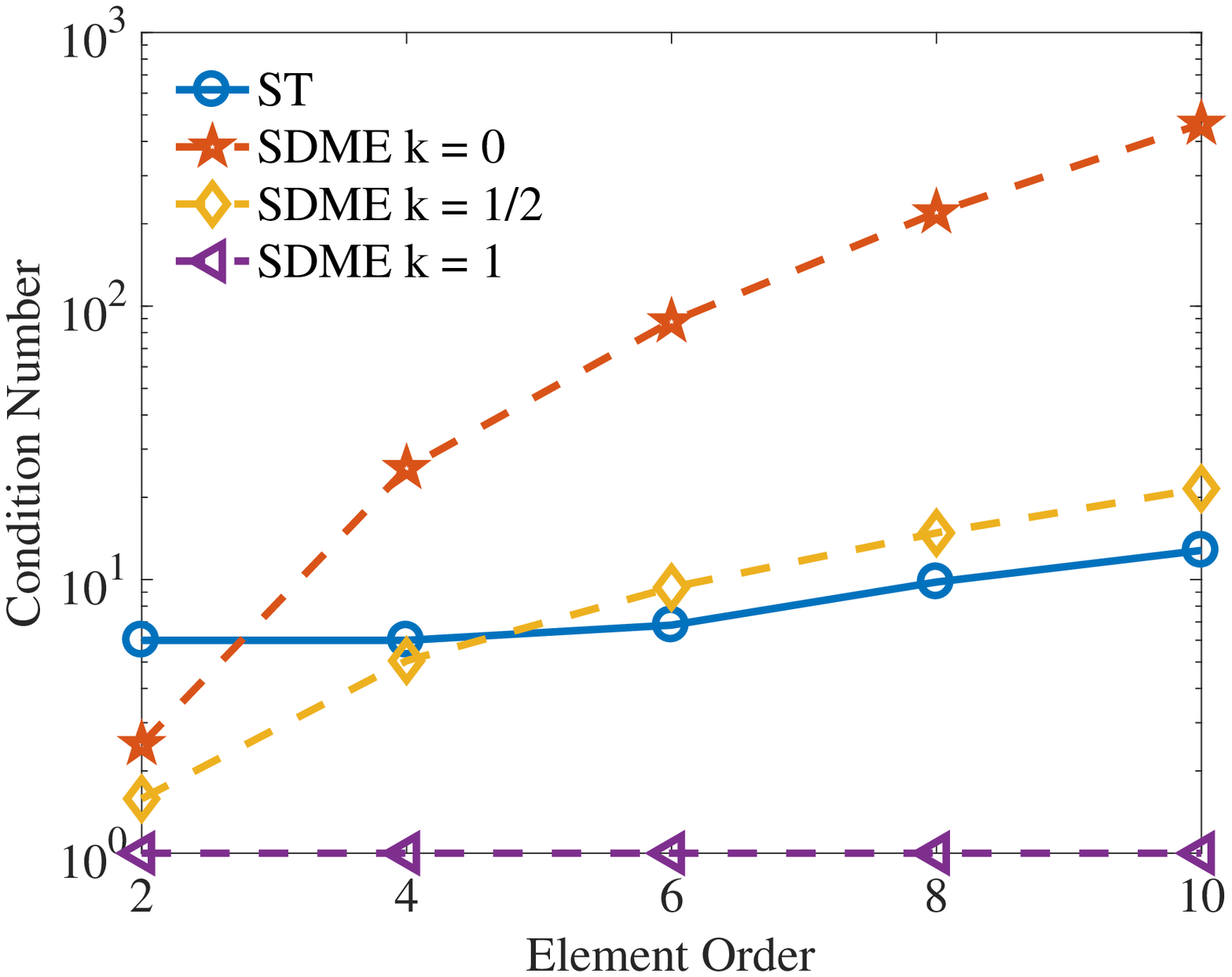}
                \caption{$\left[ K \right]$, ST and SDME-M bases.}
        \end{subfigure}%
        ~
        \begin{subfigure}[b]{0.325\textwidth}
                \includegraphics[width=\columnwidth]{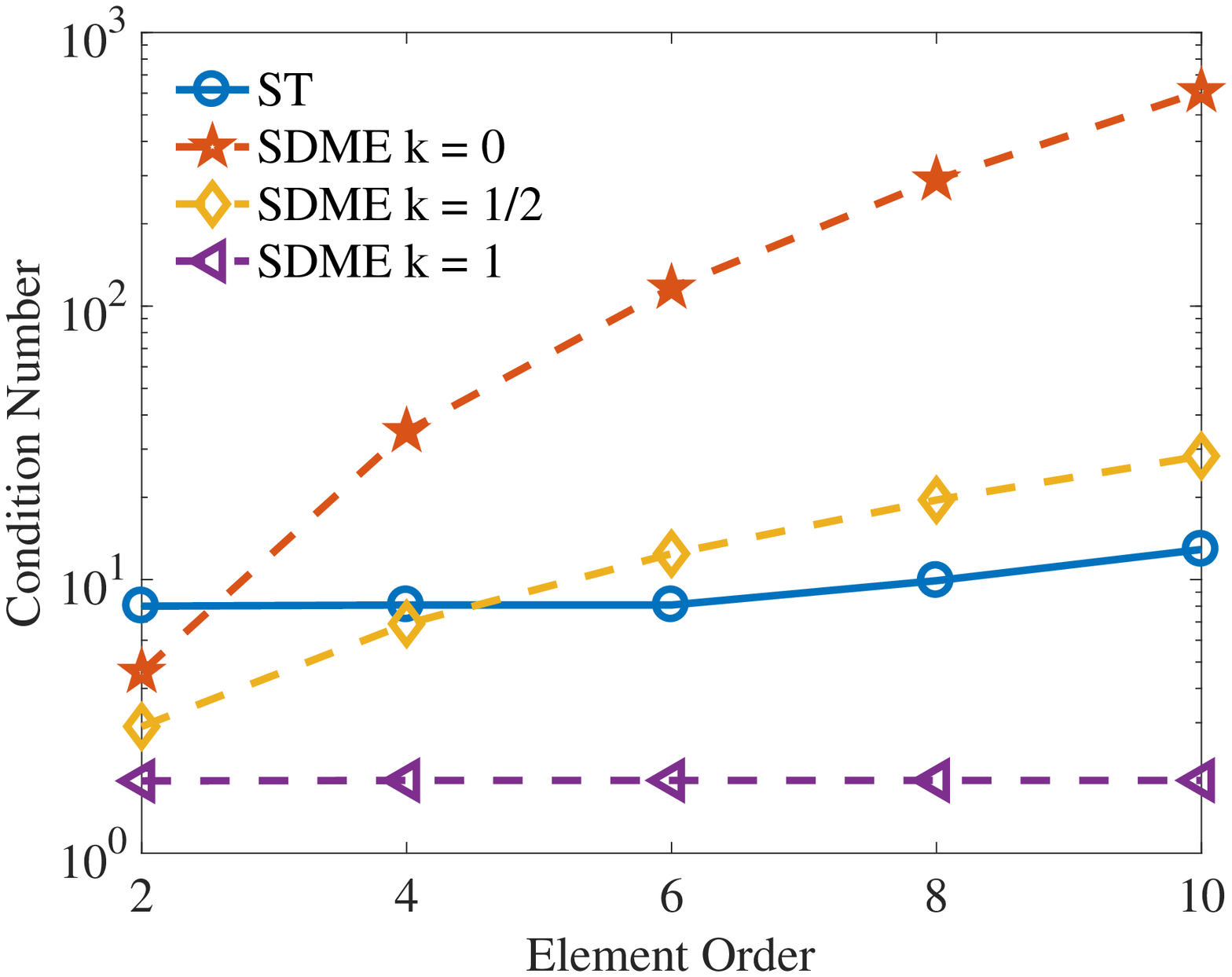}
                \caption{$\left[\hat{K}\right]$, ST and SDME-H bases.}
        \end{subfigure}
\caption{Condition numbers for the mass, stiffness and equivalent element matrices using
the standard (ST), SDME-M and SDME-H bases with $\lambda = 1$.}
\label{ConditionNumbers1D}
\end{figure}

{
Figure \ref{ConditionNumbers1D} illustrates the condition numbers for the 1D element local mass and
stiffness matrices for orders up to 10 obtained using the SDME-M basis and for the
local equivalent stiffness element matrices for the SDME-H basis. It may be observed that
parameter $k$ affects the conditioning. The element local mass matrices calculated with the SDME-M
basis and $k=\frac{1}{2}$  have condition numbers lower than the respective ones of the ST basis and
increasing slightly with the polynomial order. For $k=1$, the stiffness matrices calculated with
the SDME-M basis have constant condition numbers and equal to 1 for any polynomial order.
Again for $k=1$, the equivalent stiffness matrix with the SDME-H basis has almost constant condition
numbers for any order. These features will be similar for multidimensional elements and have
positive influence on the number of iterations for the conjugate gradient method
which will increase slightly for higher polynomial orders as shown in Section \ref{NumericalResults}.
Aspects related to the sparsity of element matrices and the influence of parameter $\lambda$ in
the conditioning of matrices calculated with the SDME-H basis are presented in \cite{CF_Rodrigues2014}.
}

The shape functions for  squares and hexahedra are obtained
using the tensor product of the {previously developed} one-dimensional functions, respectively,
in the local coordinate systems $\xi_{1}\times\xi_{2}$ and $\xi_{1}\times\xi_{2}\times\xi_{3}$
 \cite{GE_Karniadakis2005, ML_Bittencourt2007,ML_Bittencourt2014}:
\begin{eqnarray}
N_{i}(\xi_{1},\xi_{2}) &=& \varphi_{p}(\xi_{1})\varphi_{q}(\xi_{2})
\;\; (0\leq p,q \leq P), \label{ProdTensorialQuadrado} \\
N_{i}(\xi_{1},\xi_{2},\xi_{3}) &=& \varphi_{p}(\xi_{1})\varphi_{q}
(\xi_{2})\varphi_{r}(\xi_{3})\;\; (0\leq p,q,r \leq P),
\label{ProdTensorialHexaedro}
\end{eqnarray}
where $p$, $q$ and $r$ {represent the} tensor product indices associated with the
topological entities of the element; $P$ {denotes} the polynomial order in
directions $\xi_{1}$, $\xi_{2}$ and $\xi_{3}$; $i = 1,\ldots,(P+1)^{2}$
for squares and $i = 1,\ldots,(P+1)^{3}$ for hexahedra.
The SDME bases are hierarchical, conforming and continuous on the element boundaries.

{

It is possible to define procedures to construct the tensor indices $p$, $q$ and $r$ for
any polynomial order $P$. Observe that as the polynomial order increases, the number of
body shape functions of hexahedron increases very fast with the cubic power of $P$.
In this way, it is very important to construct the shape functions using the tensor
product of the one-dimensional functions, avoiding large memory demand.

The shape functions of squares are associated with the element topological entities, which include
four vertices $(V_1, V_2, V_3, V_4)$, four edges $(E_1,E_2, E_3, E_4)$, and one face $(F_1)$,
illustrated in Figure~\ref{enttop}. The indices $p$ and $q$ of Equation (\ref{ProdTensorialQuadrado})
are associated to the topological entities according to Figure~\ref{entindice}. The linear, quadratic
and cubic square standard elements are illustrated in Figure \ref{squareelems} in the local
coordinate system $\xi_1\times \xi_2$. Vertex and edge nodes/modes defines the boundary
modes while the face ones define the internal set. Schur complement of the element matrices
are calculated to condense the contribution of the internal to the boundary DOFs.

\begin{figure}[!htbp]
    \centering
    \subcaptionbox{Indices $p$ and $q$.\label{quad}}
    {\includegraphics{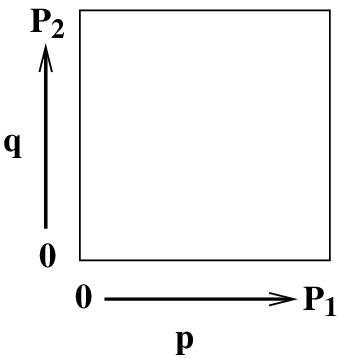}}\hspace*{5mm}%
    \subcaptionbox{Topological entities.\label{enttop}}
    {\includegraphics{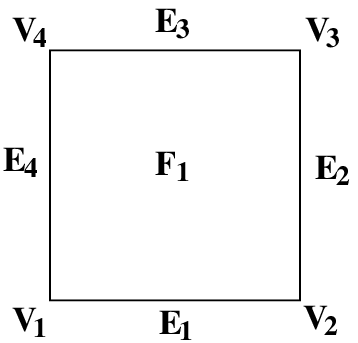}}\hspace*{5mm}
    \subcaptionbox{Entities and indices $p$ and $q$.\label{entindice}}
    {\includegraphics{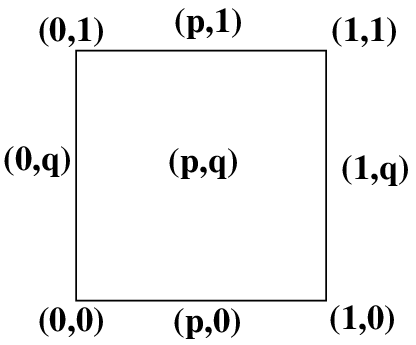}}
    \caption{Association between the topological entities and tensor indices $p$ and
    $q$ in the square  \cite{ML_Bittencourt2007}.}
    \label{Intquad}
\end{figure}

\begin{figure}[!htbp]
    \centering
    \subcaptionbox{Linear.}
    {\includegraphics{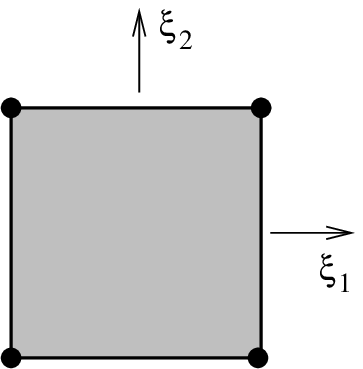}}\hspace*{5mm}%
    \subcaptionbox{Quadratic}
    {\includegraphics{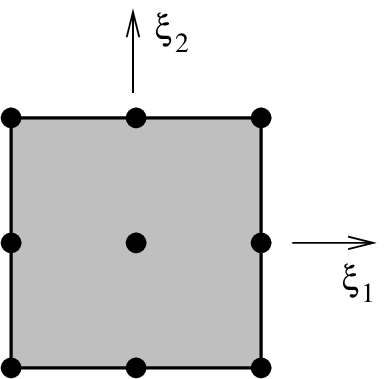}}\hspace*{5mm}
    \subcaptionbox{Cubic.}
    {\includegraphics{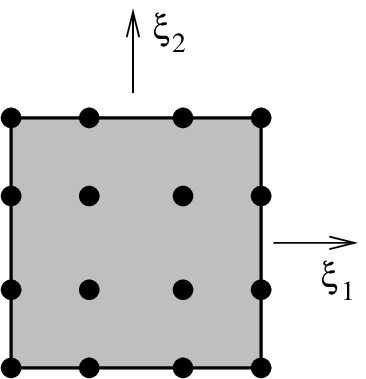}}
    \caption{Linear, quadratic and cubic squares  \cite{ML_Bittencourt2007}.}
    \label{squareelems}
\end{figure}
}

\section{Weak form of the boundary value problem}

Given a  deformation state $\varphi_k$ of a three-dimensional body and an arbitrary virtual displacement
$\delta\mathbf{u}$ kinematically admissible in the current position, we can write the total Lagrangian
description of the principle of virtual work for elastodynamics as: find the displacement vector field
$\mathbf{u}\in V_{t}$ such that for $\forall\delta \mathbf{u}\in V$ \cite{BATHELivro}
\begin{equation}
\delta\Pi_{s}\left(\varphi_k,\mathbf{u}\right)=\int_{\Omega_{0}}\rho_{0}\frac{\partial^{2}\mathbf{u}}{\partial t^{2}}\cdot\delta\mathbf{u}d\Omega+\int_{\Omega_{0}}\mathbf{S}\colon\delta\mathbf{E}d\Omega-
\int_{\Omega_{0}}\rho_{0}\mathbf{f}\cdot\delta\mathbf{u}d\Omega-
\int_{\Gamma_{\sigma}}\bar{\mathbf{t}}\cdot\delta\mathbf{u}d\Gamma_{\sigma},\label{eq:material virtual work}
\end{equation}
\noindent where $\Omega$ denotes the current domain occupied by the body with boundary
$\partial\Omega=\Gamma_{\sigma}\bigcup\Gamma_{u}$ and $\Gamma_{\sigma}\bigcap\Gamma_{u}=\emptyset$.
$\Gamma_{\sigma}$ and $\Gamma_{u}$ denote the Neumann and Dirichlet boundaries with prescribed tractions
($\bar{\mathbf{t}}$) and displacement ($\bar{\mathbf{u}}$) vector fields, respectively. {The term}
$\mathbf{S}$ is the second Piolla-Kirchhoff tensor, $\delta\mathbf{u}$ is the virtual
displacement vector field from the current position, $\delta\mathbf{E}$ is the associated virtual Green strain tensor,
$\rho_{0}$ is the mass density, $\mathbf{f}$ is the vector field of body forces. $V_{t}$ is the standard solution
space for elastodynamics and $V$ the test space, respectively, defined as
\begin{eqnarray}
V_{t}&=&\left\{ \mathbf{u}\in H^{\text{1}}\left(\varOmega\right)\colon\mathbf{u}=\bar{\mathbf{u}}\left(\mathbf{x},t\right)
\mbox{ for }\mathbf{x}\mbox{ on }\Gamma_{u}\right\} ,\label{eq:standard space of trial functions}\\
V&=&\left\{ \mathbf{w}\in H^{\text{1}}\left(\varOmega\right)\colon\mathbf{w}
\left(\mathbf{x}\right)=0\mbox{ for } \mathbf{x}\mbox{ on }\Gamma_{u}\right\} .\label{eq:space of test functions}
\end{eqnarray}
$H^{1}\left(\Omega\right)$ is the Hilbert space of all vector valued functions over $\Omega$ which
together with their first derivatives are square integrable over the domain.

We use the neo-Hookean hyperelastic material for describing the nonlinear elastic constitutive model. The
{corresponding} strain energy density function is defined as
\begin{equation}
\Psi=\frac{\mu}{2}\left[\mbox{tr}\left(\mathbf{C}\right)-3\right]-\mu\text{ln}J+\frac{\lambda}{2}\left(\text{ln}J\right)^{2},\label{eq:hyperelastic strain energy}
\end{equation}
\noindent where $\mathbf{C}$ is the right Cauchy-Green deformation tensor directly related to
$\mathbf{E}=\frac 12\left(\mathbf{C}-\mathbf{I}\right)$, $J=\det \mathbf{F}$, {where} $\mathbf{F}$  is the gradient of
deformation tensor and $\ln\left(J\right)$. {The terms} $\mu$ and $\lambda$ are the Lam\`{e}
parameters \citep{Bonet2008}. The constitutive equation for the second Piola-Kirchhoff stress is
\begin{equation}
\mathbf{S}=\frac{\partial\Psi}{\partial\mathbf{E}}=2\frac{\partial\Psi}{\partial\mathbf{C}}=\mu\left(\mathbf{I}-\mathbf{C}^{-1}\right)+\lambda\left(\text{ln}J\right)\mathbf{C}^{-1}.\label{eq:second Piolla_Kirchhoff relation}
\end{equation}

The first term in (\ref{eq:material virtual work}) is the virtual work of the inertia denoted by
$\delta\Pi^{ine}\left(\varphi_{k},\delta\mathbf{u}\right)$. We use $\delta\Pi^{int}\left(\varphi_{k},\delta\mathbf{u}\right)$
to denote the second term in (\ref{eq:material virtual work}), which represents the internal virtual work due to the
stresses and is nonlinear with respect to the displacement field. The third and fourth terms represent the
virtual work of  the external loads and denoted by $\delta\Pi^{ext}\left(\varphi_{k},\delta\mathbf{u}\right)$.
Therefore, Eq.(\ref{eq:material virtual work}) can be rewritten as
\begin{eqnarray}
\delta\Pi_{s}\left(\varphi_{k},\delta\mathbf{u}\right) & = & \delta\Pi^{ine}\left(\varphi_{k},\delta\mathbf{u}\right)+\delta\Pi^{int}\left(\varphi_{k},\delta\mathbf{u}\right)-\delta\Pi^{ext}\left(\varphi_{k},\delta\mathbf{u}\right)\label{eq:virtual work matrix 0}
 = \delta\hat{\mathbf{u}}^{T}\psi,\label{eq:virtual work matrix 2}
\end{eqnarray}
\noindent where $\psi$ is the residue vector.

{Equation} (\ref{eq:virtual work matrix 2}) can be linearized in the direction of a displacement increment
$\Delta\mathbf{u}$, using a first order {Taylor} expansion about the trial point
\begin{equation}
\delta\Pi_{s}\left(\varphi_{k},\delta\mathbf{u}\right)+D\delta\Pi_{s}\left(\varphi_{k},
\delta\mathbf{u}\right)\Delta\mathbf{u}=0.\label{eq:directional derivative}
\end{equation}

In the HOFEM, the element displacement field, virtual displacement field and material coordinates can be interpolated
in each element similarly to the standard FEM, respectively as
\begin{equation}
\mathbf{u}^{e}(\boldsymbol{\xi})=\sum_{i=1}^{N_{me}}N_{i}(\boldsymbol{\xi})\hat{\mathbf{u}}_{i},
\label{eq: interpolacao u}
\end{equation}
\begin{equation}
\delta\mathbf{u}^{e}(\boldsymbol{\xi})=\sum_{i=1}^{N_{me}}N_{i}(\boldsymbol{\xi})\delta\hat{\mathbf{u}}_{i},
\label{eq: interpolacao virtual u}
\end{equation}
\begin{equation}
\mathbf{X}^{e}(\boldsymbol{\xi})=\sum_{i=1}^{N_{me}}N_{i}(\boldsymbol{\xi})\hat{\mathbf{X}}_{i},
\label{eq: interpolacao material coord}
\end{equation}
\noindent where $N_{me}$ is the number of modes per element, $N_{i}$ are the shape functions and $\boldsymbol{\xi}$
is the local coordinate system (for square $\boldsymbol{\xi}=\xi_{1}\times\xi_{2}$) and
hexahedra $\boldsymbol{\xi}=\xi_{1}\times\xi_{2}\times\xi_{3}$.
$\hat{\mathbf{u}}_{i}$, $\delta\hat{\mathbf{u}}_{i}$ and $\hat{\mathbf{X}}_{i}$ are the expansion coefficients for the
displacement field, virtual displacement field and material coordinates, respectively.

The first term in (\ref{eq:directional derivative}) is the virtual work and the {derivative in the} second term gives rise to the
effective tangent stiffness matrix $\hat{\mathbf{K}}\left(\varphi_{k},\delta\mathbf{u}\right)$ using the previous
approximations. Therefore,
\begin{eqnarray}
D\delta\Pi_{s}\left(\varphi_{k},\delta\mathbf{u}\right)\Delta\mathbf{u} & = &
D\left[\delta\Pi^{ine}+\delta\Pi^{int}-\delta\Pi^{ext}\right]\Delta\mathbf{u} 
=\delta\hat{\mathbf{u}}^{T}\hat{\mathbf{K}}\Delta\mathbf{u}.\label{eq:effective stiffness 2}
\end{eqnarray}
Finally, we can rewrite (\ref{eq:directional derivative}) as
\begin{equation}
\delta\hat{\mathbf{u}}^{T}\left(\psi+\delta\hat{\mathbf{u}}^{T}\hat{\mathbf{K}}\Delta\right)=0.
\label{eq:directional derivative rew}
\end{equation}

Once the virtual displacements are arbitrary, we can solve (\ref{eq:directional derivative rew}) iteratively to find
$\Delta\mathbf{u}$. Detailed description on the derivation of the weak form of the BVP, including
the directional derivatives  $D\delta\Pi^{ine}\Delta\mathbf{u}$ and $D\delta\Pi^{int}\Delta\mathbf{u}$, can be
found in \cite{Bonet2008}.

\section{Non-linear elastodynamics}

In this section we present the equation for conservation of linear momentum in discrete form, and show the explicit
central difference and implicit Newmark time-integration schemes employed, with the respective expressions
after applying the Schur complement to condense the internal modes.

\subsection{Explicit time integration}

We consider the equation of motion (\ref{eq:virtual work matrix 2}) in discrete form for the current time $t_n$, neglecting
damping effects, for non-linear elastic problems{, which is} given by
\begin{equation}\label{Eq:Motion_explicit1}
	\mathbf{M}\mathbf{a}_n + \mathbf{R}\left(\mathbf{u}_n\right) = \mathbf{P}_n,
\end{equation}
with $\mathbf{M}$ denoting the global mass matrix in the reference configuration,
$\mathbf{a}_n=\ddot{\boldsymbol{u}}_n$ the global acceleration vector, $\mathbf{R}\left(\mathbf{u}_n \right)$ the global
internal force vector and $\mathbf{P}_n$ the global external load vector. The velocities $\mathbf{v}_n$
and accelerations  $\mathbf{a}_n$ can be approximated using the central-difference schemes in the following ways:
\begin{equation}\label{Eq:Vel_acel_explicit}
	\mathbf{v}_n = \frac{\mathbf{u}_{n+1} - \mathbf{u}_{n-1}}{2 \Delta t}, \quad \quad
	\mathbf{a}_n = \frac{\mathbf{u}_{n+1} - 2 \mathbf{u}_n + \mathbf{u}_{n-1}}{{\Delta t}^2}.
\end{equation}
Substituting the accelerations $\mathbf{a}_n$ from
Eq.(\ref{Eq:Vel_acel_explicit}) into Eq.(\ref{Eq:Motion_explicit1}) and rearranging the terms, we obtain
\begin{equation}\label{Eq:Motion_explicit2}
	\hat{\mathbf{M}} \mathbf{u}_{n+1} = \boldsymbol{\psi}_n + \hat{\mathbf{M}}\left( 2\mathbf{u}_n - 		\mathbf{u}_{n-1} \right),
\end{equation}
with,
\begin{equation}
	\hat{\mathbf{M}} = \frac{1}{{\Delta t}^2} \mathbf{M},
\end{equation}
\begin{equation}
	\boldsymbol{\psi}_n = \mathbf{P}_n - \mathbf{R}\left(\mathbf{u}_n\right).
\end{equation}
Therefore, we must solve Eq.(\ref{Eq:Motion_explicit2}) to determine the displacements at time $t_{n+1}$. Considering
initial conditions for the displacements and velocities ($\mathbf{u}_0$, $\mathbf{v}_0$ known), the initial condition for the
acceleration can be obtained by setting $t=t_0$ in Eq.(\ref{Eq:Motion_explicit1}). Therefore,
\begin{equation}
	\mathbf{a}_0 = \mathbf{M}^{-1} \boldsymbol{\psi}_0.
\end{equation}
The displacement for time $t=t_{-1}$ can be obtained with Eq.(\ref{Eq:Vel_acel_explicit}) and is given by
\begin{equation}
	\mathbf{u}_{-1} = \mathbf{u}_0 - \Delta t \mathbf{v}_0 + \frac{{\Delta t}^2}{2} \mathbf{a}_0.
\end{equation}

We can express Eq.(\ref{Eq:Motion_explicit2})  in terms of the boundary and internal modes in the following form:
\begin{equation} \label{Eq:Matrix_form_schur1}
\begin{bmatrix}
\hat{\mathbf{M}}_{bb} & \hat{\mathbf{M}}_{bi} \\
\hat{\mathbf{M}}^T_{bi} & \hat{\mathbf{M}}_{ii} \\
\end{bmatrix}
\begin{bmatrix}
\mathbf{u}_b \\ \mathbf{u}_i
\end{bmatrix}_{n+1} =
\begin{bmatrix}
\boldsymbol{\psi}_b \\ \boldsymbol{\psi}_i
\end{bmatrix}_n +
\begin{bmatrix}
\hat{\mathbf{M}}_{bb} & \hat{\mathbf{M}}_{bi} \\
\hat{\mathbf{M}}^T_{bi} & \hat{\mathbf{M}}_{ii} \\
\end{bmatrix}
\begin{bmatrix}
\mathbf{u}^*_b \\ \mathbf{u}^*_i
\end{bmatrix},
\end{equation}
with
\begin{equation}
\begin{bmatrix}
\mathbf{u}^*_b \\ \mathbf{u}^*_i
\end{bmatrix} =
\begin{bmatrix}
2\mathbf{u}_b \\ 2\mathbf{u}_i
\end{bmatrix}_n -
\begin{bmatrix}
\mathbf{u}_b \\ \mathbf{u}_i
\end{bmatrix}_{n-1}.
\end{equation}

Expanding Eq.(\ref{Eq:Matrix_form_schur1}) in terms of the boundary modes and dropping the subscripts for time, we obtain the following equations:
\begin{equation}\label{Eq:Schur_split1}
\hat{\mathbf{M}}_{bb} \mathbf{u}_b + \hat{\mathbf{M}}_{bi} \mathbf{u}_i = \boldsymbol{\psi}_b + \hat{\mathbf{M}}_{bb} \mathbf{u}^*_b + \hat{\mathbf{M}}_{bi} \mathbf{u}^*_i,
\end{equation}
\begin{equation} \label{Eq:Schur_split2}
\hat{\mathbf{M}}^T_{bi} \mathbf{u}_b + \hat{\mathbf{M}}_{ii} \mathbf{u}_i = \boldsymbol{\psi}_i + \hat{\mathbf{M}}^T_{bi} \mathbf{u}^*_b + \hat{\mathbf{M}}_{ii} \mathbf{u}^*_i.
\end{equation}
Multiplying Eq.(\ref{Eq:Schur_split2}) by $\hat{\mathbf{M}}^{-1}_{ii}$ and solving for
$\mathbf{u}_i$, we obtain
\begin{equation} \label{Eq:Internal_Schur_Explicit}
\mathbf{u}_i = \mathbf{u}^*_i +
\hat{\mathbf{M}}^{-1}_{ii} \left(
\boldsymbol{\psi}_i -
\hat{\mathbf{M}}^T_{bi} \mathbf{u}_b  +  \hat{\mathbf{M}}^T_{bi} \mathbf{u}^*_b \right).
\end{equation}
Substituting the above equation into Eq.(\ref{Eq:Schur_split1}) and rearranging the terms, we have
%
\begin{equation} \label{Eq:boundary_schur}
\hat{\mathbf{M}}^{sc}_{b} \mathbf{u}_b = \boldsymbol{\psi}^{sc}_b + \hat{\mathbf{M}}^{sc}_{b} \mathbf{u}^*_b,
\end{equation}
where
\begin{equation}
\hat{\mathbf{M}}^{sc}_{b} = \hat{\mathbf{M}}_{bb} - \hat{\mathbf{M}}_{bi} \hat{\mathbf{M}}^{-1}_i \hat{\mathbf{M}}^T_{bi},
\end{equation}
\begin{equation}
\boldsymbol{\psi}^{sc}_b = \boldsymbol{\psi}_b - \hat{\mathbf{M}}_{bi} \hat{\mathbf{M}}^{-1}_i \boldsymbol{\psi}_i.
\end{equation}

{Finally,} we can solve Eq.(\ref{Eq:boundary_schur}) to obtain
\begin{equation}\label{Eq:Boundary_schur_explicit}
	\mathbf{u}_b = \hat{\mathbf{M}}^{{sc}^{-1}}_{b} \boldsymbol{\psi}^{sc}_b + \mathbf{u}^*_b.
\end{equation}
Therefore, we perform the Schur complement on $\hat{\mathbf{M}}_b$ and $\boldsymbol{\psi}_b$, calculate
the coefficients of the boundary modes from Eq.(\ref{Eq:Boundary_schur_explicit}) and {then} recover the
coefficients of the internal modes using Eq.(\ref{Eq:Internal_Schur_Explicit}).

\subsection{Implicit time integration (Newmark)}

We consider the equilibrium equation (\ref{Eq:Motion_explicit1}) for the current time step $t_{n+1}$
\begin{equation}\label{Eq:Motion_implicit1}
\mathbf{M}\mathbf{a}_{n+1} + \mathbf{R}_{n+1} = \mathbf{P}_{n+1},
\end{equation}
where $\mathbf{M}$ is the global mass
matrix, $\mathbf{R}_{n+1}$ the global internal force vector dependent on the updated configuration with
coordinates $\mathbf{x}_{n+1}$, which in turn depend on the displacements $\mathbf{u}_{n+1}$. The term
$\mathbf{P}_{n+1}$ represents the global external nodal force vector. The terms $\mathbf{a}_{n+1}$ and
$\mathbf{v}_{n+1}$ respectively denote the global acceleration and velocity vectors at time step $t_{n+1}$.

We define the following residue force vector $\boldsymbol{\psi}_{n+1}$ at time step $t_{n+1}$:
\begin{equation}\label{Eq:residual1}
\boldsymbol{\psi}_{n+1} = \mathbf{M} \mathbf{a}_{n+1} + \mathbf{R}_{n+1} - \mathbf{P}_{n+1} = \mathbf{0}.
\end{equation}
The following approximations for the velocity and accelerations are used by the Newmark method \citep{Wriggers2008}:
\begin{eqnarray}\label{eq:newmark_acceleration}
\mathbf{a}_{n+1} & = b_1 \left(\mathbf{u}_{n+1}-\mathbf{u}_n \right) - b_2 \mathbf{v}_n - b_3 \mathbf{a}_n, \\
\mathbf{v}_{n+1} & = b_4 \left(\mathbf{u}_{n+1}-\mathbf{u}_n \right) - b_5 \mathbf{v}_n - b_6 \mathbf{a}_n,
\end{eqnarray}
with the coefficients
\[\arraycolsep=2.0pt\def\arraystretch{1.5}
\begin{array}{ccc}
b_1 = \frac{1}{g_1 \Delta t^2}, & b_2 = \frac{1}{g_1 \Delta t}, & b_3 = \frac{1-2 g_1}{2 g_1}, \\
b_4 = \frac{g_2}{g_1 \Delta t^2}, & b_5 = \left( 1- \frac{g_2}{g_1}\right), & b_6 = \left( 1 - \frac{g_2}{2 g_1} \right) \Delta t.
\end{array}
\]
{Here, we} choose $g_1 = 0.5$ to obtain quadratic convergence in time and $g_2 = 0.25$ for unconditional
stability. Substituting Eq.(\ref{eq:newmark_acceleration}) in Eq.(\ref{Eq:residual1}), we obtain
\begin{equation}\label{Eq:residual}
 \boldsymbol{\psi}_{n+1} = \mathbf{M} \left[ b_1 \left(\mathbf{u}_{n+1}-\mathbf{u}_n \right) - b_2 \mathbf{v}_n - b_3 \mathbf{a}_n \right] + \mathbf{R}_{n+1} - \mathbf{P}_{n+1} = \mathbf{0},
\end{equation}

The equilibrium system, Eq.(\ref{Eq:residual}), is linearized with the
Newton-Raphson method  using incremental global displacements defined as
\begin{equation}
	\mathbf{u}^{k+1}_{n+1} = \mathbf{u}^{k}_{n+1} + \Delta \mathbf{u}.
\end{equation}
Accordingly, the updated global coordinates are given by
\begin{equation}\label{Eq:Updated_coords}
	\mathbf{x}^{k+1}_{n+1} = \mathbf{x}_n + \mathbf{u}^{k+1}_{n+1},
\end{equation}
where the superscript $k+1$ refers to the current iteration of the Newton method.

The linearized form of Eq.(\ref{Eq:residual}) in {the} direction of a displacement increment $\Delta \mathbf{u}$ is given by the following system of equations:
\begin{equation} \label{eq:NR_explicit}
\left[ b_1 \mathbf{M} + \mathbf{K}^k_{T_{n+1}} \right] \Delta \mathbf{u} = -\mathbf{M}\left[b_1 (\mathbf{u}^k_{n+1} - \mathbf{u}_n) - b_2 \mathbf{v}_n - b_3 \mathbf{a}_n \right] - \mathbf{R}^k_{n+1} + \mathbf{P}_{n+1},
\end{equation}
{where the} terms $\mathbf{u}_n$, $\mathbf{v}_n$, $\mathbf{a}_n$ are {known} from the last converged time step $t_{n}$. The term $\mathbf{K}_T$ is the tangent stiffness matrix and is updated at each iteration $k$ along with the internal force vector.

Now we consider the application of the Schur complement for the system given by Eq.(\ref{eq:NR_explicit}). We will drop the scripts $n+1$ and $k$ for simplicity. Consider
Eq.(\ref{eq:NR_explicit}) rewritten in the following form:
\begin{equation}\label{Eq:Newmark_system_noschur}
\hat{\mathbf{K}} \boldsymbol{\Delta u} = \boldsymbol{\psi},
\end{equation}
where $\hat{\mathbf{K}}$ denotes the effective tangent stiffness matrix{, given by}
\begin{equation}
\hat{\mathbf{K}} = b_1 \mathbf{M} + \mathbf{K}^{T},
\end{equation}
and $\boldsymbol{\psi}$ represents the residual force vector
\begin{equation}
\boldsymbol{\psi} = -\mathbf{M} \mathbf{a} -  \mathbf{R}^k_{n+1} + \mathbf{P}_{n+1}.
\end{equation}

Differently from the explicit method, we apply the Schur complement directly on the equivalent system,
Eq.(\ref{Eq:Newmark_system_noschur}), since we work with an equivalent global matrix in this case.
The previous equation can be written in terms of boundary, internal and coupled matrix blocks as
\begin{equation} \label{Eq:Newmark_Matrix_form_schur1}
\begin{bmatrix}
\hat{\mathbf{K}}_{bb} & \hat{\mathbf{K}}_{bi} \\
\hat{\mathbf{K}}^T_{bi} &\hat{\mathbf{K}}_{ii} \\
\end{bmatrix}
\begin{bmatrix}
\Delta \mathbf{u}_b \\ \Delta \mathbf{u}_i
\end{bmatrix} =
\begin{bmatrix}
\boldsymbol{\psi}_b \\ \boldsymbol{\psi}_i
\end{bmatrix},
\end{equation}
where, after applying the Schur complement, we have
\begin{equation}
	\left( \hat{\mathbf{K}}_{bb} - \hat{\mathbf{K}}_{bi} \hat{\mathbf{K}}^{-1}_{ii}
	\hat{\mathbf{K}}^T_{bi} \right) \Delta \mathbf{u}_b = \left( \boldsymbol{\psi}_b - \hat{\mathbf{K}}_{bi} \hat{\mathbf{K}}^{-1}_{ii} \boldsymbol{\psi}_i \right),
\end{equation}
\begin{equation}
\Delta \mathbf{u}_i = \hat{\mathbf{K}}^{-1}_{ii} \left( \boldsymbol{\psi}_i - \hat{\mathbf{K}}^T_{bi} \Delta \mathbf{u}_b \right).
\end{equation}

We can rewrite the linearized equilibrium system (\ref{eq:NR_explicit}) using the nonlinear Newmark time integration scheme to {the} contact problem as \citep{AllanThesis,DiasBook1}
\begin{eqnarray} \label{eq:NR_explicit_contact}
\left[ b_1 \mathbf{M} + \mathbf{K}^k_{T_{n+1}} + \mathbf{K}_{T_{n+1}}^{c^{k}} \right] \Delta \mathbf{u} & = & -\mathbf{M}\left[b_1 (\mathbf{u}^k_{n+1} - \mathbf{u}_n) - b_2 \mathbf{v}_n - b_3 \mathbf{a}_n \right] \nonumber \\
& & - \mathbf{R}^k_{n+1} + \mathbf{P}_{n+1} - \mathbf{F}_{n+1}^{c^{k}},
\end{eqnarray}
where $\mathbf{K}^{c}_{T}$ and $\mathbf{F}^{c}$ are the global contact tangent stiffness matrix
and force vector  after the application of the high-order finite element approximations.

\section{Numerical results}
\label{NumericalResults}

{
In this section, we analyze the performance of the SDME-M and SDME-H modal bases compared to the
standard Jacobi modal  (ST) basis in terms of  the number of iterations and and computational time
for linear system solution, using the conjugate gradient method with the Gauss-Seidel (CGGS) and the
diagonal (CGD) pre-conditioners \cite{Axelsson1994}. The first example considers the static analysis
of a large strain problem with fabricated solution comparing the convergence behavior of
the bases. Section \ref{TransientFabSolution} presents analyses of transient problems
with fabricated solutions using explicit and implicit time integration to verify  spatial
convergence and second-order time rate. Section \ref{TransientFabSolution} shows the
explicit and implicit analyses of a 3D conrod submitted to a transient dynamic load calculated
from the pressure curve of a four stroke engine. The polynomial orders are increased and the
results in terms of number of iterations and speedup are presented. Sections
\ref{2D Impact problem} and \ref{3D Impact problem} presents the analyses of 2D and 3D
frictionless impact problems, respectively. In Section \ref{3D Impact problem}, results for
the Lagrange nodal basis with Gauss-Lobatto-Legendre collocation points are included.
All the examples use a Neohookean hyperelastic material model with Lagragian description.
}

\subsection{Static non-linear problem with large strains}

To verify the performance {and accuracy} of the standard and minimum energy bases, we consider
the cube domain with coordinates $0 \le x,y,z \le 1$ discretized using 8 hexahedra
and the fabricated solution with the following displacement components:
\begin{equation}
	u_x = 1.9 \sin(x) - x, \quad u_y = 0, \quad u_z = 0.
\end{equation}
The Young's modulus and Poisson ratio are respectively $E = 1000\, Pa$ and $\nu = 0.3$.

We consider the ST, SDME-M and SDME-H ($\lambda=100$) bases with $k = 0.5$. We {compute} the
average number of iterations of the conjugate gradient method with diagonal preconditioner (CGD)
and time (for linear system solution) per Newton-Raphson iteration. The CGD tolerance is {chosen as} $10^{-12}$ and
the Newton solver tolerance is {set to} $10^{-8}$. We perform the Schur complement on the tangent stiffness matrix and
residual (out-of-balance) force vector. We also considered isoparametric mapping.

{The obtained spectral accuracy results are presented in Table \ref{Tab:Static_error}, which are identical for all employed bases}. From Table \ref{Tab:NH_average_IT_CGD}, we observe that the standard basis require  $7.85$ times
more iterations when compared to the SDME-M basis for $P = 8$. The same ratio is obtained for the average time in
Table \ref{Tab:NH_average_Time_CGD}. The SDME-H basis provided better performance than the SDME-M basis for
$P<8$.

\begin{table}[H]

\centering
\caption{$L_2$ error norms for the displacement components of the static analysis using the ST, SDME-M and SDME-H bases.}
\label{Tab:Static_error}
\begin{tabular}{@{}ccccc@{}}
\toprule
\multirow{2}{*}{Order} & \multirow{2}{*}{\begin{tabular}[c]{@{}c@{}}Number \\ of DOFs\end{tabular}} & \multicolumn{3}{c}{$L_2$ error} \\ \cmidrule(l){3-5} &    & $u_x$            & $u_y$            & $u_z$           \\ \cmidrule(r){1-5}
2	&	300  & 1.09e-03 &  1.83e-04 & 1.68e-04	 \\
4   &  1944  & 7.01e-07 &  1.22e-07 & 1.34e-07   \\
6   &  6084  & 4.51e-09 &  8.47e-10 & 7.98e-10   \\
8   & 13872  & 8.15e-11 &  1.21e-11 & 1.20e-11   \\ \bottomrule
\end{tabular}
\end{table}

\begin{table}[H]
\centering
\caption{Average numbers of CGD iterations per Newton iteration, total of 5 Newton iterations for convergence,
for the static problem of fabricated solution.}
\label{Tab:NH_average_IT_CGD}
\begin{tabular}{@{}crccc@{}}
\toprule
\multirow{2}{*}{Order} & \multirow{2}{*}{\begin{tabular}[c]{@{}c@{}}Number \\ of DOFs\end{tabular}} & \multicolumn{3}{r}{Average number of CGD iterations} \\ \cmidrule(l){3-5}
                        &                                                                            & ST            & SDME-M            & SDME-H           \\ \cmidrule(r){1-5}
2                       & 300                                                                        & 75.8             & 43.2                 & 41.2                \\
4                       & 1\,944                                                                       & 249.4             & 67.2                 & 57.0                \\
6                       & 6\,084                                                                       & 444.4            & 81.8                 & 64.0                \\
8                       & 13\,872                                                                      & 660.0             & 84.0                 & 95.0                \\ \bottomrule
\end{tabular}
\end{table}

\begin{table}[H]
\centering
\caption{Average time for CGD solution per Newton iteration for the static problem of
fabricated solution.}
\label{Tab:NH_average_Time_CGD}
\begin{tabular}{@{}crccc@{}}
\toprule
\multirow{2}{*}{Order} & \multirow{2}{*}{\begin{tabular}[c]{@{}c@{}}Number \\ of DOFs\end{tabular}} & \multicolumn{3}{c}{Average time for CGD solution [s]} \\ \cmidrule(l){3-5}
                        &                                                                            & ST            & SDME-M            & SDME-H           \\ \cmidrule(r){1-5}
2                       & 300                                                                        & 0.0042             & 0.0024                 & 0.0023                \\
4                       & 1\,944                                                                       & 0.2230             & 0.0598                 & 0.0501                \\
6                       & 6\,084                                                                       & 2.0236            & 0.3760                & 0.2464                \\
8                       & 13\,872                                                                      & 9.7028             & 1.2363                 & 1.4031                \\ \bottomrule
\end{tabular}
\end{table}

The same analysis was performed using the conjugate gradient method with the
Gauss-Seidel preconditioner (CGGS). The results are presented in Tables
\ref{Tab:NH_average_IT_CGGS} and \ref{Tab:NH_average_Time_CGGS}, showing a
better performance of the SDME-H basis for all polynomial orders.

\begin{table}[H]
\centering
\caption{Average number of CGGS iterations per Newton iteration, total of 5 Newton iterations for convergence,
for the static problem of fabricated solution.}
\label{Tab:NH_average_IT_CGGS}
\begin{tabular}{@{}crccc@{}}
\toprule
\multirow{2}{*}{Order} & \multirow{2}{*}{\begin{tabular}[c]{@{}c@{}}Number \\ of DOFs\end{tabular}} & \multicolumn{3}{c}{Average number of CGGS iterations} \\ \cmidrule(l){3-5}
                        &                                                                            & ST           & SDME-M            & SDME-H           \\ \cmidrule(r){1-5}
2                       & 300                                                                        & 61.2         &    49.6              & 47.8                \\
4                       & 1\,944                                                                       & 132.0        &  60.4                & 52.6                \\
6                       & 6\,084                                                                       & 218.2        &  66.4             & 51.6                \\
8                       & 13\,872                                                                      & 313.0             &  65.4                & 56.6                 \\ \bottomrule
\end{tabular}
\end{table}

\begin{table}[H]
\centering
\caption{Average time for CGGS solution per Newton iteration for the static problem of
fabricated solution.}
\label{Tab:NH_average_Time_CGGS}
\begin{tabular}{@{}crccc@{}}
\toprule
\multirow{2}{*}{Order} & \multirow{2}{*}{\begin{tabular}[c]{@{}c@{}}Number \\ of DOFs\end{tabular}} & \multicolumn{3}{c}{Average time for CGGS solution [s]} \\ \cmidrule(l){3-5}
                        &                                                                            & ST            & SDME-M            & SDME-H           \\ \cmidrule(r){1-5}
2               & 300                                                                        & 0.0031        &    0.0026              & 0.0023                \\
4                       & 1\,944                                                                       & 0.1037        &  0.0496                &  0.0419               \\
6                       & 6\,084                                                                       & 0.8952        &   0.2856              &   0.2389              \\
8                       & 13\,872                                                                      & 4.2829              &  0.9190                & 0.7569                \\ \bottomrule
\end{tabular}
\end{table}

\subsection{Transient nonlinear problems with large strains}
\label{TransientFabSolution}

\subsubsection{Explicit time integration}

In this case, we consider the same mesh of the previous example, time interval $t= [0,\, 0.25]\,s$ and
the following fabricated solution:
\begin{equation} \label{eq:large_field}
	u_x = \sin \left( \frac{\pi}{2} x \right) \sin(2\pi t), \quad u_y = 0, \quad u_z = 0,
\end{equation}
which gives $u_x = 1.0$\, m for $x = 1.0$ and $t = 0.25\, s$. The material properties  are
$E = 1000\, Pa$, $\nu = 0.3$ and $\rho = 1\, kg/m^3$.
Homogeneous Dirichlet boundary conditions are applied in all displacement directions of the face with
coordinate ${x} = 0$.

We used the CFL condition to estimate the number of time steps as \cite{Wriggers2008}:
\begin{equation}
	\Delta t \le \delta \frac{h}{c_L},
\end{equation}
where $h$ represents the element size (considered as the edge length in the initial configuration for this problem), $\delta$ denotes a constant in the range $0.2 < \delta < 0.9$ and
\begin{equation}
	c_L = \frac{3K(1-\nu)}{\rho (1 + \nu)}, \qquad K = \frac{E}{3(1-2\nu )}.
\end{equation}
Considering the given material properties and $\delta = 0.85$, we obtained
$\Delta t = 3.33 \times 10^{-4}$ and $N = 800$ times steps for the analysis
with a single load step.

{
Initially, we performed spatial convergence analysis by increasing the polynomial approximation
orders. Tables \ref{Tab:Transient_NH_L2errorHexa_ST} and \ref{Tab:Transient_NH_L2errorHexa_SDME_M}
show the $L_2$ error norms for the displacement components using  the ST and SDME-M bases, respectively.
We may observe spectral spatial convergence rates using the two bases. Furthermore, we compare the employed bases with a standard Lagrange basis with a diagonal mass matrix in Table \ref{Tab:Transient_NH_L2errorHexa_Lagrange}. We observe that the Lagrange basis with a diagonal mass matrix yields slightly higher errors for the $x$-direction due to the decreased number of integration points.}
\begin{table}[H]

\centering
\caption{$L_2$ error norms for the displacement components of explicit analysis using the  ST basis.}
\label{Tab:Transient_NH_L2errorHexa_ST}
\begin{tabular}{@{}ccccc@{}}
\toprule
\multirow{2}{*}{Order} & \multirow{2}{*}{\begin{tabular}[c]{@{}c@{}}Number \\ of DOFs\end{tabular}} & \multicolumn{3}{c}{$L_2$ error} \\ \cmidrule(l){3-5}
                        &                                                                            & $u_x$            & $u_y$            & $u_z$           \\ \cmidrule(r){1-5}
2					  &	276 	      &	2.08e-3														&	2.45e-4				&	2.26e-4        					 \\

4                     & 1246                                                                       & 3.38e-6             & 5.16e-7         & 5.61e-7         \\

6                   & 3084                                                                       &              4.21e-8				& 6.09e-9	      & 5.76e-9
\\ \bottomrule
\end{tabular}
\end{table}
\begin{table}[H]
\centering

\caption{$L_2$ error norms for the displacement components of explicit analysis using the  SDME-M basis.}
\label{Tab:Transient_NH_L2errorHexa_SDME_M}
\begin{tabular}{@{}ccccc@{}}
\toprule
\multirow{2}{*}{Order} & \multirow{2}{*}{\begin{tabular}[c]{@{}c@{}}Number \\ of DOFs\end{tabular}} & \multicolumn{3}{c}{$L_2$ error} \\ \cmidrule(l){3-5}
                        &                                                                            & $u_x$            & $u_y$            & $u_z$           \\ \cmidrule(r){1-5}
2					  &	276 	      &	2.09e-3														&	2.53e-4				&	2.34e-4        					 \\

4                     & 1246                                                                       & 3.36e-6             & 5.06e-7         & 5.51e-7         \\

6                   & 3084                                                                       &              4.21e-8				& 6.08e-9	      & 5.76e-9
\\ \bottomrule
\end{tabular}
\end{table}

\begin{table}[H]

\centering
\caption{$L_2$ error norms for the displacement components of the explicit analysis utilizing a Lagrange basis yielding a diagonal mass matrix.}
\label{Tab:Transient_NH_L2errorHexa_Lagrange}
\begin{tabular}{@{}ccccc@{}}
\toprule
\multirow{2}{*}{Order} & \multirow{2}{*}{\begin{tabular}[c]{@{}c@{}}Number \\ of DOFs\end{tabular}} & \multicolumn{3}{c}{$L_2$ error} \\ \cmidrule(l){3-5} &    & $u_x$            & $u_y$            & $u_z$           \\ \cmidrule(r){1-5}
2	&	276  & 2.12e-03 &  2.78e-04 & 2.78e-04	 \\
4   &  1246  & 3.75e-06 &  4.44e-07 & 4.44e-07   \\
6   &  3084  & 5.29e-08 &  5.50e-09 & 5.50e-09  \\
 \bottomrule
\end{tabular}
\end{table}

We compared the performance of the ST and SDME-M bases for the average number of iterations
and average time using the conjugate gradient method with the Gauss-Seidel
(CGGS) preconditioner  to solve the linear system of equations with tolerance of $10^{-12}$.
The results for the number of iterations are presented in Table
\ref{Tab:Transient_NH_NIT_Hexa}  and the computational times per time step
are given in Table \ref{Tab:Transient_NH_time_Hexa}. The SDME-M basis improved
the standard Jacobi modal basis about 18 times.
\begin{table}[H]
\centering
\caption{Average number of iterations for convergence using the CGGS method for the fabricated solution and explicit time integration.}
\label{Tab:Transient_NH_NIT_Hexa}
\begin{tabular}{@{}crcc@{}}
\toprule
Order & ST 	& SDME-M & Ratio ST/SDME \\ \midrule
2      & 68.37  & 9.99        &  6.84  \\
4      & 172.19 & 7.97        &  21.60  \\
6      & 297.48 & 16.48       &  18.05  \\ \bottomrule
\end{tabular}
\end{table}
\begin{table}[H]
\centering
\caption{Average time per time step using the CGGS methods for the ST and
SDME-M bases for the fabricated solution and explicit time integration. }
\label{Tab:Transient_NH_time_Hexa}
\begin{tabular}{@{}cccc@{}}
\toprule
Order & ST [s] 	& SDME-M [s] & Speedup \\ \midrule
2      & 0.0153 	&  0.0026    & 5.88   \\
4      & 0.6279 	&  0.0363    & 17.30   \\
6      & 5.7101 	&  0.3537    & 16.43   \\ \bottomrule
\end{tabular}
\end{table}

%

\subsubsection{Implicit time integration}

\begin{figure}[H]
        \centering
        \begin{subfigure}[b]{0.45\textwidth}
				\includegraphics[width=\columnwidth]{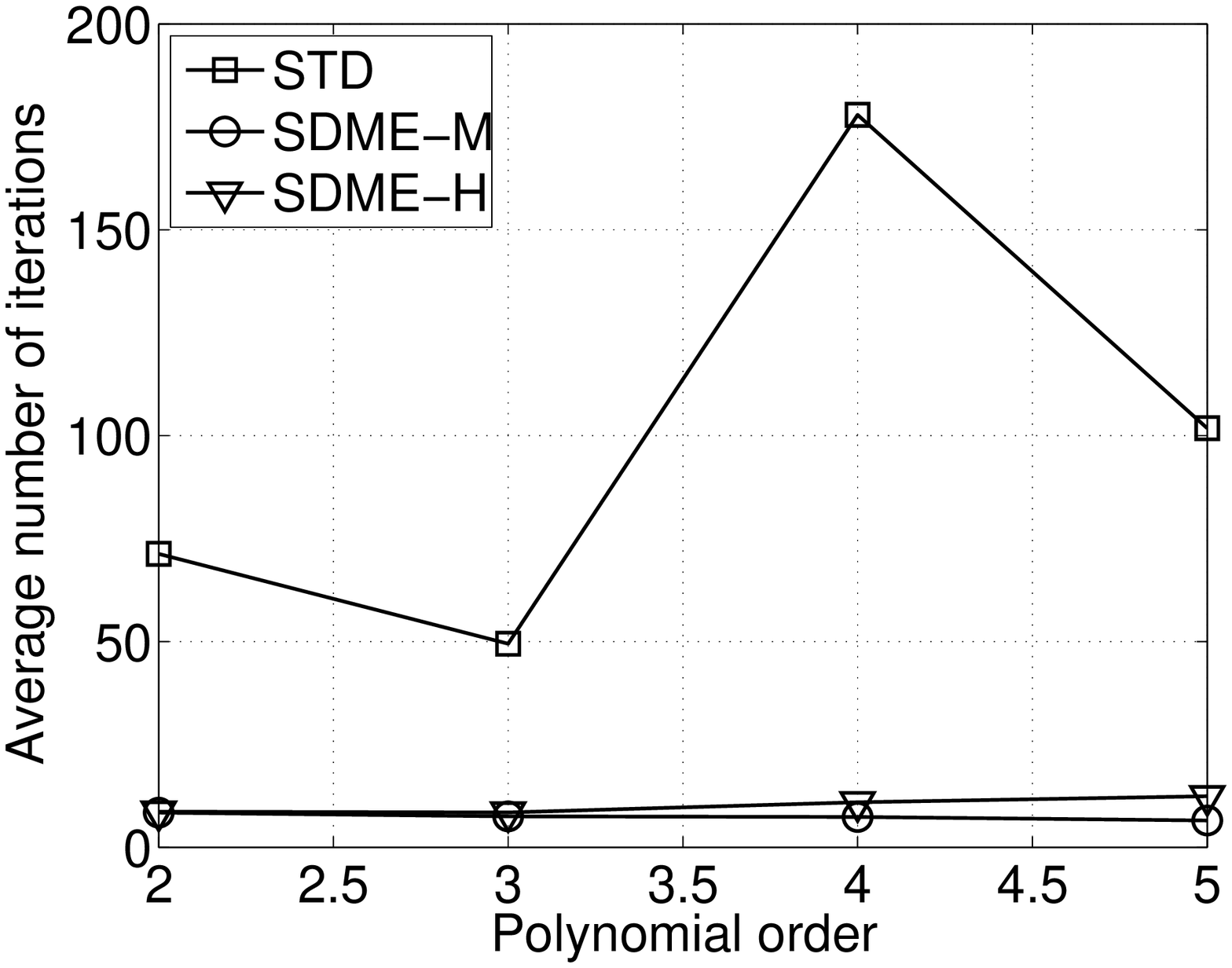}
				\caption{Average number of iterations.}

        \end{subfigure}%
        ~
        \begin{subfigure}[b]{0.45\textwidth}
                \includegraphics[width=\columnwidth]{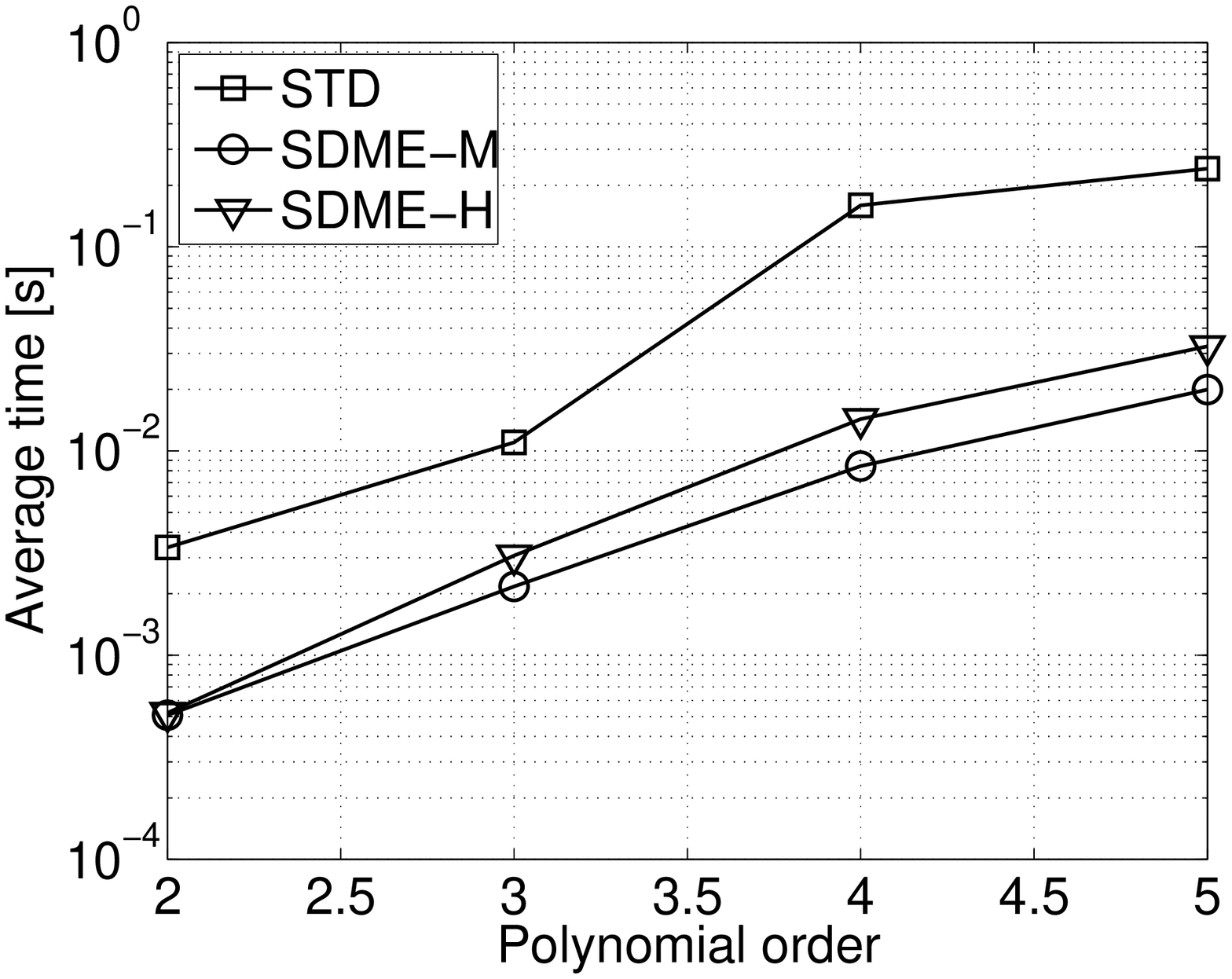}
                \caption{Average time.}
        \end{subfigure}

\caption{Average number of iterations for linear system solution using CGGS method in terms of the polynomial order, $\Delta t = 3.90 \times 10^{-4}\, s$ (a); Average time for linear system solution in terms of the polynomial order (b).\label{fig:Newmark_performance}}
\end{figure}

\begin{figure}[H]
        \centering
		\includegraphics[width=0.45\columnwidth]{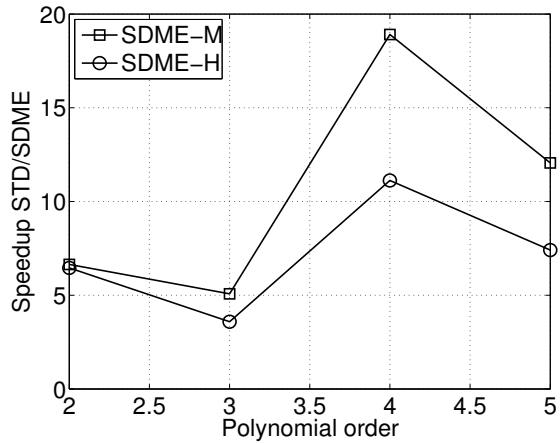}

\caption{Speedup ratio for the computation time to solve the linear system of equations between the standard Jacobi basis and the minimum energy bases SDME-M and SDME-H with $k = 0.5$ and $\lambda = 100$.\label{fig:Newmark_speedup}}
\end{figure}

For the implicit Newmark integration method, we first consider the following fabricated solution for $u(x,t)$:
\begin{equation}
	u_x = x^4 \sin(2\pi t), \quad u_y = 0, \quad u_z = 0.
\end{equation}
The total time is $T=0.025\, s$ and the solution gives $u_x = 0.157\, m$ for $x = 1.0\, m$ and $t = 0.025\, s$. The material properties and boundary conditions are the same of the previous example. The tolerance for the residue norm in the Newton method is $10^{-12}$.
%
%
We tested the performance of bases ST, SDME-M, SDME-H with $k=0.5$ and $\lambda = 100$ in terms of the average
number of iterations, average times and speedup. using the CGGS method. The results are presented in Figs.
\ref{fig:Newmark_performance} and \ref{fig:Newmark_speedup}. We observe that the SDME-M basis performed
better than the SDME-H basis, with a speedup up to 19 with polynomial order $P=4$. The SDME-H basis
achieved at least a speedup ratio of 3 compared to the ST basis as illustrated in Fig. \ref{fig:Newmark_speedup}.

\begin{figure}[H]
        \centering
        \begin{subfigure}[b]{0.45\textwidth}
				\includegraphics[width=\columnwidth]{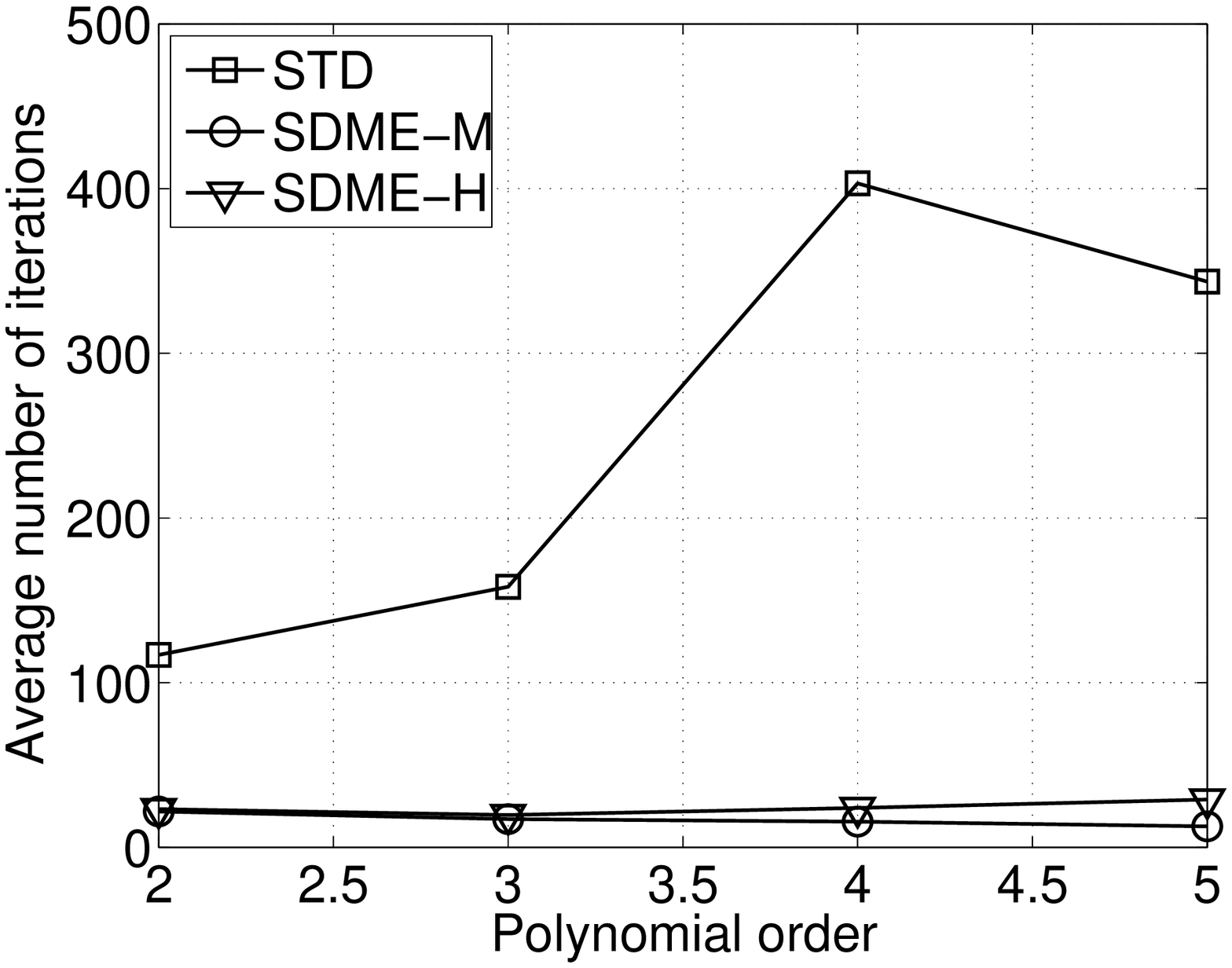}
				\caption{Average number of iterations.}
        \end{subfigure}%
        ~
        \begin{subfigure}[b]{0.45\textwidth}
                \includegraphics[width=\columnwidth]{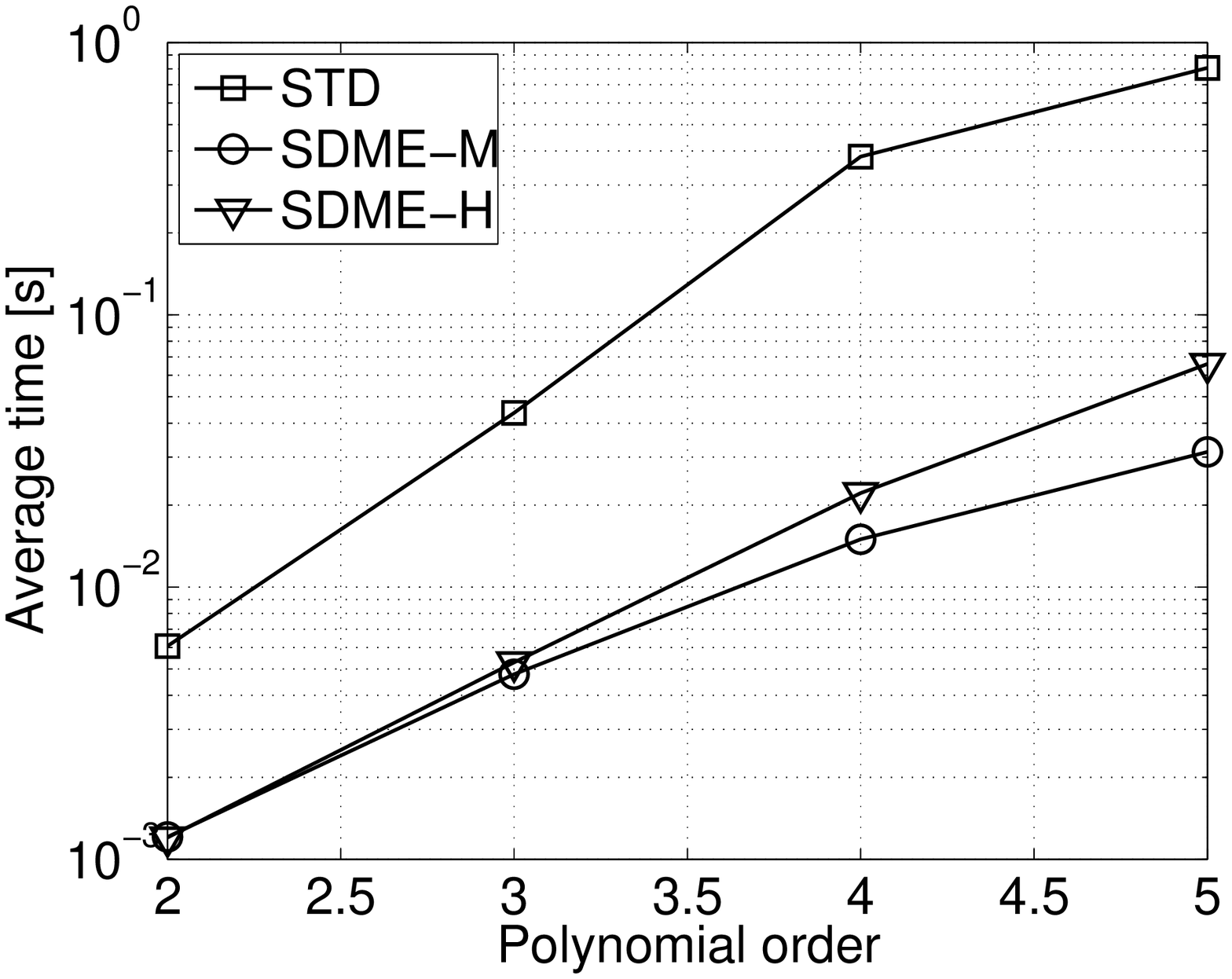}
                \caption{Average time.}
        \end{subfigure}

\caption{Average number of iterations for linear system solution using the CGD method in terms of the polynomial order, $\Delta t = 3.90 \times 10^{-4}\, s$ (a); Average time for linear system solution in terms of the polynomial order (b).\label{fig:Newmark_performance_Poly_CGD}}
\end{figure}

\begin{figure}[H]
        \centering
		\includegraphics[width=0.45\columnwidth]{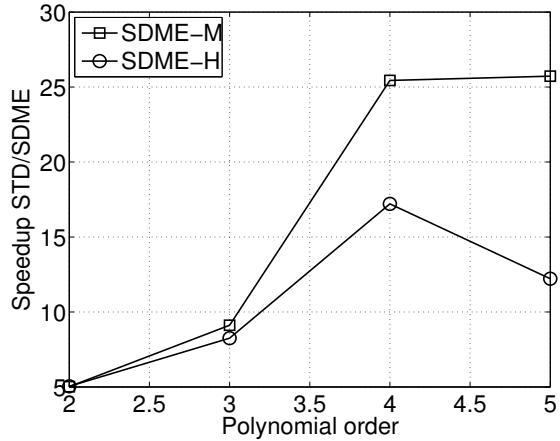}
\caption{Speedup ratio for the computation time to solve the linear system of equations between the standard Jacobi basis and the minimum energy bases SDME-M  and SDME-H with $k = 0.5\, \lambda = 100$.\label{fig:Newmark_speedup_Poly_CGD}}
\end{figure}

We also  considered the solution using the CGD method. The results are presented in Figs.
\ref{fig:Newmark_performance_Poly_CGD} and \ref{fig:Newmark_speedup_Poly_CGD}. Similarly to the CGGS
preconditioner, both minimum energy bases performed much better than the ST basis, with speedups up to 26 for
the SDME-M basis. In general, the speedup achieved by the SDME-M basis was also larger than the SDME-H basis.

The solution given in Eq.(\ref{eq:large_field}) was also considered and the results are presented in
Tables \ref{Tab:TransientNmk_NH_L2errorHexa_SDME_H}, \ref{Tab:TransientNmk_NH_NIT_Hexa} and
\ref{Tab:TransientNmk_NH_time_Hexa}. {The $L_2$ error norms for the displacement
components in Table \ref{Tab:TransientNmk_NH_L2errorHexa_SDME_H} are with the same order
of magnitude to those ones of Tables \ref{Tab:Transient_NH_L2errorHexa_ST} and
\ref{Tab:Transient_NH_L2errorHexa_SDME_M} }. We observe that the performance of the SDME-H basis
is slightly better than the SDME-M for $P=2,4$ and the SDME-M basis is better for $P=6$.

\begin{table}[H]
\centering

\caption{$L_2$ error norm using the Newmark method and SDME-H basis for 8 hexahedra, $\Delta t = 2\times 10^{-3}\, s$.}
\label{Tab:TransientNmk_NH_L2errorHexa_SDME_H}
\begin{tabular}{@{}ccccc@{}}
\toprule
\multirow{2}{*}{Order} & \multirow{2}{*}{\begin{tabular}[c]{@{}c@{}}Number \\ of DOFs\end{tabular}} & \multicolumn{3}{c}{$L_2$ error} \\ \cmidrule(l){3-5}
                        &                                                                            & $u_x$            & $u_y$            & $u_z$           \\ \cmidrule(r){1-5}
2					  &	276 	      &	2.09e-3														&	2.54e-4				&	2.35e-4        					 \\

4                     & 1246                                                                       & 3.40e-6             & 5.09e-7         & 5.54e-7         \\

6                   & 3084                                                                       &              9.17e-8				& 2.23e-8	      & 2.22e-8         \\ \bottomrule
\end{tabular}
\end{table}

\begin{table}[H]
\centering

\caption{Average number of iterations for convergence using the CGD method for the ST, SDME-M ($k=0.5$), and SDME-H ($k=0.5,\, \lambda = 100$) bases and Newmark integration,
$\Delta t = 2\times 10^{-3}\, s$. 
}
\label{Tab:TransientNmk_NH_NIT_Hexa}
\begin{tabular}{@{}cccccc@{}}
\toprule
Order & ST 	 & SDME-M & SDME-H & Ratio ST/SDME-M & Ratio ST/SDME-H \\ \midrule
2      &  119.65 & 21.40  & 20.86 & 5.59 & 5.74 \\
4      &  352.05 & 19.45  & 17.94 & 18.10 & 19.62 \\
6      &  610.83 & 32.28  & 42.01 & 18.92 & 14.54 \\ \bottomrule
\end{tabular}
\end{table}

\begin{table}[H]
\centering

\caption{Average time per linear system solution using the CGD method for
the ST, SDME-M ($k=0.5$), and SDME-H ($k=0.5,\, \lambda = 100$) bases and Newmark integration,
$\Delta t = 2\times 10^{-3}\, s$.}
\label{Tab:TransientNmk_NH_time_Hexa}
\begin{tabular}{@{}cccccc@{}}
\toprule
Order & ST (s)    & SDME-M (s) & SDME-H (s) & Speedup           & Speedup \\ \midrule
2      & 0.0066    & 0.0011     & 0.0012     & 6.000             & 5.500  \\
4      & 0.3167	   & 0.0185     & 0.0176     & 17.119            & 17.994 \\
6      & 2.8062    & 0.1561     & 0.1963     & 17.977            & 14.295 \\ \bottomrule
\end{tabular}
\end{table}

\begin{figure}[H]
        \centering
        \begin{subfigure}[b]{0.49\textwidth}
				\includegraphics[width=\columnwidth]{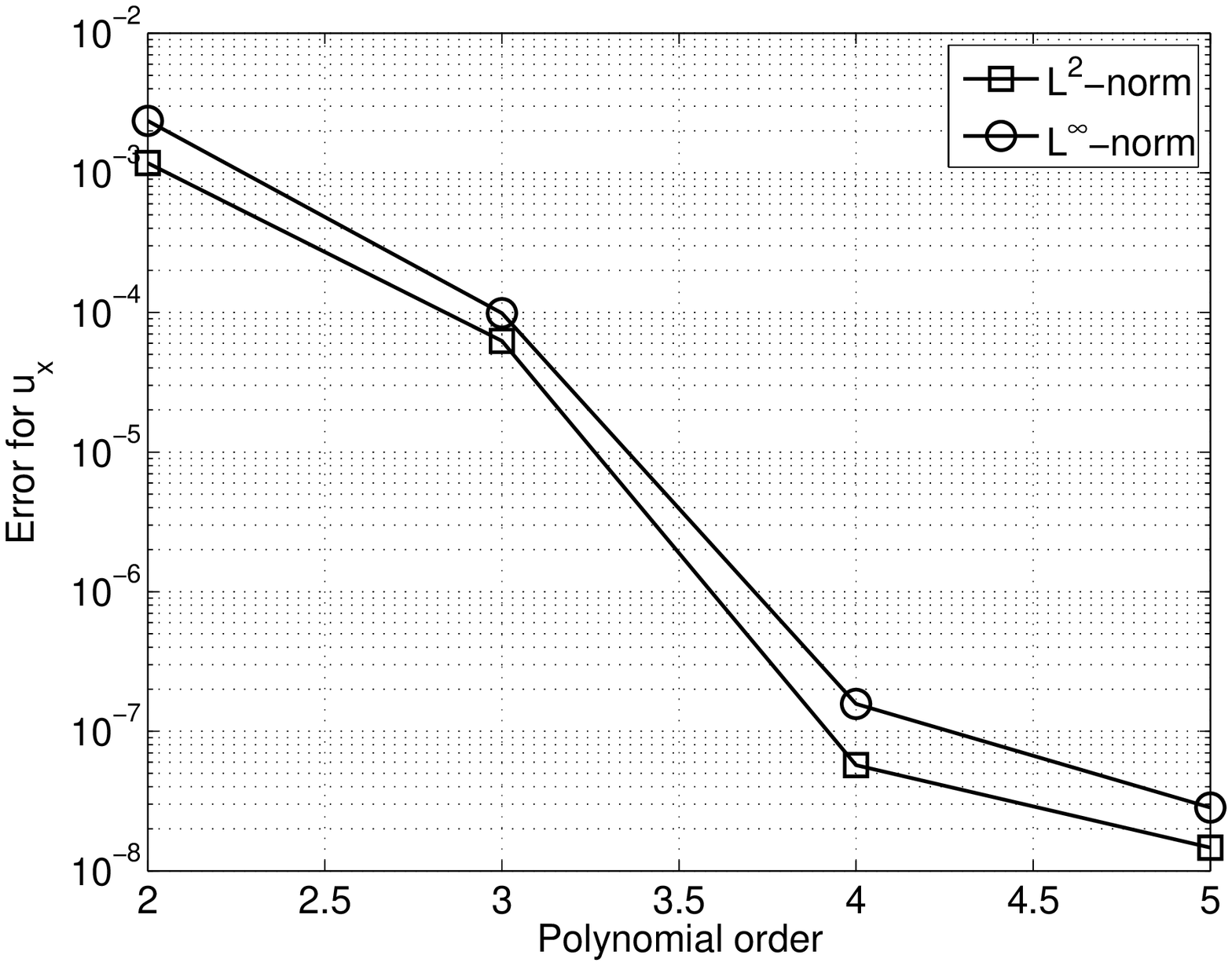}
				\caption{Spatial convergence.}
        \end{subfigure}%
        ~
        \begin{subfigure}[b]{0.49\textwidth}
                \includegraphics[width=\columnwidth]{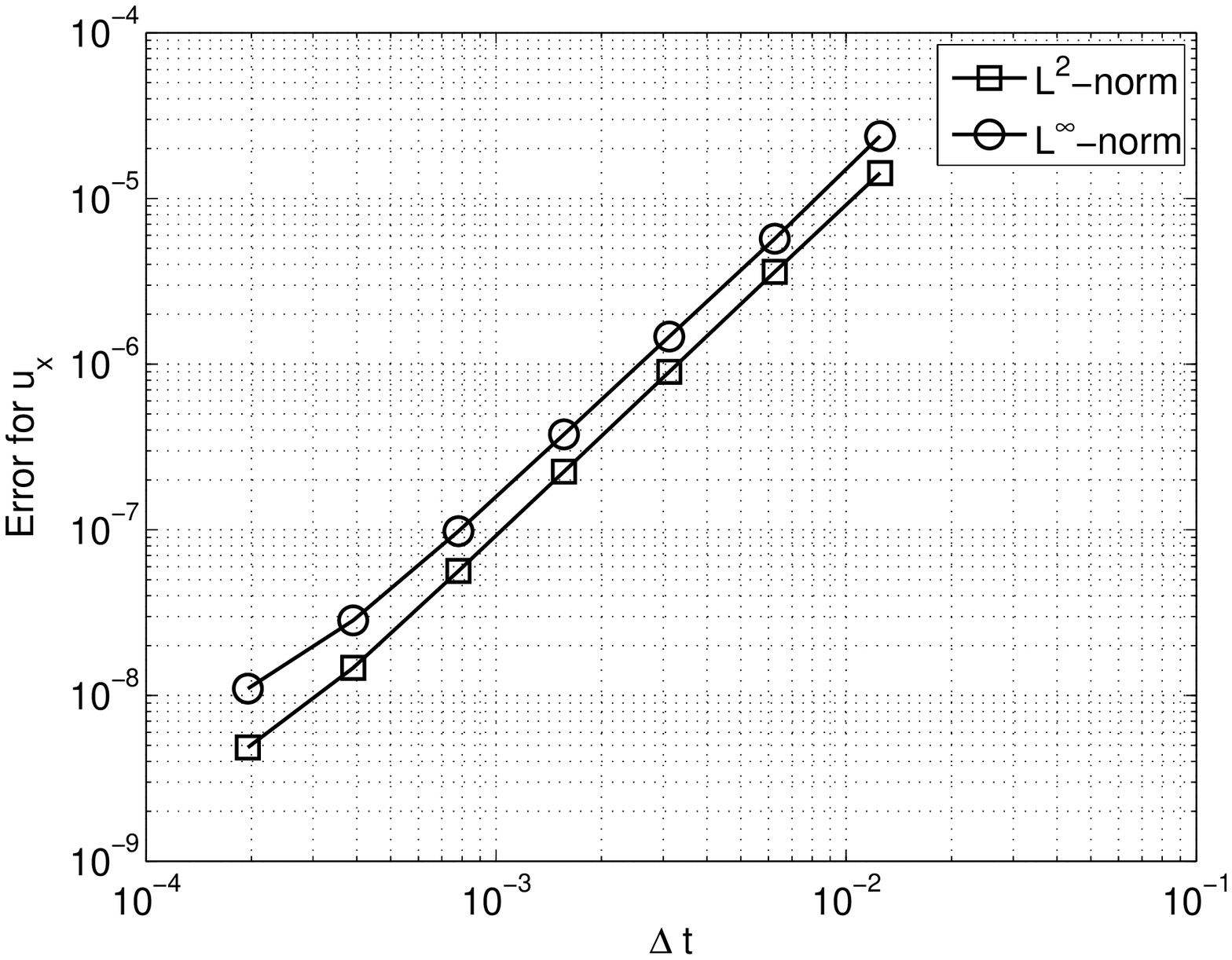}
                \caption{Time convergence.}
        \end{subfigure}
\caption{Spatial convergence for $\Delta t = 3.90 \times 10^{-4}\, s$ and polynomial orders $P=2,3,4,5$ (a); Convergence in time with $P=5$ and decreasing the time increment $\Delta t$ (b).\label{fig:Newmark_poly_conv}}
\end{figure}

{
We now consider spatial convergence using fixed $\Delta t = 3.90 \times 10^{-4}\,s$ and polynomial orders 2 to
6 and results for both cases are shown in Figure \ref{fig:Newmark_poly_conv} for the $L_2$ and $L_\infty$ norms. The $L_2$ error norm dropped about 1000 times when increasing the polynomial order from 2 to 3 and from 3 to 4, indicating exponential spatial convergence. The time convergence is analysed with fixed $P=5$ and decreasing the time increments $\Delta t$. The time rate was quadratic as expected for the Newmark method.}

\subsection{Conrod}
\label{Conrod}

\begin{figure}[H]
	\centering
	\includegraphics[width=0.6\columnwidth]{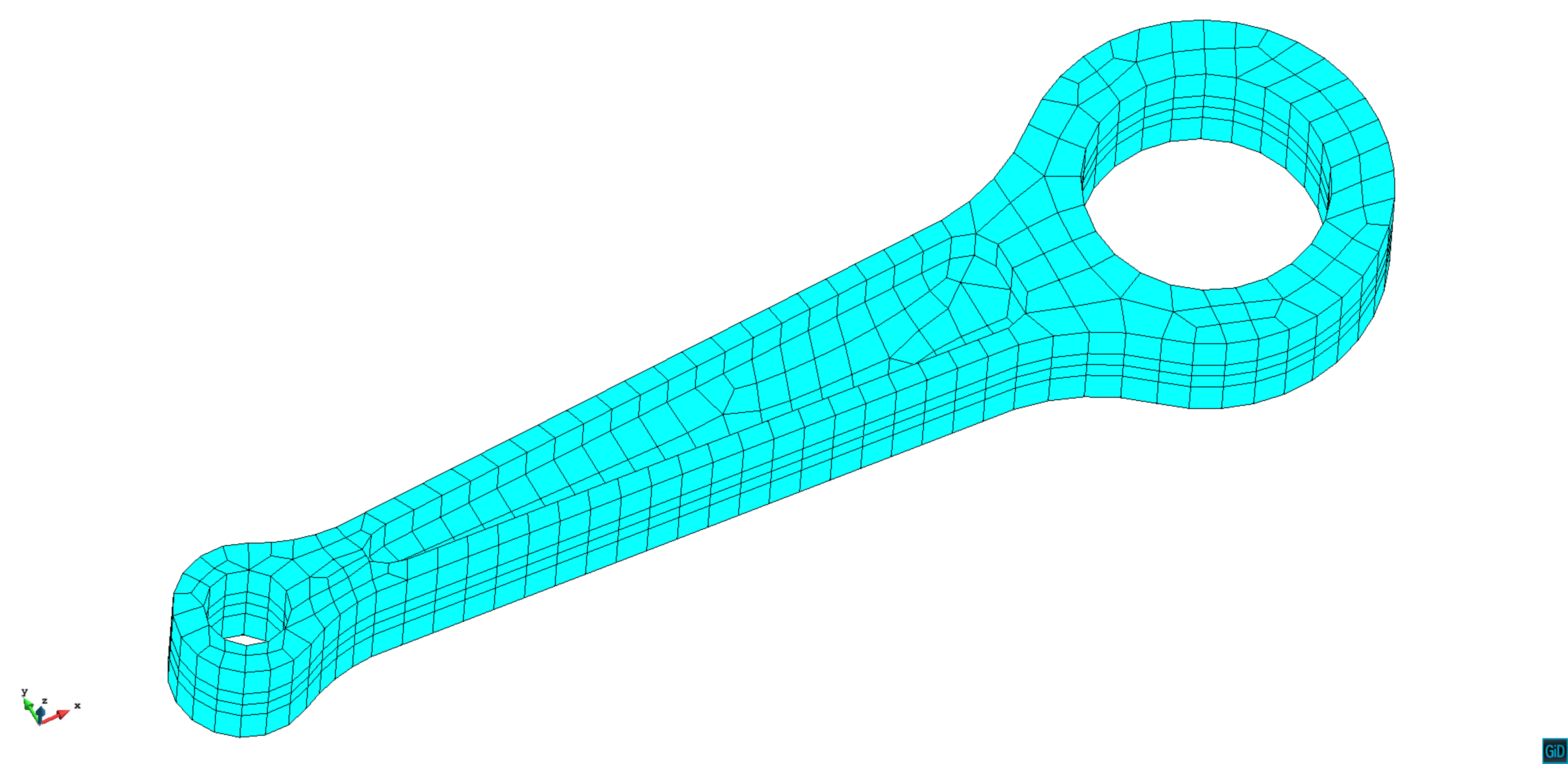}
	\caption{Mesh for the conrod discretized with 768 hexahedra and $P=1$. The material properties used
	for the compressible Neo-Hookean  are $E=500\, Pa$, $\nu = 0.3$,
	$\rho = 9.5 \times 10^{-7}\, Kg/m^3$.\label{fig:conrod}}
\end{figure}

We  consider now the transient dynamic analysis of the conrod illustrated in Fig. \ref{fig:conrod}
and discretized using 768 hexahedra. Homogeneous Dirichlet boundary conditions
($u_x=u_y=u_z=0$) are applied on the nodes of the internal surface of the small end. The conrod is
subjected to a time dependent distributed load on the internal element faces  of the big end in
directions $x$ and $y$. The material model is the compressible Neo-Hookean  with
$E=500\, Pa$, $\nu = 0.3$, $\rho = 9.5 \times 10^{-7}\, Kg/m^3$.  For the SDME bases, we
used $\lambda = 100$ and $k=0.5$.

The total simulation time is $T=5.454 \times 10^{-3}\,s$ which corresponds to one engine cycle for
the rotational speed of $2200\ RPM$. In the explicit analysis, we used the time step
$\Delta t = 6.06 \times 10^{-6}$ for polynomial orders $P=2,\,3$;
$\Delta t = 3.03 \times 10^{-6}$ for $P=4$;
$\Delta t = 1.51 \times 10^{-6}$ for $P=5,\,6,\,7$; and
$\Delta t = 7.58 \times 10^{-7}$ for $P=8$. For the implicit analysis, we considered
$\Delta t = 1.515 \times 10^{-5}$ for all polynomial orders. We solved the linear system of
equations using a parallel, element-by-element diagonal preconditioned conjugate gradient
method (PCG), with tolerances of $10^{-8}$ and $10^{-4}$ respectively for the explicit and
implicit analyses. We also considered the tolerance $10^{-4}$ for the convergence of
the Newton sub-iterations in the implicit case. The initial conditions are
$\mathbf{u}_0=\mathbf{0}\, m$ and $\mathbf{v}_0=\mathbf{0}\, m/s$.

Table \ref{Tab:Conrod_explicit_iterations} presents the average number of PCG iterations per time step for the explicit time integration. We observe that the best results are achieved using the SDME-M basis, which is expected, since in the explicit time integration the operator in the left-hand-side is the mass matrix ({\textit{see Eq.(\ref{Eq:boundary_schur})}}). However, we remark that such behavior can be recovered by the SDME-H basis when we set higher values for $\lambda$ \cite{JCP2015}. The results for the speedup ratios are presented in Table \ref{Tab:Conrod_explicit_time}. Similarly to the number of iterations, we obtained the largest speedups (up to $41.71$) for the SDME-M basis when compared to the ST basis.

\begin{table}[H]
	\centering
	\caption{Average number of iterations per time step for the conrod with explicit time integration using the ST and SDME bases. We observe that the SDME-M basis has the best performance, with a ratio (ST/SDME-M) up to $142.3$. The SDME-H basis increases the number of iterations with the polynomial order with the choice of $\lambda = 100$.
}
	\renewcommand{\arraystretch}{1.2}
	\begin{tabular}{@{}crrcc@{}}
		\toprule
		Order    & Number      & \multicolumn{3}{c}{Average number of iterations per time step} \\ \cline{3-5}
		          & of DOFs     &    ST       & SDME-M & SDME-H   \\ \midrule
		2         & 24\,651     &	283.60    &	22.90  &	30.80  \\
		3         & 75\,984     &	205.90    &	15.70  &	17.20 \\
		4         &	171\,753    &	688.85    &	13.70  &	20.80 \\
		5         &	325\,782    &	403.57    & 11.88  &	26.48 \\
		6         &	551\,895    &	1\,145.08 & 12.22  &	40.65 \\
		7         & 863\,916    &	759.30	  & 11.57  &	54.23 \\
		8         & 1\,275\,669 &   1\,706.09 & 11.99  &	72.55 \\ \bottomrule
	\end{tabular}
	\label{Tab:Conrod_explicit_iterations}
\end{table}

\begin{table}[H]
	\centering
	\caption{Computational speedup per time step for the conrod with explicit time integration. We observe that the highest speedup obtained is $41.71$ for the SDME-M basis using polynomial order $P=6$. \label{Tab:Conrod_explicit_time}}
	\renewcommand{\arraystretch}{1.2}
	\begin{tabular}{@{}crcc@{}}
		\toprule
		Order & Number    & \multicolumn{2}{l}{Speedup $\bar{t}_{ST}/\bar{t}_{SDME}$} \\ \cline{3-4}
		& of DOFs & SDME-M & SDME-H        \\ \midrule
		2 & 24\,651  & 6.48 &	2.43        \\
		3 & 75\,984  & 16.40 & 14.97       \\
		4 & 171\,753 & 25.59 & 17.02        \\
		5 &	325\,782 & 27.13 & 16.42	 \\
		6 &	551\,895 & 41.71 & 16.28	 \\
		7 & 863\,916 & 35.43 & 14.90	 \\
		8 & 1\,275\,669 & 34.34 & 16.28  \\ \bottomrule
	\end{tabular}
\end{table}

Table \ref{Tab:Conrod_implicit_iterations} presents the average number of PCG iterations per time step for the implicit time integration. We observe that the SDME-H basis obtained a smaller number of iterations for $P \le 4$. However, for higher polynomial orders tested, the SDME-M basis provided the best results, with ratio up to $39.20$ when compared to the ST basis. All bases had a smaller number of iterations compared to the explicit case. However, we observe that although we use higher values of $\Delta t$ for the implicit time integration, we need to perform Newton sub-iterations and recalculate the tangent stiffness matrix at every iteration, which is still more time-consuming than the explicit case.

\begin{table}[H]
	\centering
	\caption{Average number of iterations per time step for the conrod with implicit time integration using the ST and SDME bases. We observe that the SDME-M basis has the best performance.\label{Tab:Conrod_implicit_iterations}}
	\renewcommand{\arraystretch}{1.2}
	\begin{tabular}{@{}crccc@{}}
		\toprule
		Order    & Number      & \multicolumn{3}{c}{Average number of iterations per time step} \\ \cline{3-5}
		          & of DOFs     &    ST     & SDME-M & SDME-H   \\ \midrule
		2         & 24\,651     & 122.75   	& 10.17 & 9.75	  \\
		3         & 75\,984     & 90.25	    & 7.25  & 6.92 \\
		4         &	171\,753    & 269.33    & 8.33  & 7.42 \\
		5         &	325\,782    & 211.75    & 8.83  & 9.00	 \\
		6         &	551\,895    & 401.17    & 10.67 & 13.08 \\
		7         & 863\,916    & 383.00	& 12.08 & 16.50 \\
		8         & 1\,275\,669 & 539.08    & 13.75 & 21.08 \\ \bottomrule
	\end{tabular}
\end{table}

{
Increasing the polynomial orders from 2 to 8, the numbers of DOFs increased over 51 times
while the numbers of iterations of the PCG methods for convergence had a slight increase
mainly for the SDME-M basis as can be seen in Tables \ref{Tab:Conrod_explicit_iterations} and \ref{Tab:Conrod_implicit_iterations}. This aspect means that the conditions numbers of the
global matrices had also slight increase as illustrated for the 1D matrices in Figure \ref{ConditionNumbers1D}.
}

\begin{figure}[H]
	\begin{center}$
		\begin{array}{ccc}
		\includegraphics[width=0.3\columnwidth]{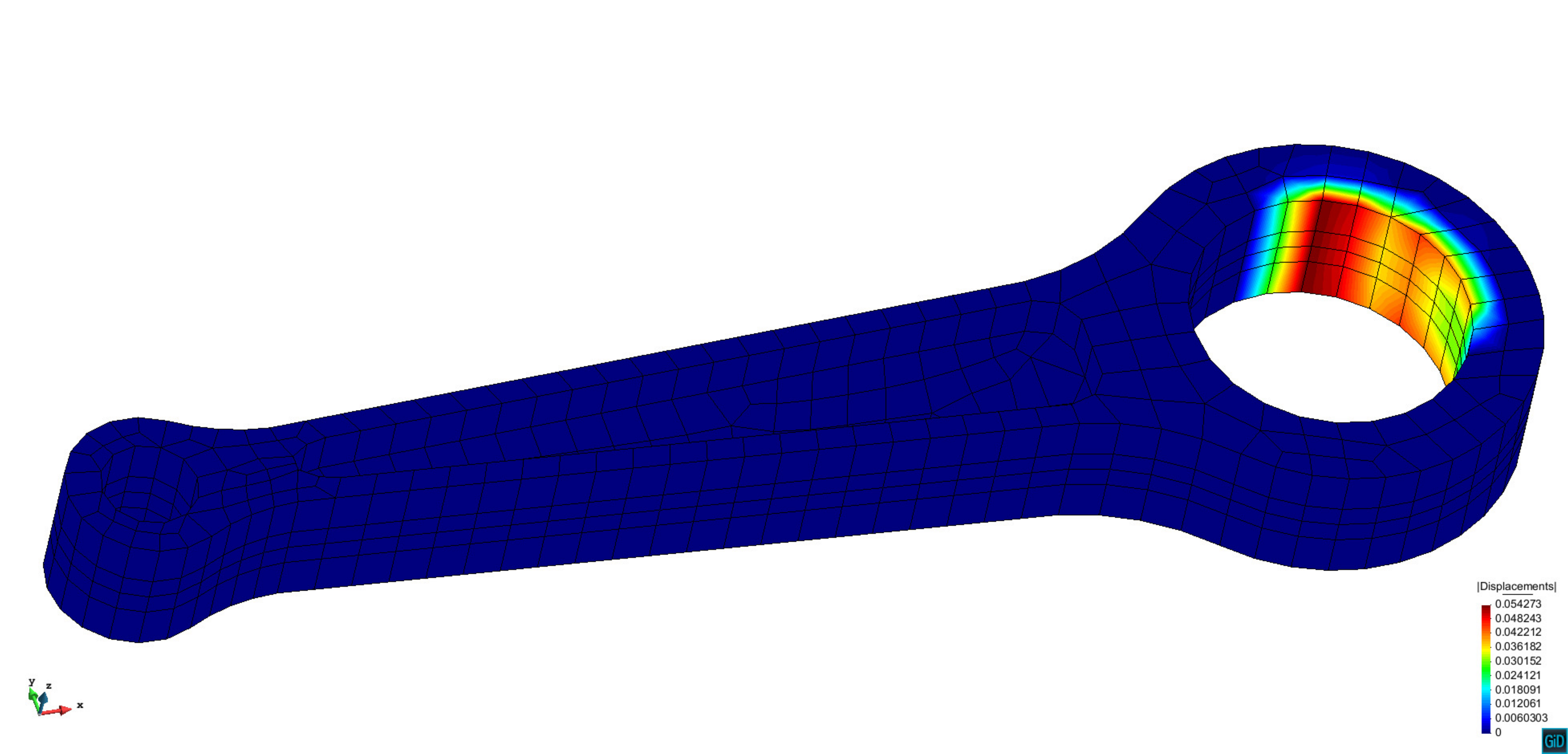} & \includegraphics[width=0.3\columnwidth]{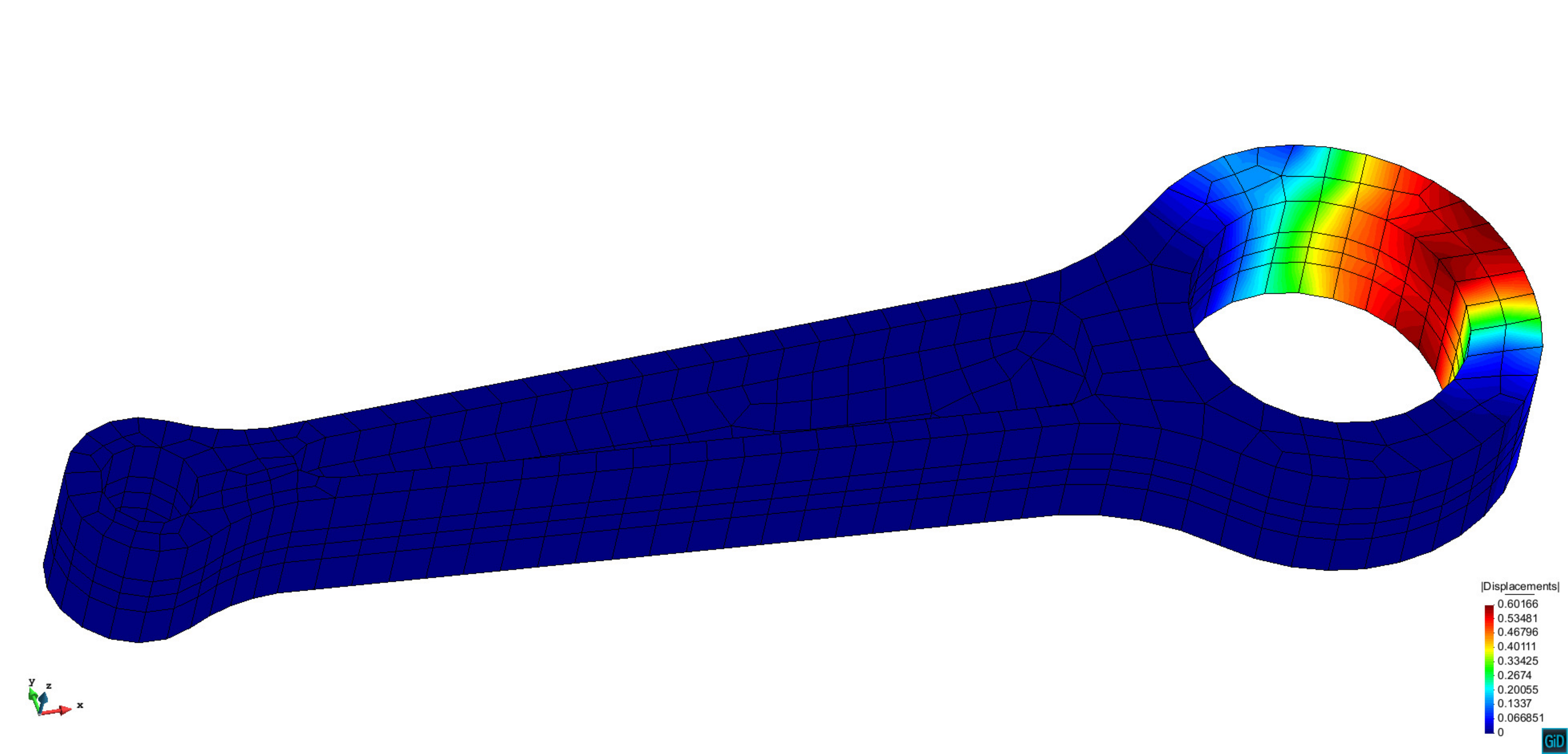} &
		\includegraphics[width=0.3\columnwidth]{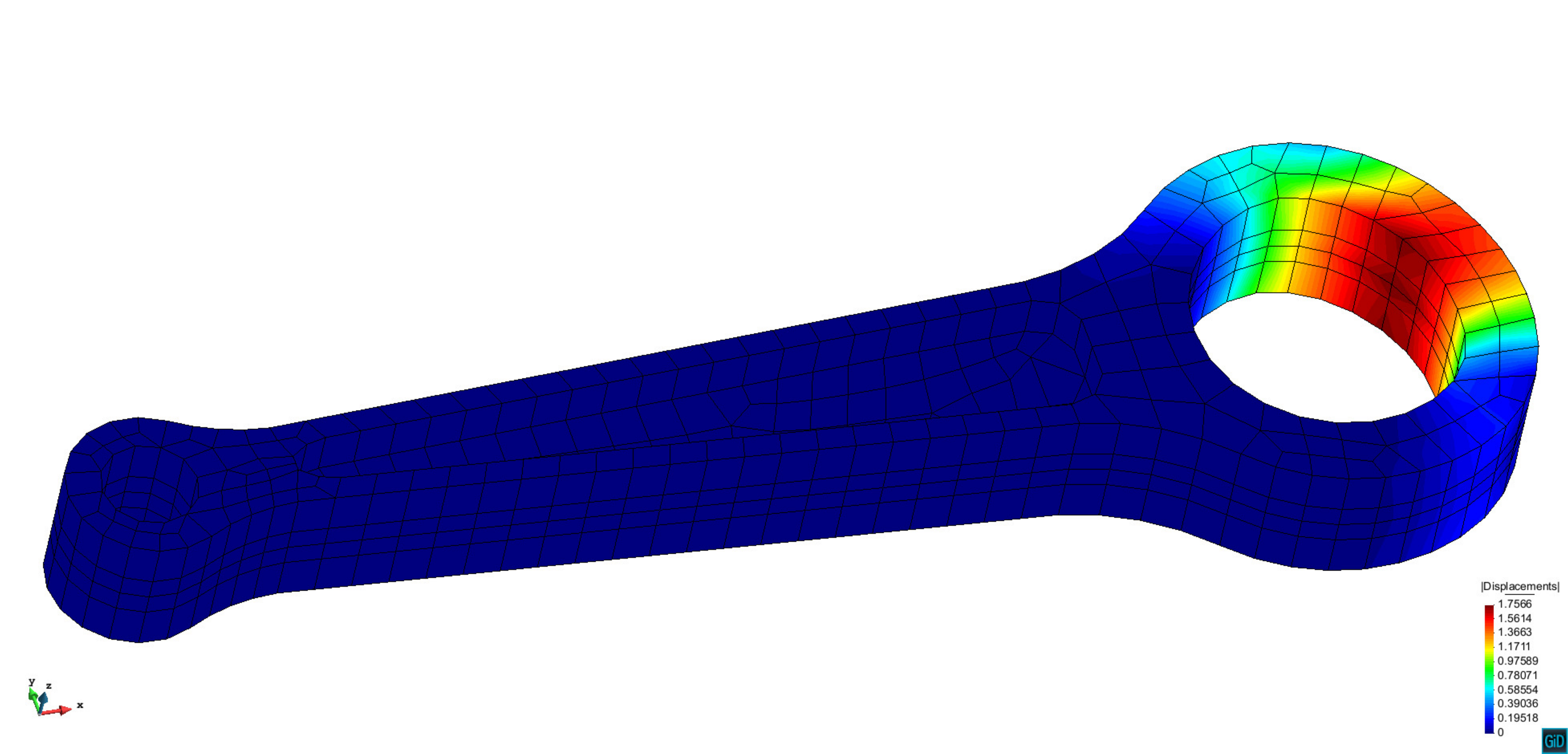} \\
		\includegraphics[width=0.3\columnwidth]{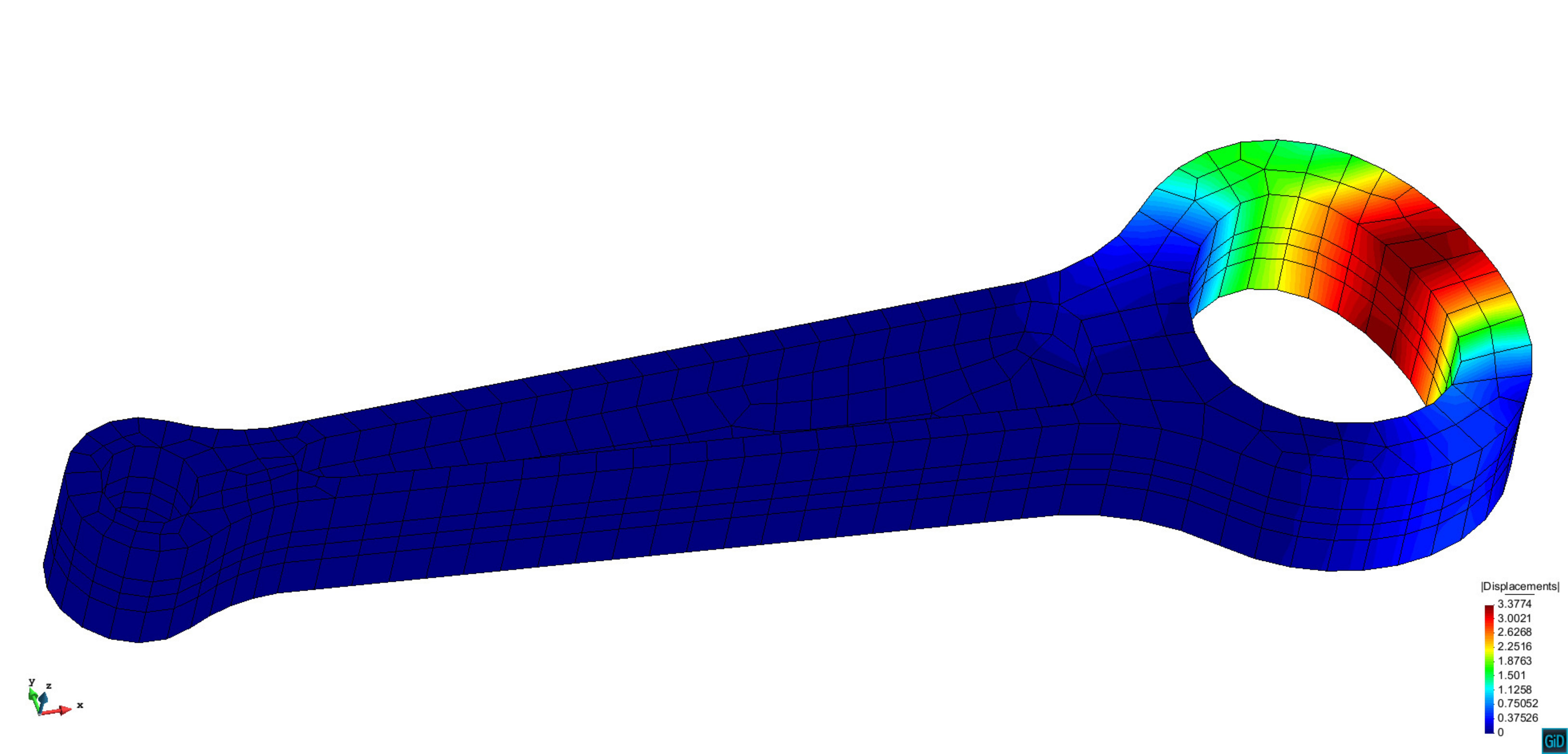} & 	\includegraphics[width=0.3\columnwidth]{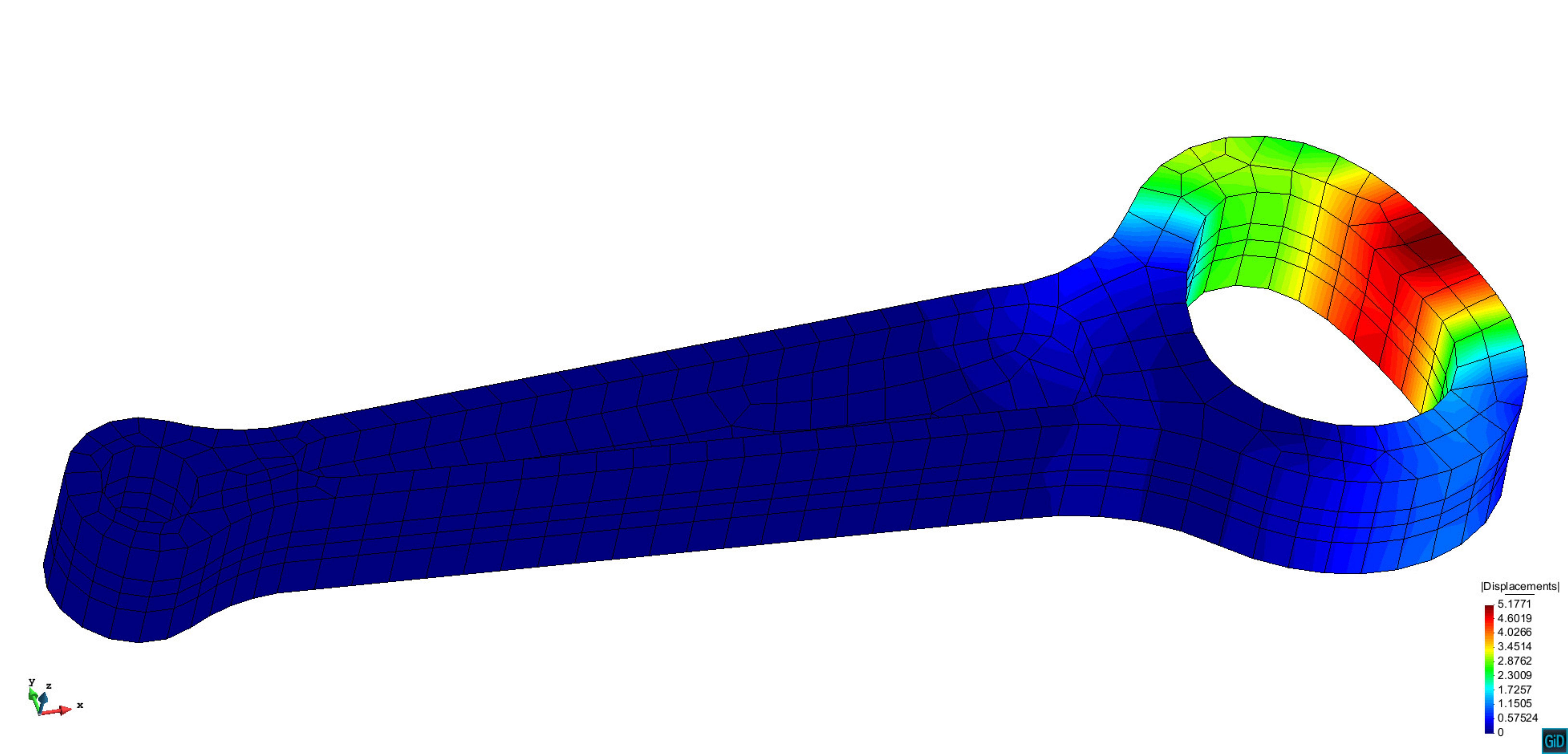} &
		\includegraphics[width=0.3\columnwidth]{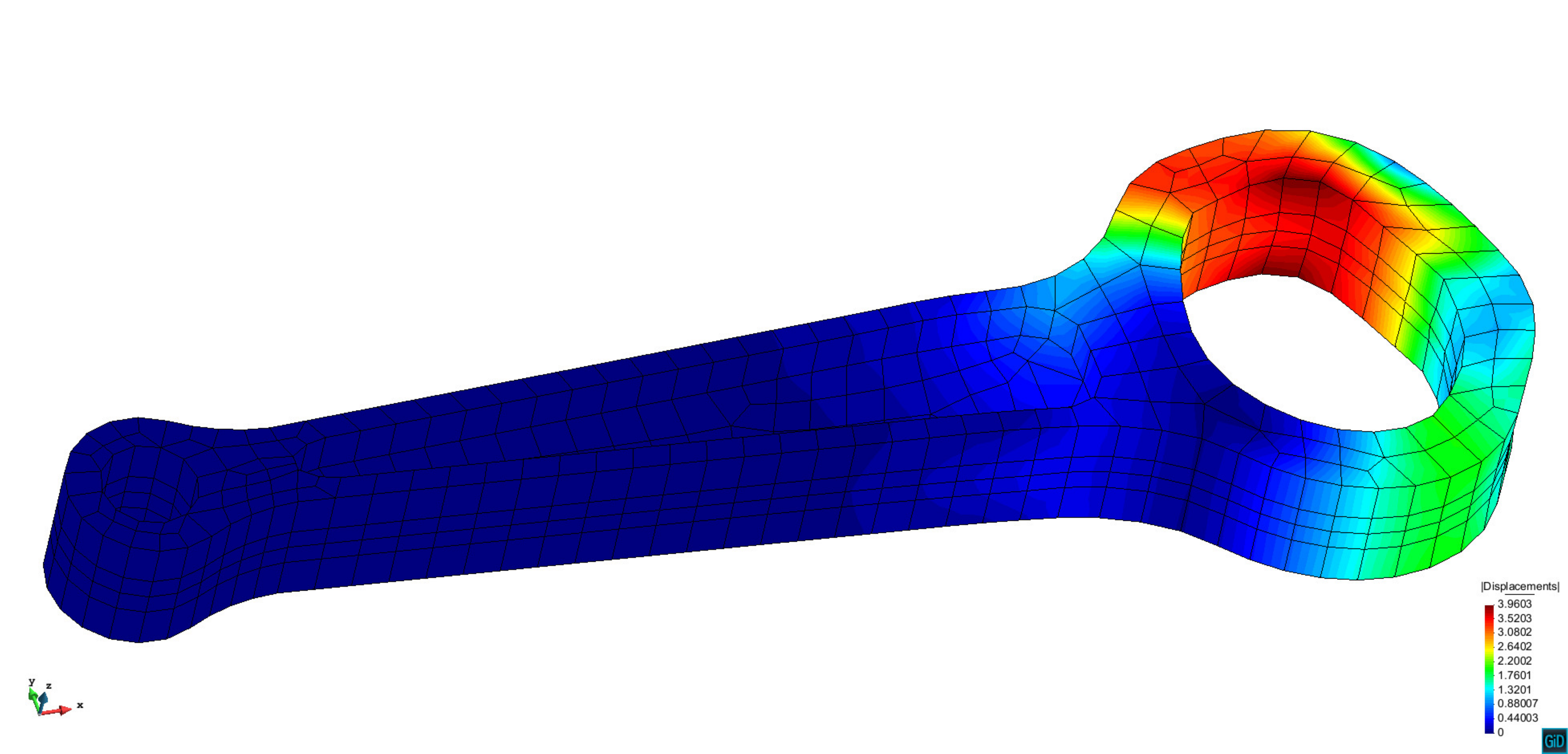} \\
		\includegraphics[width=0.3\columnwidth]{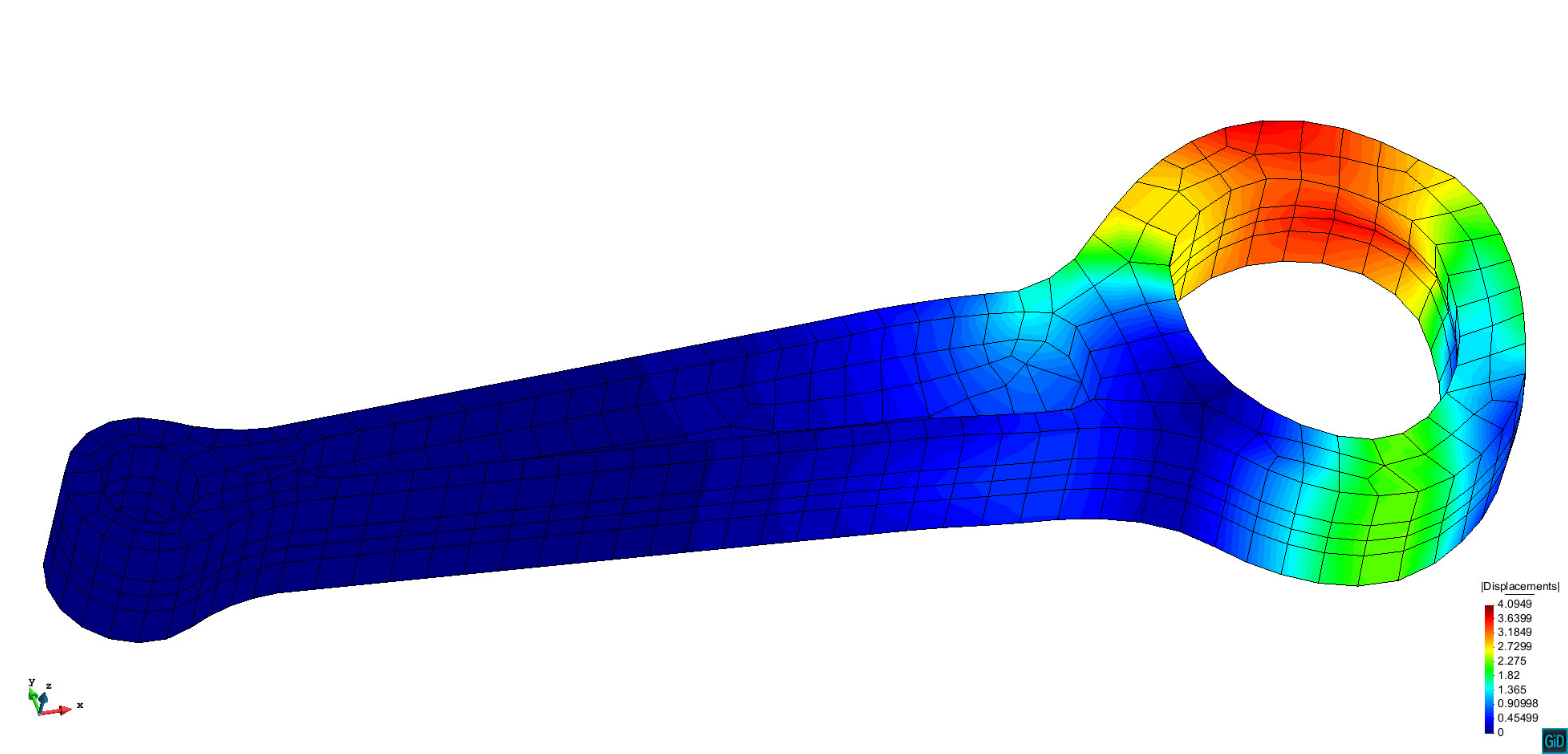} & 	\includegraphics[width=0.3\columnwidth]{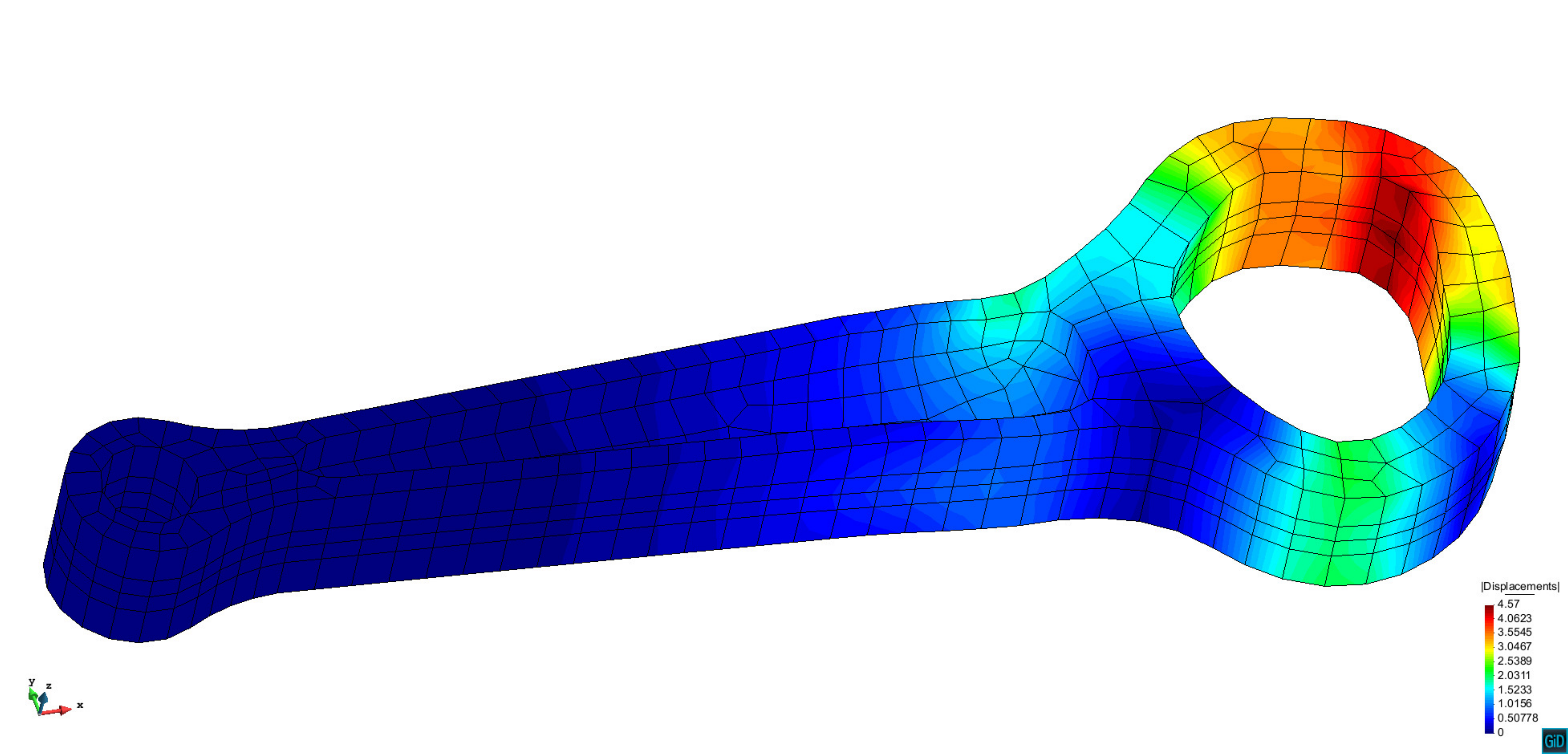} &
		\includegraphics[width=0.3\columnwidth]{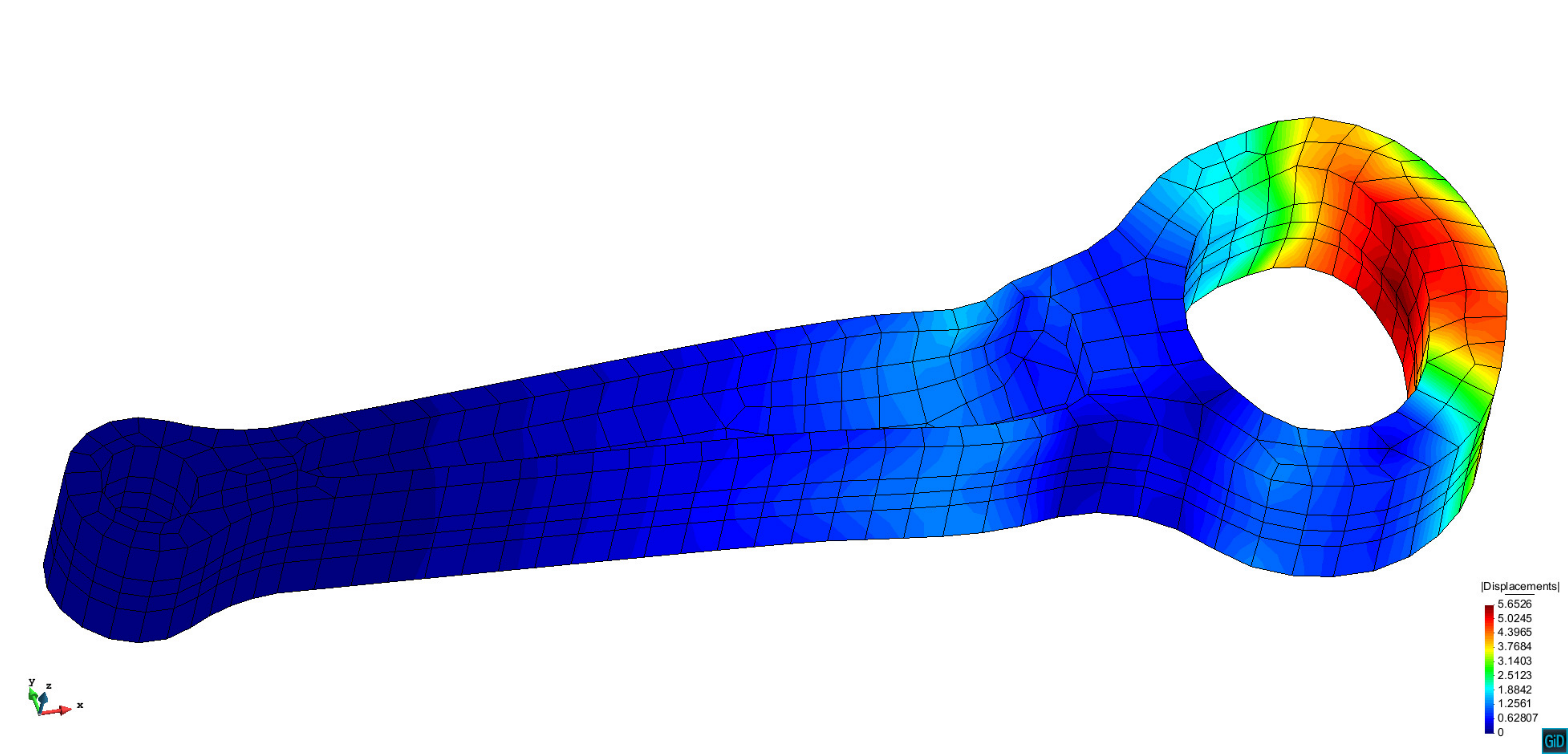} \\
		\end{array}$
	\end{center}
	\caption{Solution for the displacements and deformed geometry of the conrod over time for
	the transient analysis obtained with the SDME-M basis. \label{fig:conrod_time}}
\end{figure}

The solution for the displacements over time is shown in Fig.\ref{fig:conrod_time}. We observe that most of the deformation occurs at the bigger end where the loads are applied.

\subsection{Two-dimensional disk impact problem}\label{2D Impact problem} 

\begin{figure}[!htbp]
\begin{centering}
\includegraphics[scale=0.35]{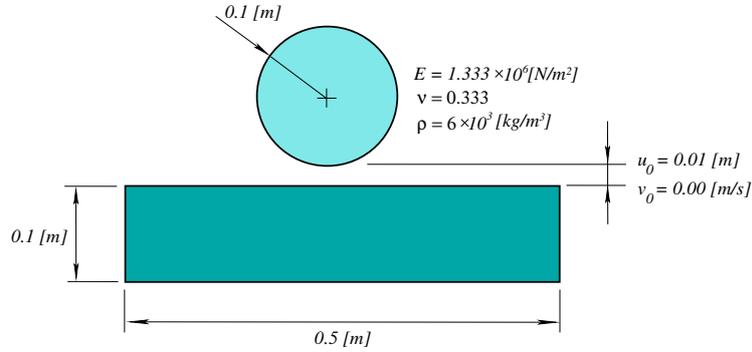}
\par\end{centering}
\protect\caption{Disk impact problem representation with domain dimensions, material properties and
initial conditions.\label{fig:impact_ilustracao}}
\end{figure}

\begin{figure}[!htbp]
\begin{centering}
\includegraphics[scale=0.25]{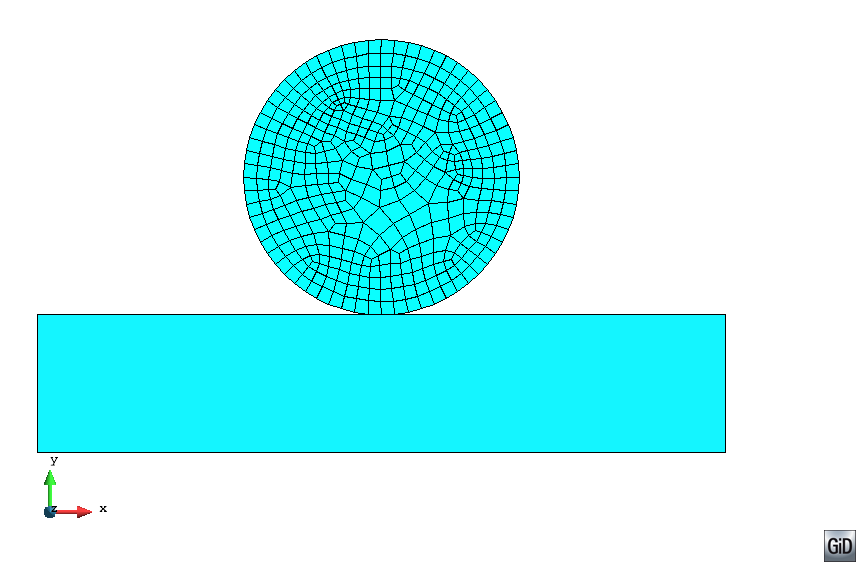}
\par\end{centering}
\protect\caption{Mesh for the disc impact problem. The disk was meshed with quadrangular elements.
The foundation is meshed by one rectangular element and its edges are fully constrained
 \label{fig:impact_mesh}}
\end{figure}

The next example is the small deformation frictionless impact of a linear elastic disk on a foundation as shown
in Fig. \ref{fig:impact_ilustracao}  \citep{Renardmassredis}.
The implicit Newmark time integration scheme was used. The convergence tolerances used for the CGGS and
Newton-Raphson procedures were both $10^{-6}$. The Schur complement was taken for the tangent
stiffness matrix and residue vector. The mesh used is illustrated in Fig. \ref{fig:impact_mesh} and the penalty
parameter is $\epsilon_{N}=1.0\times10^{6}$ and $\Delta\epsilon_{N}=0.0$. The geometry and material properties are
presented in Fig \ref{fig:impact_ilustracao}. We used $P+1$ Gauss-Legendre integration points for the contact elements
that were enough to achieve good results. The tolerances for the gap function and contact stress were $10^{-3}$
and $10^{-2}$, respectively .The integration time is $T=0.1\, s$ and $\Delta t=10^{-4}\, s$. The initial conditions are
$\mathbf{u}_0=\{0.00 \; 0.01\; 0.00\}^T\, m$ and $\mathbf{v}_0=\mathbf{0}\, m/s$. The foundation was discretized with one finite element with
fully constrained edges.

\begin{figure}[!htbp]
\includegraphics[scale=0.2]{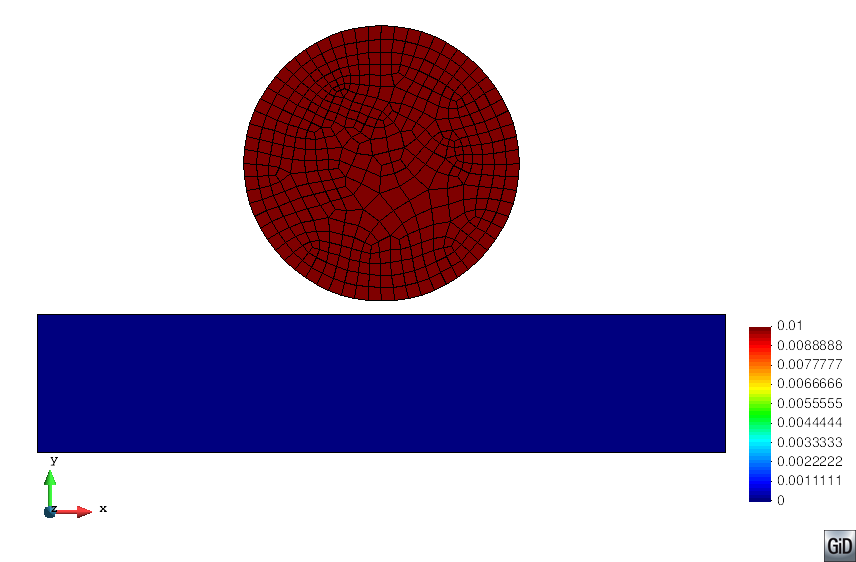}\includegraphics[scale=0.2]{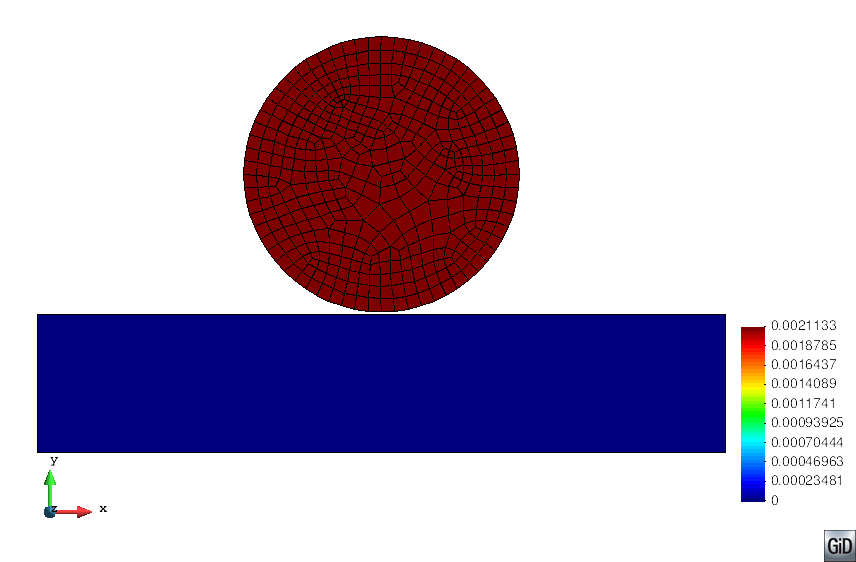}\\\includegraphics[scale=0.2]{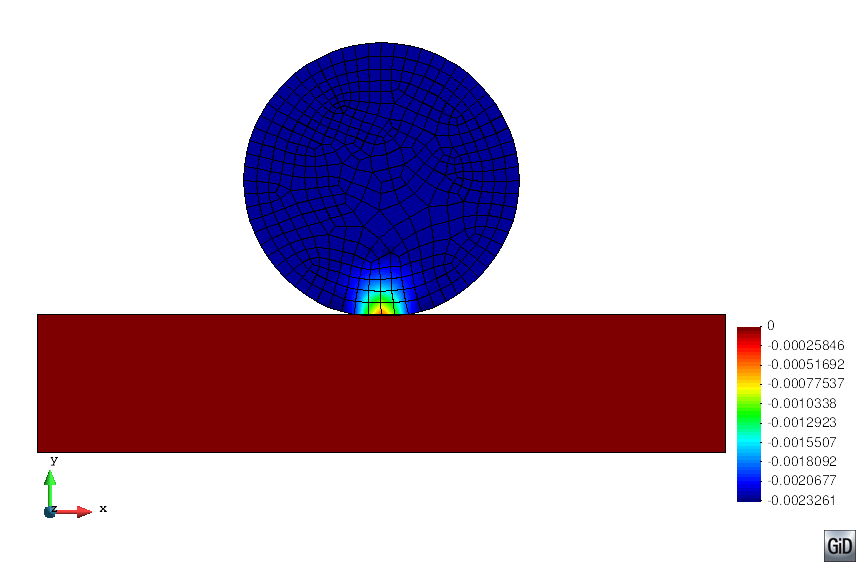}\includegraphics[scale=0.2]{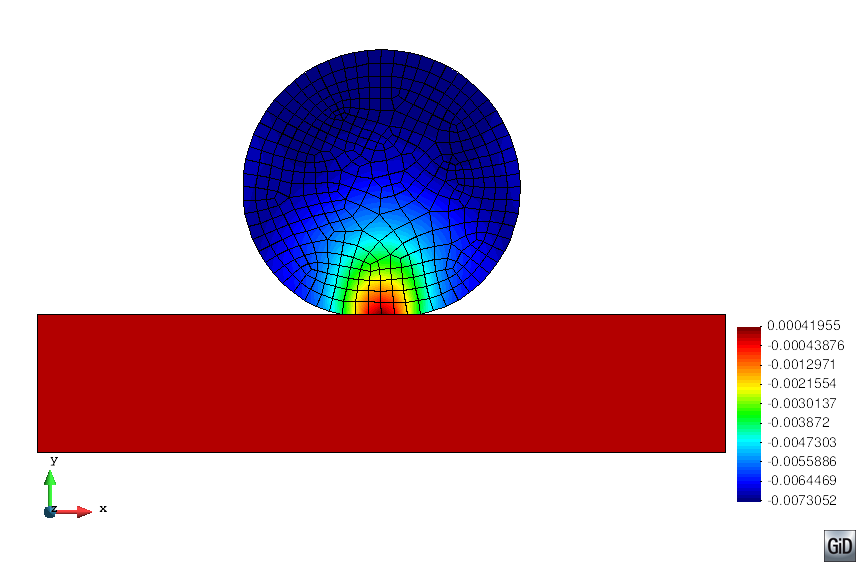}\\\includegraphics[scale=0.2]{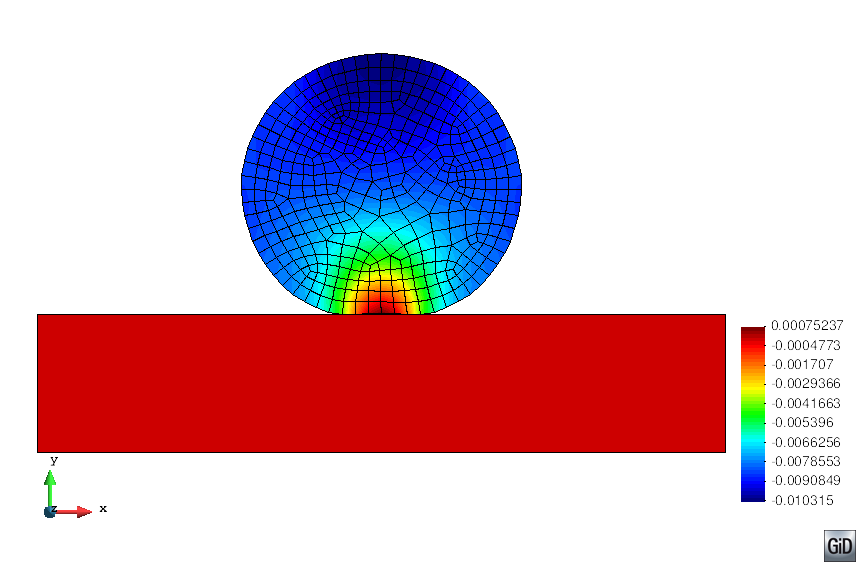}\includegraphics[scale=0.2]{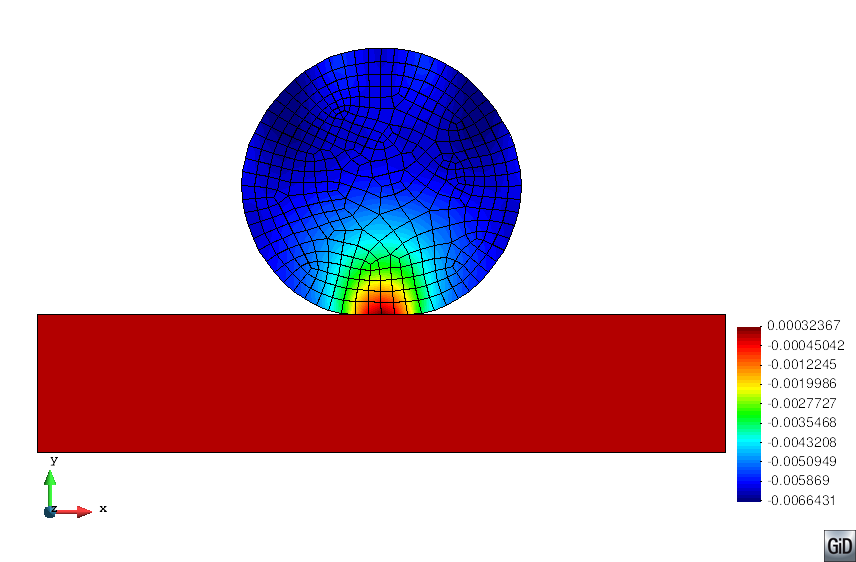}\\
\includegraphics[scale=0.2]{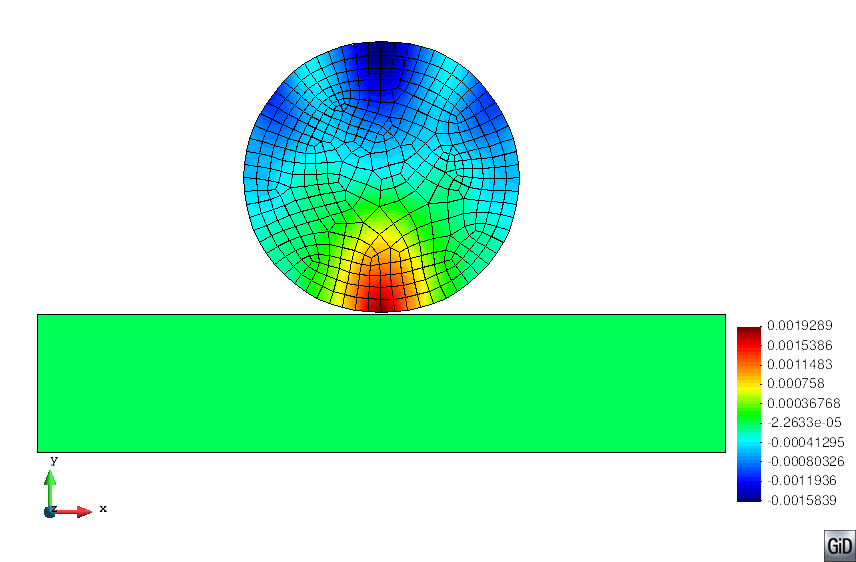}
\includegraphics[scale=0.2]{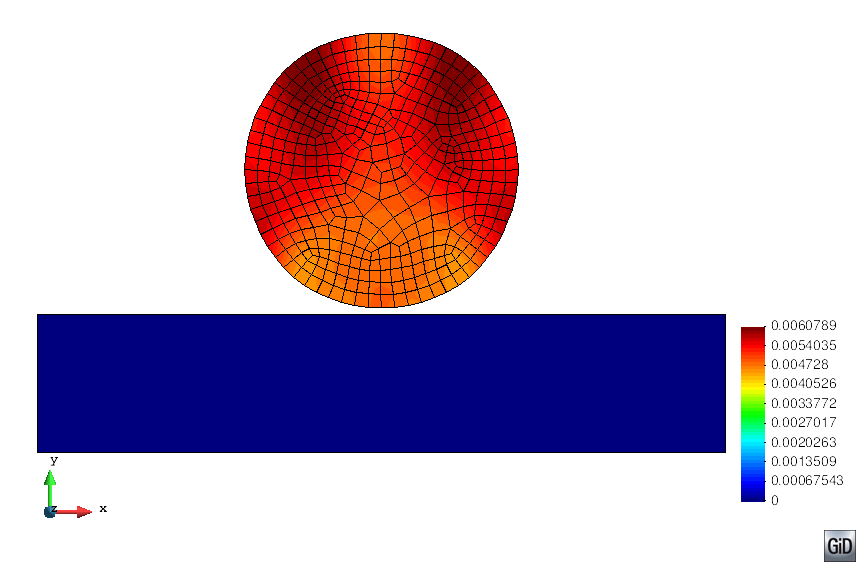}
\protect\caption{Vertical displacement field $u_{y}$ for the deformed geometry at $t=0$, $0.04$, $0.05$, $0.06$,
$0.07$, $0.08$, $0.09$ and $0.1s$, respectively.\label{fig:impact_1_deformation}}
\end{figure}

\begin{figure}[!htbp]
\begin{centering}
\includegraphics[scale=0.47]{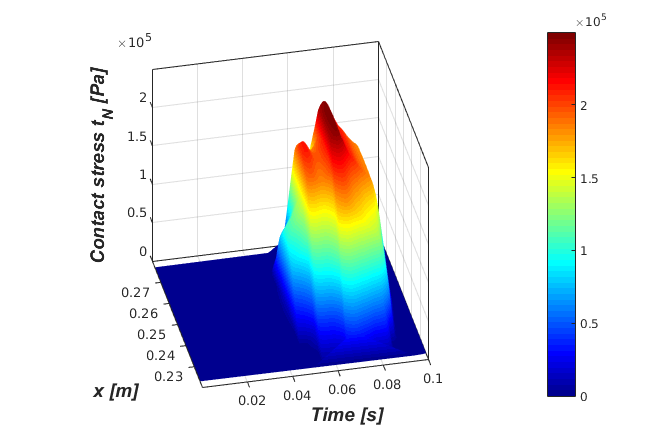}%
\includegraphics[scale=0.47]{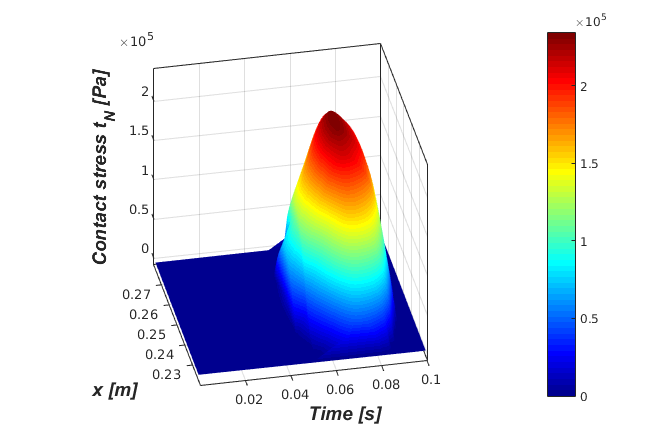}
\par\end{centering}
\protect\caption{Contact stress distributions for $P=1$ and $P=2$. The stresses are much smoother
for $P=2$ when compared with the oscillatory  distribution for $P=1$. \label{fig:impact_stress_p1_p2}}
\end{figure}

Figure \ref{fig:impact_1_deformation} shows the displacement field $u_{y}$ in the deformed geometry for
different time steps of the solution for interpolation order $P=1$.
Figure \ref{fig:impact_stress_p1_p2} shows the comparison of the contact stress field $t_N$ for $P=1$ and $P=2$.
There is more oscillation in the contact stress field with $P=1$. The oscillation is reduced by increasing the interpolation order, as can be seen in Fig. \ref{fig:impact_stress_p1_p2}. The increase of the interpolation order induces a better distribution of the mass inside the finite element, reducing the effect caused by the kinetic energy loss in the contact area (velocity became instantly zero at the contact surface).

%
%

Tables \ref{Tab:CGGS_IMPACT}, \ref{Tab:time_IMPACT} and \ref{Tab:SpeedUp_IMPACT} present the average
numbers of CGGS iterations, average time per time step and speedup, respectively, using the standard
modal Jacobi, the nodal Lagrange and the SDME bases. The best results were achieved for the SDME-H basis.

\begin{table}[H]
\centering
\caption{Average number of CGGS iterations per time step. We observe that the SDME-H basis has the best performance}
\label{Tab:CGGS_IMPACT}
\begin{tabular}{@{}cccccc@{}}
\toprule
\multirow{2}{*}{Order} & \multirow{2}{*}{\begin{tabular}[c]{@{}c@{}}Number \\ of DOFs\end{tabular}} & \multicolumn{3}{c}{Average number of CGGS iterations} \\ \cmidrule(l){3-6}
                        &                                                                            & ST      & Lagrange       & SDME-M            & SDME-H  \\ \cmidrule(r){1-6}
1 & 962 & 27.74 & 19.68 & - & -\tabularnewline
2 & 3722 & 75.05 & 20.17  & 19.53  & 19.52 \tabularnewline
3 & 8282 & 64.47 & 19.04  & 21.57  & 20.20\tabularnewline
4 & 14642 & 131.03 & 22.24  & 25.60 & 22.32 \tabularnewline
5 & 22802 & 129.05 & 25.38  & 27.61  & 22.52 \tabularnewline
6 & 32762 & 185.15 & 28.35  & 30.32  & 22.18 \tabularnewline
7 & 44522 & 199.64 & 30.94  & 31.07  & 22.99 \tabularnewline
8 & 58082 & 228.65 & 32.99  & 33.63  & 25.55 \tabularnewline
9 & 73442 & 270.68 & 35.95  & 34.88  & 28.10 \tabularnewline
10 & 90602 & 279.37 & 38.47  & 36.97  & 30.50 \tabularnewline
\bottomrule
\end{tabular}
\end{table}

\begin{table}[H]
\centering
\caption{Average time per time step using the CGGS method for the ST, Lagrange and SDME bases.
We observe that the SDME-H basis has the best performance}
\label{Tab:time_IMPACT}
\begin{tabular}{@{}cccccc@{}}
\toprule
\multirow{2}{*}{Order} & \multirow{2}{*}{\begin{tabular}[c]{@{}c@{}}Number \\ of DOFs\end{tabular}} & \multicolumn{3}{c}{Average time for CGGS solution [s]} \\ \cmidrule(l){3-6}
                        &                                                                            & ST      & Lagrange       & SDME-M            & SDME-H  \\ \cmidrule(r){1-6}
1 & 962 & 0.0073 & 0.0057 & - & -\tabularnewline
2 & 3722 & 0.0326 & 0.0112 & 0.0104 & 0.0104\tabularnewline
3 & 8282 & 0.0623 & 0.0224 & 0.0247 & 0.0236\tabularnewline
4 & 14642 & 0.2013 & 0.0436 & 0.0483 & 0.0437\tabularnewline
5 & 22802 & 0.3034 & 0.0733 & 0.0783 & 0.0671\tabularnewline
6 & 32762 & 0.6009 & 0.1127 & 0.1180 & 0.0933\tabularnewline
7 & 44522 & 0.8590 & 0.1610 & 0.1610 & 0.1283\tabularnewline
8 & 58082  & 1.2552 & 0.2156 & 0.2211 & 0.1784\tabularnewline
9 & 73442 & 1.8583 & 0.2927 & 0.2866 & 0.2429\tabularnewline
10 & 90602 & 2.3456 & 0.3850 & 0.3698 & 0.3183\tabularnewline
\bottomrule
\end{tabular}
\end{table}

\begin{table}[H]
\centering
\caption{Speedup ratio for the computation time to solve the linear system of equations between the standard Jacobi basis and Lagrange basis and between the standard Jacobi basis and the minimum energy bases.}
\label{Tab:SpeedUp_IMPACT}
\begin{tabular}{@{}ccccc@{}}
\toprule
\multirow{2}{*}{Order} & \multirow{2}{*}{\begin{tabular}[c]{@{}c@{}}Number \\ of DOFs\end{tabular}} & \multicolumn{3}{c}{Speedup} \\ \cmidrule(l){3-5}
                        &                                                                           & Lagrange       & SDME-M            & SDME-H  \\ \cmidrule(r){1-5}
1 & 962 & 1.28 & - & -\tabularnewline
2 & 3722 & 2.91 & 3.15 & 3.13\tabularnewline
3 & 8282 & 2.78 & 2.52 & 2.64\tabularnewline
4 & 14642 & 4.62 & 4.17 & 4.60\tabularnewline
5 & 22802 & 4.14 & 3.87 & 4.52\tabularnewline
6 & 32762 & 5.33 & 5.09 & 6.44\tabularnewline
7 & 44522 & 5.33 & 5.33 & 6.69\tabularnewline
8 & 58082  & 5.82 & 5.68 & 7.04\tabularnewline
9 & 73442 & 6.35 & 6.48 & 7.65\tabularnewline
10 & 90602 & 6.09 & 6.34 & 7.37\tabularnewline
\bottomrule
\end{tabular}
\end{table}

\subsection{Three-dimensional cylinder impact problem}
\label{3D Impact problem}

\begin{figure}[!htbp]
\begin{centering}
\includegraphics[scale=0.38]{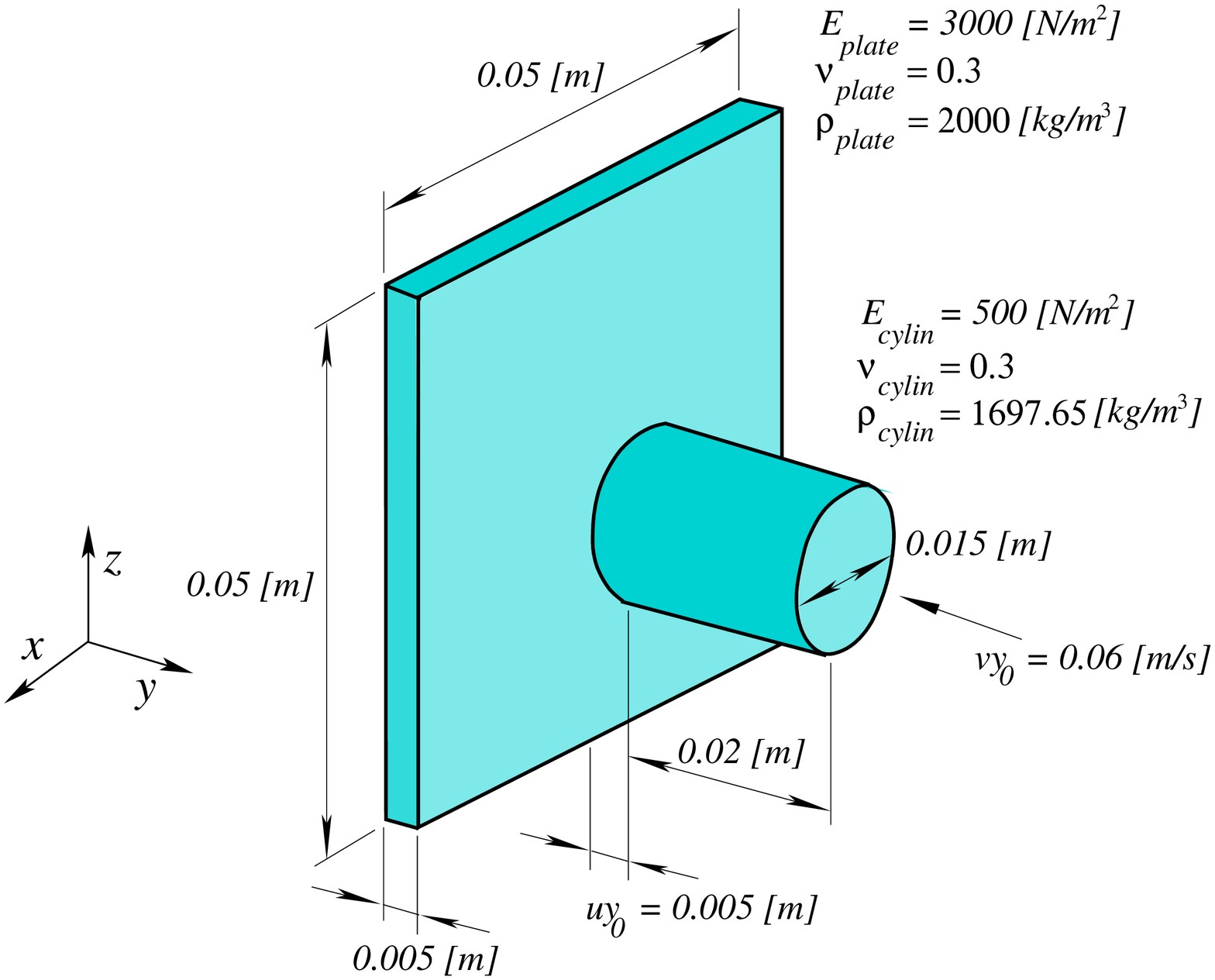}\includegraphics[scale=0.25]{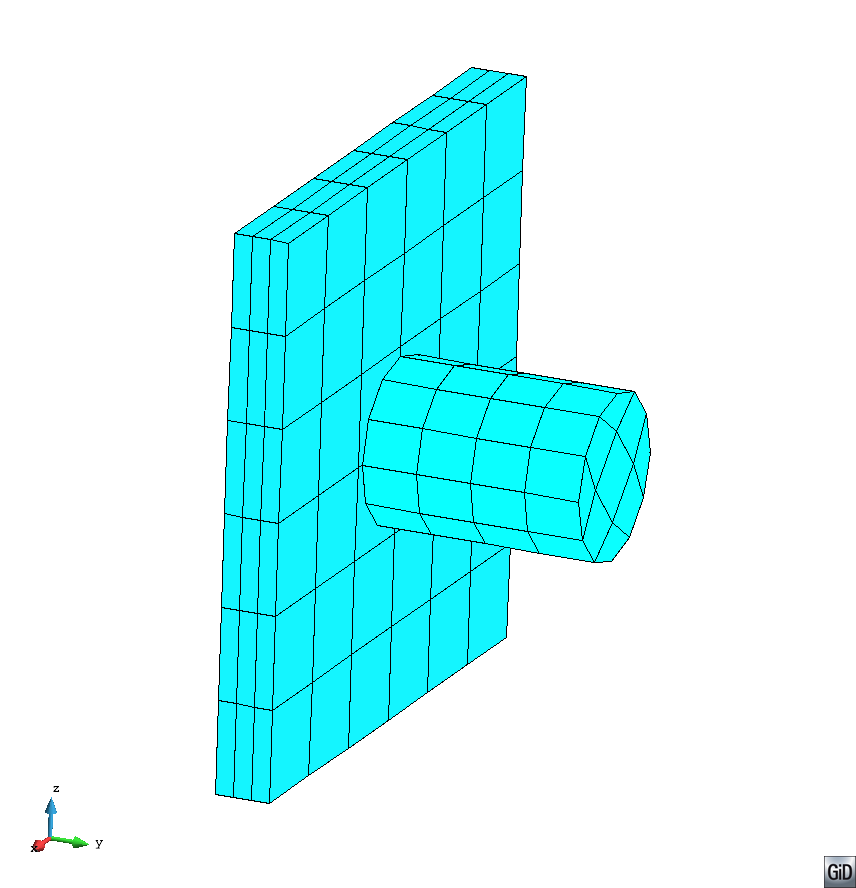}
\par\end{centering}
\protect\caption{Domain geometry, material properties and  initial conditions andesh of
hexahedra for the cylinder impact problem.\label{fig:impact_ilustracao_3D}}
\end{figure}

We present now the results for a frictionless impact of a hyperelastic cylinder on a plate, Fig.
\ref{fig:impact_ilustracao_3D}. The implicit Newmark time integration scheme was used.
The convergence tolerances for the CGGS and Newton-Raphson procedures were both $10^{-6}$.
The Schur complement was taken for the tangent stiffness matrix and the residue vector.
The penalty parameter was $\epsilon_{N}=1.0\times10^{4}$ and $\Delta\epsilon_{N}=1.0\times10^{3}$.
Large deformation was considered and the geometry and material properties are presented in Fig
\ref{fig:impact_ilustracao_3D}. We used $P+1$ Gauss-Legendre integration points for the contact elements.
The tolerances for the gap function and contact stress were $10^{-3}$ and $10^{-2}$, respectively.
The integration time is $T=0.2\, s$ and $\Delta t=10^{-3}\, s$. The initial conditions are
$\mathbf{u}_0=\{0.000\; 0.005\; 0.000\}^T\, m$ and $\mathbf{v}_0=\{0.000\; -0.060\; 0.000\}^T\, m/s$.
The two faces of the plate in the $xy$-plane were completely fixed.

\begin{figure}[!htbp]
\begin{centering}
\includegraphics[scale=0.22]{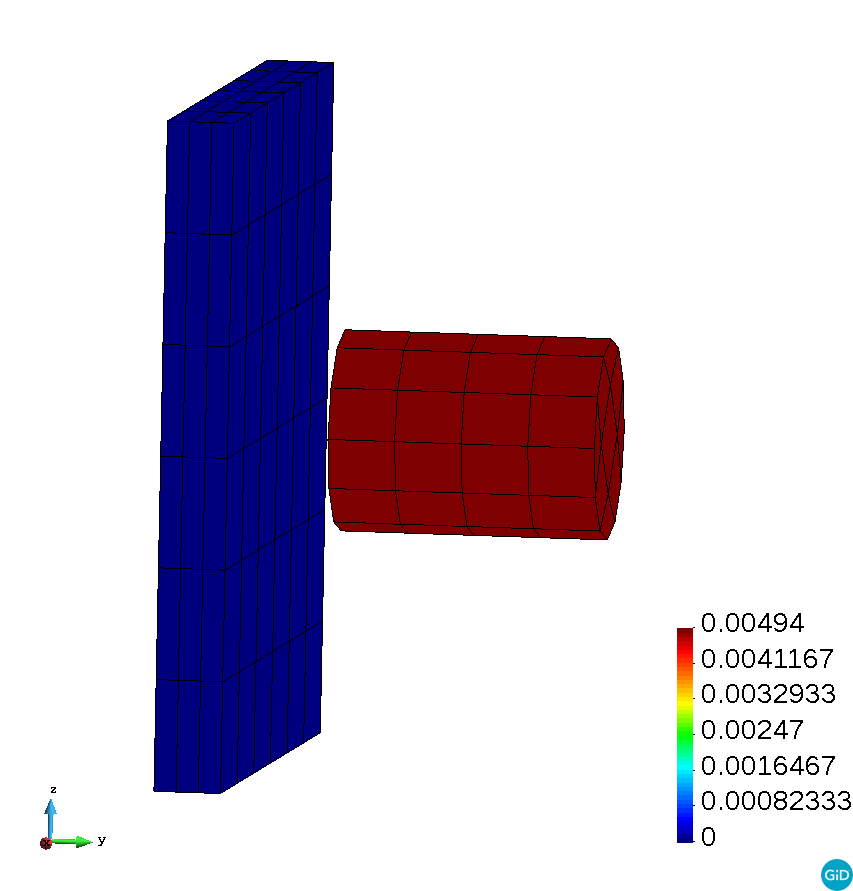}\includegraphics[scale=0.22]{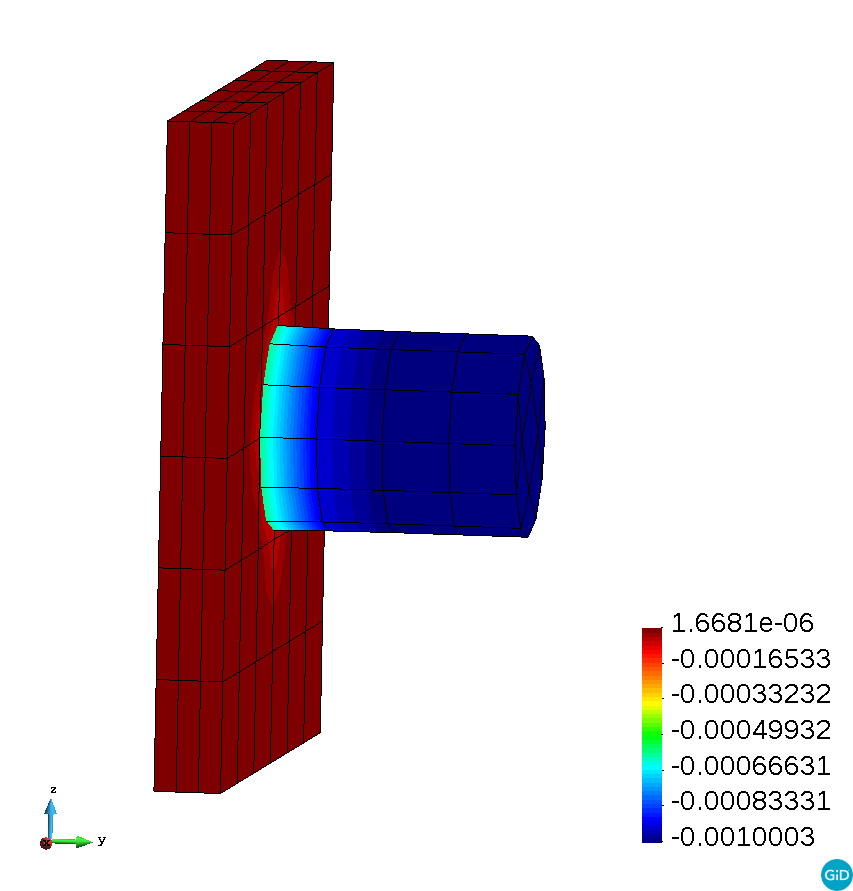}\\\includegraphics[scale=0.22]{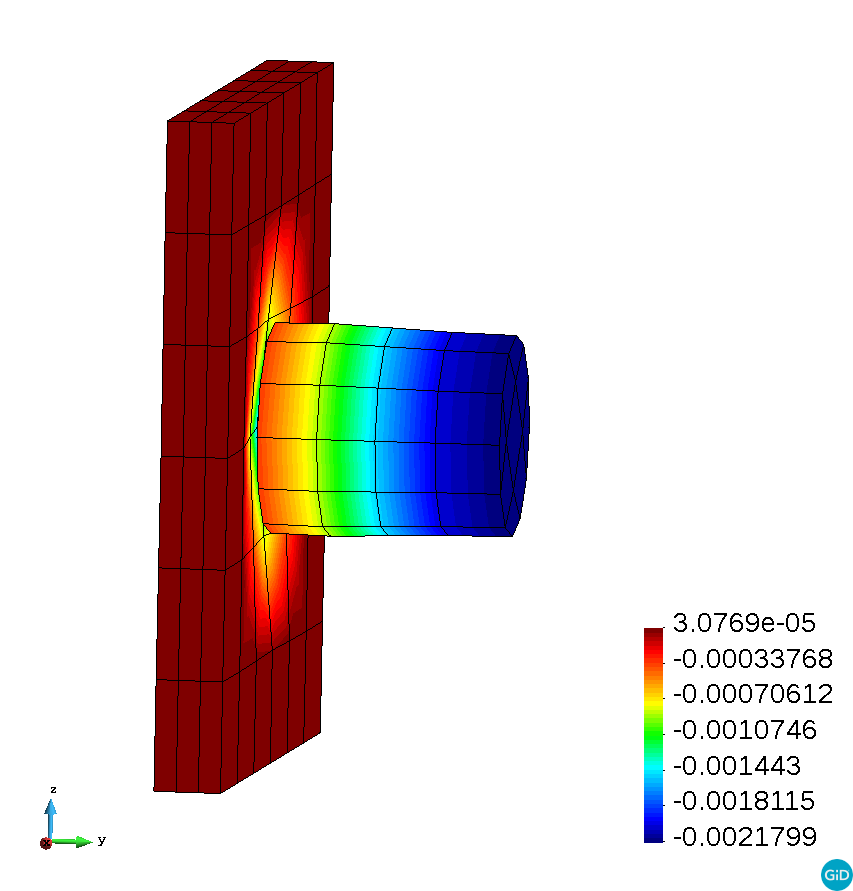}\includegraphics[scale=0.22]{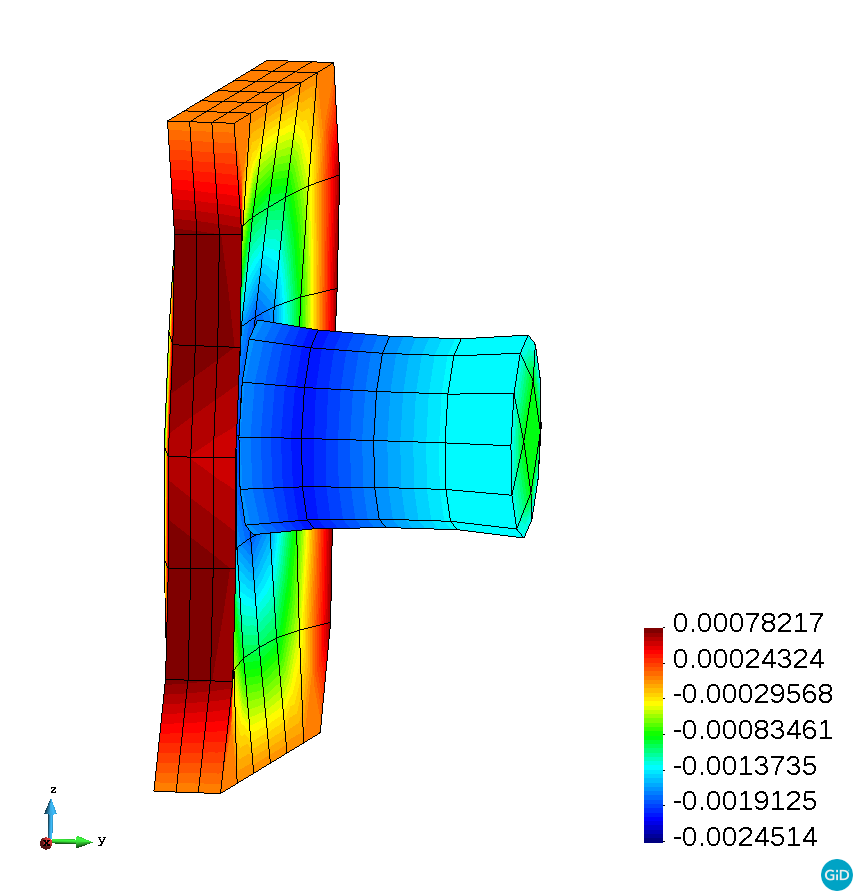}\\\includegraphics[scale=0.22]{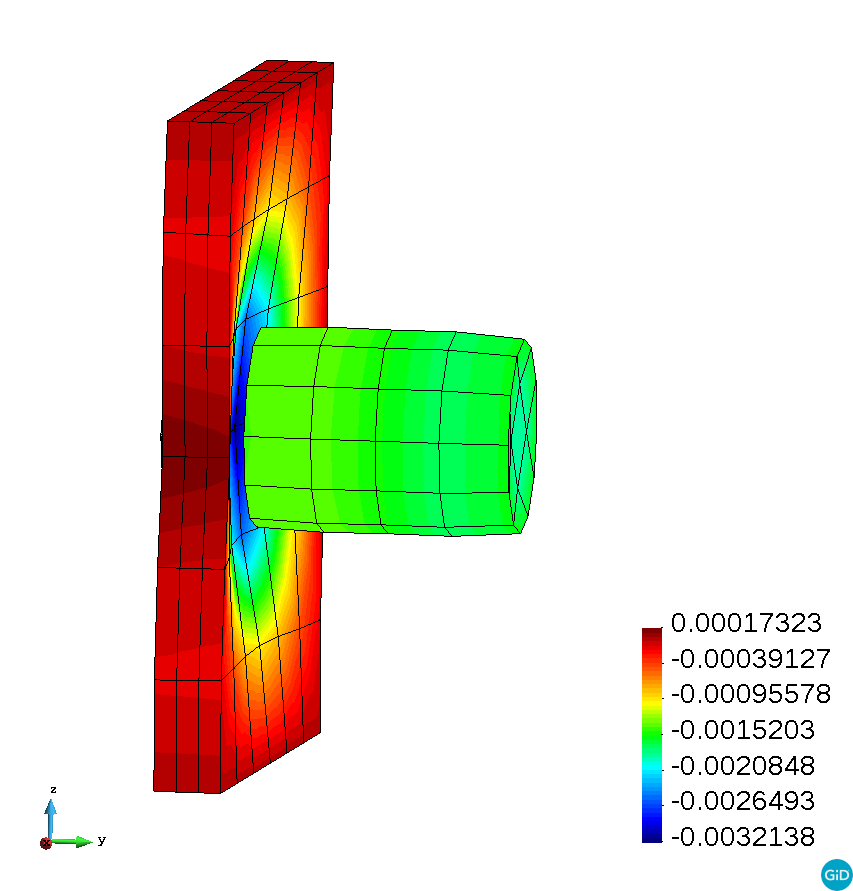}\includegraphics[scale=0.22]{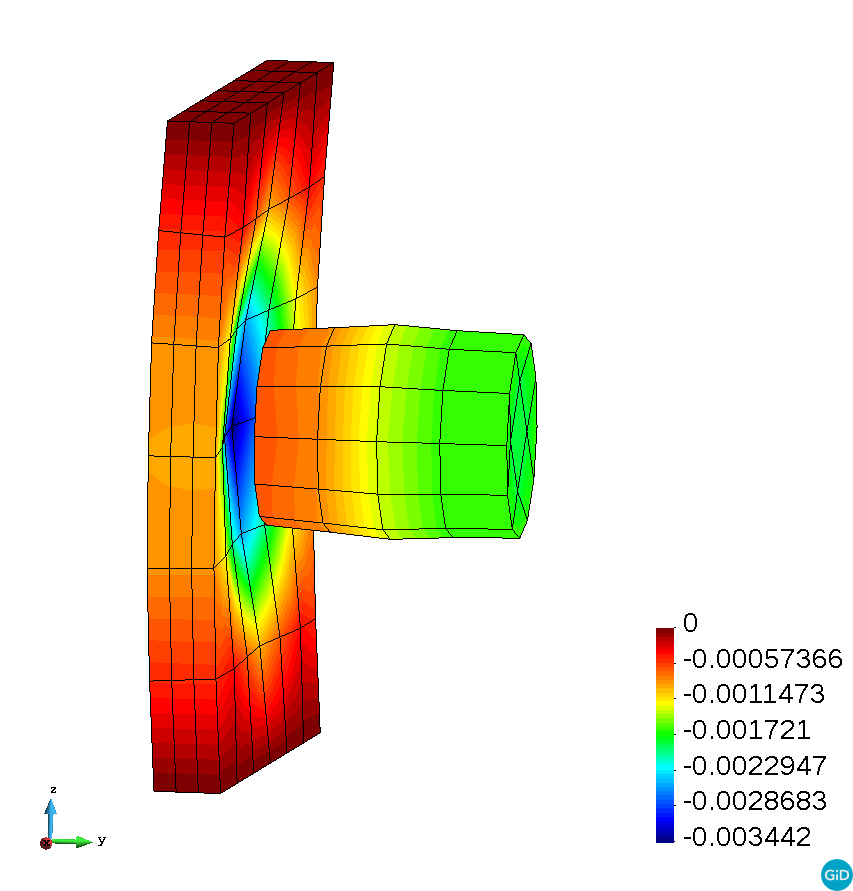}
\par\end{centering}
\protect\caption{Displacement field $u_{y}$ for the deformed structure in $t=0$, $0.10$, $0.12$, $0.16$, $0.18$ and $0.20s$, respectively, for the cylinder impact problem.\label{fig:impact_deformation_3D}}
\end{figure}

Figure \ref{fig:impact_deformation_3D} shows the displacement field $u_{y}$ in the deformed structure at different time steps of the solution for interpolation order $P=3$.
%
%
Tables \ref{Tab:CGGS_IMPACT_3D}, \ref{Tab:time_IMPACT_3D} and \ref{Tab:SpeedUp_IMPACT_3D} present the average
numbers of CGGS iterations, average time per time step and speedup, respectively, using the standard
modal Jacobi, the nodal Lagrange and SDME bases. The best results were achieved for the SDME-H basis.
%
\begin{table}[H]
\centering
\caption{Average number of CGGS iterations per time step. The SDME-H basis performed better
for all polynomial orders.}
\label{Tab:CGGS_IMPACT_3D}
\begin{tabular}{@{}cccccc@{}}
\toprule
\multirow{2}{*}{Order} & \multirow{2}{*}{\begin{tabular}[c]{@{}c@{}}Number \\ of DOFs\end{tabular}} & \multicolumn{3}{c}{Average number of CGGS iterations} \\ \cmidrule(l){3-6}
                        &                                                                            & ST      & Lagrange       & SDME-M            & SDME-H  \\ \cmidrule(r){1-6}
1 & 652 & 27.55 & 27.64 & - & -\tabularnewline
2 & 4318 & 220.62 & 39.47 & 32.49 & 32.31\tabularnewline
3 & 14299  & 169.56 & 37.67 & 35.15 & 31.78\tabularnewline
4 & 31753 & 777.90 & 43.17 & 57.43 & 42.93\tabularnewline
5 & 60922 & 785.29 & 50.28 & 52.15 & 41.32\tabularnewline
\bottomrule
\end{tabular}
\end{table}
\begin{table}[H]
\centering
\caption{Average time per time step using the CGGS method for the ST, Lagrange and SDME bases.
The SDME-H basis performed better for all polynomial orders.}
\label{Tab:time_IMPACT_3D}
\begin{tabular}{@{}cccccc@{}}
\toprule
\multirow{2}{*}{Order} & \multirow{2}{*}{\begin{tabular}[c]{@{}c@{}}Number \\ of DOFs\end{tabular}} & \multicolumn{3}{c}{Average time for CGGS solution [s]} \\ \cmidrule(l){3-6}
                        &                                                                            & ST      & Lagrange       & SDME-M            & SDME-H  \\ \cmidrule(r){1-6}
1 & 652 & 0.0092 & 0.0094 & - & -\tabularnewline
2 & 4318 & 0.4651 & 0.1085 & 0.0870 & 0.0866\tabularnewline
3 & 14299 & 1.8288 & 0.4987 & 0.4446 & 0.4274\tabularnewline
4 & 31753 & 24.7392 & 1.7316 & 2.2726 & 1.6837\tabularnewline
5 & 60922 & 62.5139 & 4.9748 & 4.8295 & 4.0466\tabularnewline
\bottomrule
\end{tabular}
\end{table}
\begin{table}[H]
\centering
\caption{Speedup ratio for the computation time to solve the linear system of equations between the standard Jacobi basis and Lagrange basis and between the standard Jacobi basis and the minimum energy bases.}
\label{Tab:SpeedUp_IMPACT_3D}
\begin{tabular}{@{}ccccc@{}}
\toprule
\multirow{2}{*}{Order} & \multirow{2}{*}{\begin{tabular}[c]{@{}c@{}}Number \\ of DOFs\end{tabular}} & \multicolumn{3}{c}{Speedup} \\ \cmidrule(l){3-5}
                        &                                                                           & Lagrange       & SDME-M            & SDME-H  \\ \cmidrule(r){1-5}
1 & 652 & 0.98 & - & -\tabularnewline
2 & 4318 & 4.29 & 5.35 & 5.37\tabularnewline
3 & 14299 & 3.67 & 4.11 & 4.28\tabularnewline
4 & 31753 & 14.29 & 10.89 & 14.69\tabularnewline
5 & 60922 & 12.57 & 12.94 & 15.85\tabularnewline
\bottomrule
\end{tabular}
\end{table}
%

\section{Conclusions}

In this work, we applied high-order finite element bases to solve two and three-dimensional
transient nonlinear  structural and impact problems. The one-dimensional bases were constructed
by performing simultaneous diagonalization of the internal modes and Schur complement of the boundary modes.
Fabricated smooth solutions involving large displacements and strains were used to test the bases in static
and transient (explicit and implicit) analyses.

The SDME bases performed significantly much better than the standard Jacobi basis for all nonlinear problems
tested. For the static nonlinear test, the SDME-H basis had a better performance than the SDME-M basis for
all polynomial orders when using the Gauss-Seidel preconditioner for the linear system solution. The same
was observed when using the diagonal preconditioner for polynomial orders up to 6.

In the case of transient nonlinear problems with explicit time integration, the SDME-M basis had a speedup up
of 27.7 when compared to the standard Jacobi basis. For the implicit time integration, the best results f
or speedup were achieved by the SDME-M basis as well, with speedup up to 26.

For the conrod simulations, we observed that the SDME-H basis with the chosen set of parameters had a
much better performance compared to the ST basis. Moreover, in the case of implicit time integration, there
were less PCG iterations compared to the SDME-M basis for lower polynomial orders $(P<5)$. However, for
higher polynomial orders, the number of PCG iterations significantly increase for the SDME-H basis, and
the lowest number of iterations was achieved with the SDME-M basis.

The SDME-H basis had the best performance for the impact problems when compared to the standard
Jacobi basis. For the two-dimensional disk impact problem, the SDME-H basis had better performance in the
average by 13\%, when compared to the nodal Lagrange basis with Gauss-Lobatto collocation points. For the
higher polynomial orders ($P=5$ to $P=10$), the improvement was close to 22\%. For the three-dimensional
cylinder impact problem, the improvement for the same comparison and $P=5$ was close to 23\%. When
compared to the nodal Lagrange basis, the performance of the SDME-H was over 25\% better.

{
The Gauss-Seidel preconditioner performed much better than the diagonal preconditioner in terms
of the number of iterations for convergence for a given tolerance. However, it is difficult to
implement it in an element wise fashion and requires the assemble of the global matrix. The diagonal
preconditioner is simpler to implement in parallel without the need of the global matrix of the system of
equations.
}

In general, the SDME bases had an outstanding performance when applied to non-linear structural problems
including large deformation and strain and impact problems. For meshes with larger degrees of freedom, it is
expected a better performance of the SDME bases.

\section*{Acknowledgements}
The authors gratefully acknowledge the support of the National Council for Scientific and Technological Development (CNPq), grant numbers 164733/2017-5 and 310351/2019-7, and the University of Campinas (UNICAMP).

 \bibliographystyle{elsarticle-num}
 \bibliography{test_ref}

\end{document}